\documentclass{amsart}

\usepackage[utf8]{inputenc}
\usepackage{a4}
\usepackage{graphics}
\usepackage{color}
\usepackage{amssymb}
\usepackage[all]{xy}
\usepackage{amsmath}
\usepackage{amsthm}
\usepackage{stmaryrd}
\usepackage{subfigure}
\usepackage{longtable}
\usepackage{setspace}
\usepackage{booktabs}

\usepackage{pst-plot}
\usepackage{pstricks}
\usepackage{color}
\usepackage{nicefrac}

\CompileMatrices
\newtheorem*{introthm1}{Theorem~\ref{thm:factrings}}
\newtheorem*{introthm2}{Theorem~\ref{thm:3fano}}

\newtheorem{theorem}{Theorem}[section]
\newtheorem{lemma}[theorem]{Lemma}
\newtheorem{proposition}[theorem]{Proposition}
\newtheorem{corollary}[theorem]{Corollary}

\theoremstyle{definition}
\newtheorem{definition}[theorem]{Definition}
\newtheorem{example}[theorem]{Example}
\newtheorem{construction}[theorem]{Construction}
\theoremstyle{remark}
\newtheorem{remark}[theorem]{Remark}
\numberwithin{equation}{section}
\def\Chi{{\mathbb X}}

\def\div{{\rm div}}

\def\mal{\! \cdot \!}

\def\rq#1{\widehat{#1}}
\def\t#1{\widetilde{#1}}
\def\b#1{\overline{#1}}
\def\bangle#1{\langle #1 \rangle}

\def\KK{{\mathbb K}}
\def\TT{{\mathbb T}}
\def\ZZ{{\mathbb Z}}

\def\QQ{{\mathbb Q}}
\def\PP{{\mathbb P}}

\def\WDiv{\operatorname{WDiv}}

\def\ord{{\rm ord}}

\def\Aut{\operatorname{Aut}}

\def\Cl{\operatorname{Cl}}

\def\Pic{\operatorname{Pic}}

\def\Hom{{\rm Hom}}

\def\grad{{\rm grad}}

\def\Spec{{\rm Spec}}

\def\cone{{\rm cone}}
\def\lin{{\rm lin}}

\def\Aut{{\rm Aut}}
\def\rk{{\rm rk}\,}

\newcommand{\D}{\mathcal{D}}

\newcommand{\fan}{\Xi}

\newcommand{\myrule}[2]{
\par
\noindent
\begin{minipage}{\textwidth}
#1 \rule[0.3\baselineskip]{.95\textwidth}{1pt}
\par
\vspace{-.1em}
#2 
\par
\vspace{.4em}
\end{minipage}}

\newgray{grayi}{0.68}
\newgray{grayii}{0.60}
\newgray{grayiii}{0.55}
\newgray{grayiv}{0.50}
\newgray{grayv}{0.45}

\def \marks(#1)(#2)(#3){}
\def \definerays (#1,#2)(#3,#4)(#5,#6)(#7,#8){
\renewcommand{\ra}[2]{! #1 100 mul ##1 add #2 100 mul ##2 add}
\renewcommand{\rb}[2]{! #3 100 mul ##1 add #4 100 mul ##2 add}
\renewcommand{\rc}[2]{! #5 100 mul ##1 add #6 100 mul ##2 add}
\renewcommand{\rd}[2]{! #7 100 mul ##1 add #8 100 mul ##2 add}
}

\def \defineextrarays (#1,#2)(#3,#4)(#5,#6)(#7,#8){
\renewcommand{\re}[2]{! #1 20 mul ##1 add #2 20 mul  ##2 add}
\renewcommand{\rf}[2]{! #3 20 mul ##1 add #4 20 mul  ##2 add}
\renewcommand{\rg}[2]{! #5 20 mul ##1 add #6 20 mul  ##2 add}
\renewcommand{\rh}[2]{! #7 20 mul ##1 add #8 20 mul  ##2 add}
}

\def \oput (#1,#2){\uput{3.2}[(#1,#2)](0,0){\tiny (#1,#2)}}
\def \oputm(#1,#2){\uput{3}[(#1,#2)](0,0){\tiny (#1,#2)}}
\def \oputmm(#1,#2){\uput{2.6}[(#1,#2)](0,0){\tiny (#1,#2)}}
\def \cput [#1](#2,#3){\uput[#1](#2,#3){(#2,#3)} \psdot[dotstyle=o](#2,#3)}
\newcommand{\rg}{}
\newcommand{\rh}{}
\newcommand{\ra}[2]{! 10 #1 add 10 #2 add}
\def \raya(#1,#2){\psline(#1,#2)(\ra{#1}{#2})}
\newcommand{\rb}[2]{! 20 #1 add -10 #2 add}
\def \rayb(#1,#2){\psline(#1,#2)(\rb{#1}{#2})}
\newcommand{\rc}[2]{! -10 #1 add -10 #2 add}
\def \rayc(#1,#2){\psline(#1,#2)(\rc{#1}{#2})}
\newcommand{\rd}[2]{! -10 #1 add 20 #2 add}
\def \rayd(#1,#2){\psline(#1,#2)(\rd{#1}{#2})}
\newcommand{\re}[2]{! 0 #1 add -10 #2 add}
\def \raye(#1,#2){\psline(#1,#2)(\re{#1}{#2})}
\newcommand{\rf}[2]{! -10 #1 add 0 #2 add}
\def \rayf(#1,#2){\psline(#1,#2)(\rf{#1}{#2})}
\def \rayg(#1,#2){\psline(#1,#2)(\rg{#1}{#2})}
\def \rayh(#1,#2){\psline(#1,#2)(\rh{#1}{#2})}

\def \beginpsclip{\psclip{\psccurve[linestyle=none](2,2)(2.8,0)(2.5,-2)(-3.2,0)(-2.2,2.5)(0,2.8)(2,2)}}

\def \stara (#1,#2){
\psset{fillstyle=solid, fillcolor=grayii} \cab(#1,#2) \psset{fillcolor=grayiii} \cbc(#1,#2) \psset{fillcolor=grayiv} \cca(#1,#2)
\psdot[dotstyle=*](0,0)
}

\def \starb (#1,#2){
\psset{fillstyle=solid, fillcolor=grayii} \cab(#1,#2) \psset{fillcolor=grayiii} \cbc(#1,#2) \psset{fillcolor=grayiv} \ccd(#1,#2) \psset{fillcolor=grayi} \cda(#1,#2)
\psdot[dotstyle=*](0,0)
}

\def \subdiva (#1,#2)(#3,#4)(#5,#6){
\psset{fillstyle=solid, fillcolor=grayii} \lab(#1,#2)(#3,#4) \psset{fillcolor=grayiii} \lbc(#3,#4)(#5,#6) \psset{fillcolor=grayiv}  \lca(#5,#6)(#1,#2) \psset{fillcolor=grayi} \psline(#1,#2)(#3,#4)(#5,#6)
\psdot[dotstyle=*](0,0)
}

\def \subdivh (#1,#2)(#3,#4){
\psset{fillstyle=solid, fillcolor=grayi} \lda (#1,#2)(#3,#4) \psset{fillcolor=grayii} \cab(#3,#4) \psset{fillcolor=grayiii}  \lbc(#3,#4)(#1,#2) \psset{fillcolor=grayiv} \lcd(#1,#2)(#1,#2)
\psdot[dotstyle=*](0,0)
}

\def \subdivv (#1,#2)(#3,#4){
\psset{fillstyle=solid, fillcolor=grayi} \cda(#1,#2) \psset{fillcolor=grayii} \lab(#1,#2)(#3,#4)  \psset{fillcolor=grayiii} \cbc(#3,#4)  \psset{fillcolor=grayiv} \lcd(#3,#4)(#1,#2)
\psdot[dotstyle=*](0,0)

}

\def \cab (#1,#2){\psline(\ra{#1}{#2})(#1,#2)(\rb{#1}{#2})}
\def \cbc (#1,#2){\psline(\rb{#1}{#2})(#1,#2)(\rc{#1}{#2})}
\def \ccd (#1,#2){\psline(\rc{#1}{#2})(#1,#2)(\rd{#1}{#2})}
\def \cda (#1,#2){\psline(\rd{#1}{#2})(#1,#2)(\ra{#1}{#2})}
\def \cef (#1,#2){\psline(\re{#1}{#2})(#1,#2)(\rf{#1}{#2})}
\def \cde (#1,#2){\psline(\rd{#1}{#2})(#1,#2)(\re{#1}{#2})}
\def \cca (#1,#2){\psline(\rc{#1}{#2})(#1,#2)(\ra{#1}{#2})}
\def \cce (#1,#2){\psline(\rc{#1}{#2})(#1,#2)(\re{#1}{#2})}
\def \cgh (#1,#2){\psline(\rg{#1}{#2})(#1,#2)(\rh{#1}{#2})}

\def \lab (#1,#2)(#3,#4){\psline(\ra{#1}{#2})(#1,#2)(#3,#4)(\rb{#3}{#4})}
\def \lcd (#1,#2)(#3,#4){\psline(\rc{#1}{#2})(#1,#2)(#3,#4)(\rd{#3}{#4})}
\def \lef (#1,#2)(#3,#4){\psline(\re{#1}{#2})(#1,#2)(#3,#4)(\rf{#3}{#4})}
\def \lda (#1,#2)(#3,#4){\psline(\rd{#1}{#2})(#1,#2)(#3,#4)(\ra{#3}{#4})}
\def \lfa (#1,#2)(#3,#4){\psline(\rf{#1}{#2})(#1,#2)(#3,#4)(\ra{#3}{#4})}
\def \lde (#1,#2)(#3,#4){\psline(\rd{#1}{#2})(#1,#2)(#3,#4)(\re{#3}{#4})}
\def \lbc (#1,#2)(#3,#4){\psline(\rb{#1}{#2})(#1,#2)(#3,#4)(\rc{#3}{#4})}
\def \lca (#1,#2)(#3,#4){\psline(\rc{#1}{#2})(#1,#2)(#3,#4)(\ra{#3}{#4})}

\SpecialCoor
\psset{unit=0.43cm}

\def \threefoldA {
\definerays(1,1)(1,-1)(-1,-1)(-1,1)
\defineextrarays(0,0)(0,0)(0,0)(0,0)
\begin{pspicture}(-3,-2.6)(3,3.5)%
\beginpsclip 
 \psset{linewidth=0.5pt,fillstyle=solid}%
 \subdivh(-1,0)(0,0)
 \psset{linestyle=dashed,fillstyle=none}
 \cef(0,0) 
 \uput[150](-1,0){$\scriptstyle (-1,0)$}
\endpsclip
\end{pspicture}
\begin{pspicture}(-3.2,-2.6)(3,3.5)%
\beginpsclip 
 \psset{linewidth=0.5pt,fillstyle=solid}%
 \subdivv(0,0)(0,-1)
 \psset{linestyle=dashed,fillstyle=none}
 \cef(0,0) 
 \uput[190](0,-1){$\scriptstyle (0,-1)$}
\endpsclip
\end{pspicture}
\begin{pspicture}(-3.2,-2.6)(3,3.5)%
\beginpsclip 
 \psset{linewidth=0.5pt,fillstyle=solid}%
 \starb(0.5,0.5)
 \psset{linestyle=dashed,fillstyle=none}
 \cef(0,0)
 \uput[0](0.5,0.5){$\scriptstyle (\frac{1}{2},\frac{1}{2})$}
\endpsclip
\end{pspicture}
\begin{pspicture}(-4.2,-2.6)(3.9,3.5)%
\beginpsclip 
 \psset{linewidth=0.6pt,fillstyle=solid}%
 \starb(0,0)
 \psset{linestyle=dashed,fillstyle=none}
 \cef(0,0)
 \psset{linewidth=1pt,linecolor=white}
 \psdot(0,0)
\endpsclip
 \oput(-1,1)
 \oput(1,1)
 \oput(1,-1)
 \oput(-1,-1)
\end{pspicture}}

\def \threefoldB{
\definerays(1,1)(2,-1)(-1,-1)(-1,2)
\defineextrarays(0,-1)(-1,0)(0,0)(0,0)
\begin{pspicture}(-3,-2.6)(3,3.5)%
\beginpsclip 
 \psset{linewidth=0.5pt,fillstyle=solid}%
 \subdivh(-0.5,0)(0,0)
\psset{linestyle=dashed,fillstyle=none}
 \cef(0,0) 
 \uput[150](-0.5,0){$\scriptstyle (-\frac{1}{2},0)$}
\endpsclip
\end{pspicture}
\begin{pspicture}(-3.2,-2.6)(3,3.5)%
\beginpsclip 
 \psset{linewidth=0.5pt,fillstyle=solid}%
 \subdivv(0,0)(0,-0.5)
 \psset{linestyle=dashed,fillstyle=none}
 \cef(0,0) 
 \uput[190](0,-0.5){$\scriptstyle (0,-\frac{1}{2})$}
\endpsclip
\end{pspicture}
\begin{pspicture}(-3.2,-2.6)(3,3.5)%
\beginpsclip 
 \psset{linewidth=0.5pt,fillstyle=solid}%
 \starb(0.33,0.33)
 \psset{linestyle=dashed,fillstyle=none}
 \cef(0,0)
 \uput{0.5}[0](0.33,0.33){$\scriptstyle (\frac{1}{3},\frac{1}{3})$}
\endpsclip
\end{pspicture}
\begin{pspicture}(-4.2,-2.6)(3.9,3.5)%
\beginpsclip 
 \psset{linewidth=0.5pt,fillstyle=solid}%
 \starb(0,0)
 \psset{linestyle=dashed,fillstyle=none}
 \cef(0,0)
 \psset{linewidth=1pt,linecolor=white}
 \psdot(0,0)
\endpsclip
 \oput(-1,2)
 \oput(1,1)
 \oput(2,-1)
 \oput(-1,-1)
 \oputmm(0,-1)
 \oput(-1,0)
\end{pspicture}}

\def \threefoldC{
\definerays(1,1)(3,-1)(-1,-1)(-1,3)
\defineextrarays(0,-1)(-1,0)(0,0)(0,0)

\begin{pspicture}(-3,-2.6)(3,3.5)%
\beginpsclip 
 \psset{linewidth=0.5pt,fillstyle=solid}%
 \subdivh(-0.5,0)(0,0)
\psset{linestyle=dashed,fillstyle=none}
 \cef(0,0) 
 \uput[150](-0.5,0){$\scriptstyle (-\frac{1}{3},0)$}
\endpsclip
\end{pspicture}
\begin{pspicture}(-3.2,-2.6)(3,3.5)%
\beginpsclip 
 \psset{linewidth=0.5pt,fillstyle=solid}%
 \subdivv(0,0)(0,-0.5)
 \psset{linestyle=dashed,fillstyle=none}
 \cef(0,0) 
 \uput[190](0,-0.5){$\scriptstyle (0,-\frac{1}{3})$}
\endpsclip
\end{pspicture}
\begin{pspicture}(-3.2,-2.6)(3,3.5)%
\beginpsclip 
 \psset{linewidth=0.5pt,fillstyle=solid}%
 \starb(0.33,0.33)
 \psset{linestyle=dashed,fillstyle=none}
 \cef(0,0)
 \uput[10](0.33,0.33){$\scriptstyle (\frac{1}{4},\frac{1}{4})$}
\endpsclip
\end{pspicture}
\begin{pspicture}(-4.2,-2.6)(3.9,3.5)%
\beginpsclip 
 \psset{linewidth=0.6pt,fillstyle=solid}%
 \starb(0,0)
 \psset{linestyle=dashed,fillstyle=none}
 \cef(0,0)
 \psset{linewidth=1pt,linecolor=white}
 \psdot(0,0)
\endpsclip
 \oput(-1,3)
 \oput(1,1)
 \oput(3,-1)
 \oput(-1,-1)
 \oputmm(0,-1)
 \oput(-1,0)
\end{pspicture}}

\def \threefoldD{
\definerays(1,2)(1,0)(-5,-2)(-1,0)
\defineextrarays(1,0)(1,0)(3,1)(0,-1)

\begin{pspicture}(-3,-2.6)(3,3.5)%
\beginpsclip 
 \psset{linewidth=0.5pt,fillstyle=solid}%
 \subdivh(0,0)(0.6,0.2)
\psset{linestyle=dashed,fillstyle=none}
 \rayb(0,0) 
 \rayh(0,0) 
 \rayg(0.6,0.2) 
 \uput[130](0.6,0.2){$\scriptstyle (\frac{3}{5},\frac{1}{5})$}
\endpsclip
\end{pspicture}
\begin{pspicture}(-3.2,-2.6)(3,3.5)%
\beginpsclip 
 \psset{linewidth=0.5pt,fillstyle=solid}%
 \subdivv(0,0)(0,-0.2)
 \psset{linestyle=dashed,fillstyle=none}
 \rayb(0,0) 
 \rayg(0,0) 
 \rayh(0,-0.2) 
 \uput[-50](0,-0.2){$\scriptstyle (0,-\frac{1}{5})$}
\endpsclip
\end{pspicture}
\begin{pspicture}(-3.2,-2.6)(3,3.5)%
\beginpsclip 
 \psset{linewidth=0.5pt,fillstyle=solid}%
 \starb(-0.5,0)
 \psset{linestyle=dashed,fillstyle=none}
 \rayg(0,0) 
 \rayh(0,0) 
 \uput[135](-0.5,0){$\scriptstyle (-\frac{1}{2},0)$}
\endpsclip
\end{pspicture}
\begin{pspicture}(-4.2,-2.6)(3.9,3.5)%
\beginpsclip 
 \psset{linewidth=0.6pt,fillstyle=solid}%
 \starb(0,0)
 \psset{linestyle=dashed,fillstyle=none}
 \rayg(0,0) 
 \rayh(0,0) 
 \psset{linewidth=1pt,linestyle=dashed,fillstyle=none,linecolor=white}
 \rayb(0,-0.06) 
 \psdot(0,0)
\endpsclip
\oput(1,2)\oput(1,0)\oput(-5,-2)\oput(-1,0)\oput(3,1)\oputmm(0,-1)
\end{pspicture}}

\def \threefoldDplain{
\definerays(1,2)(1,0)(-5,-2)(-1,0)
\defineextrarays(1,0)(1,0)(3,1)(0,-1)

\begin{pspicture}(-3,-2.6)(3,3.5)%
\beginpsclip 
 \psset{linewidth=0.5pt,fillstyle=solid}%
 \subdivh(0,0)(0.6,0.2)
 \uput[130](0.6,0.2){$\scriptstyle (\frac{3}{5},\frac{1}{5})$}
\endpsclip
\rput(-2.7,-1.8){$\fan_0$}
\end{pspicture}
\begin{pspicture}(-3.2,-2.6)(3,3.5)%
\beginpsclip 
 \psset{linewidth=0.5pt,fillstyle=solid}%
 \subdivv(0,0)(0,-0.2)
 \uput[-50](0,-0.2){$\scriptstyle (0,-\frac{1}{5})$}
\endpsclip
\rput(-2.7,-1.8){$\fan_1$}
\end{pspicture}
\begin{pspicture}(-3.2,-2.6)(3,3.5)%
\beginpsclip 
 \psset{linewidth=0.5pt,fillstyle=solid}%
 \starb(-0.5,0)
 \uput[135](-0.5,0){$\scriptstyle (-\frac{1 }{2},0)$}
\endpsclip
\rput(-2.7,-1.8){$\fan_{\infty}$}
\end{pspicture}
\begin{pspicture}(-4.2,-2.6)(3.9,3.5)%
\beginpsclip 
 \psset{linewidth=0.6pt,fillstyle=solid}%
 \starb(0,0)
\psset{linewidth=1pt,linestyle=dashed,fillstyle=none,linecolor=white}
  \psdot(0,0)
\endpsclip
\rput(-2.7,-1.8){$\fan_Y$}
\oput(1,2)\oput(1,0)\oput(-5,-2)\oput(-1,0)
\end{pspicture}}
\def \threefoldE {
\definerays(-4,3)(1,0)(2,-2.6)(1,0)
\defineextrarays(2,-1)(1,-1)(0,1)(-1,1)

\begin{pspicture}(-3,-2.6)(3,3.5)%
\beginpsclip 
 \psset{unit=1cm}
 \psset{linewidth=0.5pt,fillstyle=solid}%
 \subdiva(-1,1)(-0.25,0.5)(0,0)
 \psset{linestyle=dashed,fillstyle=none} 
 \cef(0,0)
 \cgh(-1,1)
 \rayd(0,0)
 \rayd(-1,1)
 \psline(0,0.5)(0,0)
 \psline(0,0.5)(-1,1)

 \uput[230](-0.25,0.5){$\scriptstyle (-\frac{1}{4},\frac{1}{2})$}
\uput[30](0,0.5){$\scriptstyle (0,\frac{1}{2})$}
\endpsclip
\psset{unit=1cm}
\uput[30](-1,1){$\scriptscriptstyle (-1,1)$}
\end{pspicture}
\begin{pspicture}(-3.2,-2.6)(3,3.5)%
\beginpsclip 
 \psset{unit=1cm}
 \psset{linewidth=0.5pt,fillstyle=solid}%
 \stara(-0.66,0)
 \psset{linestyle=dashed,fillstyle=none}
  \psset{linestyle=dashed,fillstyle=none} 
 \rayf(-.5,0)
 \raye(0,0)
 \rayh(-.5,0)
 \rayg(0,0)

\endpsclip
 \psset{unit=1cm}
 \uput{0}[220](-0.66,0){$\scriptstyle (-\frac{2}{3},0)$}
 \uput{0.5}[90](-0.5,0){$\scriptstyle (-\frac{1}{2},0)$}
 \psset{unit=.5cm}
\end{pspicture}
\begin{pspicture}(-3.2,-2.6)(3,3.5)%
\beginpsclip 
 \psset{unit=1cm}
 \psset{linewidth=0.5pt,fillstyle=solid}%
 \stara(1,-.5)
 \psset{linestyle=dashed,fillstyle=none}
 \uput[210](1,-0.5){$\scriptstyle (1,-\frac{1}{2})$}
 \cef(1,-0.5)
 \cgh(1,-0.5)
 \psset{unit=.5cm}
\endpsclip
\end{pspicture}
\begin{pspicture}(-4.2,-2.6)(3.9,3.5)%
\beginpsclip 
 \psset{linewidth=0.6pt,fillstyle=solid}%
 \stara(0,0)
 \psset{linestyle=dashed,fillstyle=none}
 \cef(0,0)
 \cgh(0,0)
 \psset{linewidth=1pt,linestyle=dashed,fillstyle=none,linecolor=white}
 \rayb(0,-0.06) 
 \psdot(0,0)
 \endpsclip
 \oput(-4,3)
 \oput(1,0)
 \oput(2,-2.6)
 \oput(2,-1)
 \oput(1,-1)
 \oputm(0,1)
 \oput(-1,1)
\end{pspicture}}

\def \threefoldF{
\definerays(0,1)(1,0)(-3,-2)(1,0)
\defineextrarays(-1,0)(1,1)(0,-1)(0,0)

\begin{pspicture}(-3,-2.6)(3,3.5)%
\beginpsclip 
 \psset{unit=1cm}
 \psset{linewidth=0.5pt,fillstyle=solid}%
 \subdiva(0.5,0.5)(0.66,0.33)(0,0)
 \psset{linestyle=dashed,fillstyle=none} 
  \rayb(0.5,0.5)
  \raye(0,0)
  \rayf(0.5,0.5)
  \rayb(0,0)
  \rayg(0,0)
  \uput[135](0.5,0.5){$\scriptstyle (\frac{1}{2},\frac{1}{2})$}
  \uput{0.4}[285](0.66,0.33){$\scriptstyle (\frac{2}{3},\frac{1}{3})$}
\psset{unit=.5cm}
\endpsclip
\end{pspicture}
\begin{pspicture}(-3.2,-2.6)(3,3.5)%
\beginpsclip 
 \psset{unit=1cm}
 \psset{linewidth=0.5pt,fillstyle=solid}%
 \stara(0,-0.33)
 \psset{linestyle=dashed,fillstyle=none}
\rayf(0,0)
\raye(0,0)
\rayd(0,0)
\rayg(0,-0.33)
\uput[310](0,-0.33){$\scriptstyle (0,-\frac{1}{3})$}
 \psset{unit=.5cm}
\endpsclip
\end{pspicture}
\begin{pspicture}(-3.2,-2.6)(3,3.5)%
\beginpsclip 
 \psset{unit=1cm}
 \psset{linewidth=0.5pt,fillstyle=solid}%
 \stara(-0.5,0)
 \psset{linestyle=dashed,fillstyle=none}
 \raye(-0.5,0)
 \rayf(0,0)
 \rayg(0,0)
\uput{0.1}[140](-0.5,0){$\scriptstyle (-\frac{1}{2},0)$}
 \psset{unit=.5cm}
\endpsclip
\end{pspicture}
\begin{pspicture}(-4.2,-2.6)(3.9,3.5)%
\beginpsclip 
 \psset{linewidth=0.6pt,fillstyle=solid}%
 \stara(0,0)
 \psset{linestyle=dashed,fillstyle=none}
 \cef(0,0)
 \rayg(0,0) 
 \psset{linewidth=1pt,linestyle=dashed,fillstyle=none,linecolor=white}
 \rayb(0,-0.06) 
 \psdot(0,0)
 \endpsclip
 \oput(-3,-2)
 \oput(1,0)
 \oputm(0,1)
 \oput(-1,0)
 \oput(1,1)
 \oputmm(0,-1)
\end{pspicture}}

\def \threefoldG{
\definerays(0,-1)(-1,6)(1,0)(0,-1)
\defineextrarays(0,1)(0,0)(0,0)(0,0)

\begin{pspicture}(-3,-2.6)(3,3.5)%
\beginpsclip 
 \psset{xunit=1.2cm}
 \psset{linewidth=0.5pt,fillstyle=solid}%
 \subdivh(-0.8,0)(-1,1)
 \psset{linestyle=dashed,fillstyle=none} 
  \raya(-0.75,0)
  \raya(-0.66,0)
  \raya(-0.5,0)
  \raya(0,0)
  \raye(-0.75,0)
  \raye(-0.66,0)
  \raye(-0.5,0)
  \raye(0,0)  
  \raye(-0.8,0)  
  \raye(-1,1)  
\psset{unit=.5cm}
\endpsclip
 \psset{xunit=1.2cm}
 \uput{0.2}[150](-1,1){$\scriptstyle (-1,1)$}
 \psline[linestyle=dotted](-1.1,-2.7)(-0.8,-1.5)
 \rput(-1.1,-2.6){$\scriptstyle -\frac{4}{5}$}
 \psline[linestyle=dotted](-0.75,-2.7)(-0.75,-1.5)
 \rput(-0.75,-2.6){$\scriptstyle -\frac{3}{4}$}
 \psline[linestyle=dotted](-0.5,-2.7)(-0.66,-1.8)
 \rput(-0.5,-2.6){$\scriptstyle -\frac{2}{3}$}
 \psline[linestyle=dotted](-0.2,-2.7)(-0.5,-2)
 \rput(-0.2,-2.6){$\scriptstyle -\frac{1}{2}$}
\end{pspicture}
\begin{pspicture}(-3.2,-2.6)(3,3.5)%
\beginpsclip 
 \psset{xunit=0.9cm}
 \psset{linewidth=0.5pt,fillstyle=solid}%
 \stara(0.33,0)
 \psset{linestyle=dashed,fillstyle=none}
 \raya(0.5,0)
 \raya(1,0)
 \raye(0.5,0)
 \raye(1,0)
 \raye(0.33,0)  
\endpsclip
 \psset{xunit=0.9cm}
 \psline[linestyle=dotted](0.33,-2.7)(0.33,-2.0)
 \rput(0.33,-2.6){$\scriptstyle \frac{1}{3}$}
 \psline[linestyle=dotted](0.5,-2.7)(0.5,-2.0)
 \rput(0.5,-2.6){$\scriptstyle \frac{1}{2}$}
 \psline[linestyle=dotted](1,-2.7)(1,-2.0)
 \rput(1,-2.6){$\scriptstyle 1$}
\end{pspicture}
\begin{pspicture}(-3.2,-2.6)(3,3.5)%
\beginpsclip 
 \psset{xunit=0.9cm}
 \psset{linewidth=0.5pt,fillstyle=solid}%
 \stara(0.5,0)
 \psset{linestyle=dashed,fillstyle=none}
 \raya(1,0)
 \raye(1,0)
 \raye(0.5,0)  
\endpsclip
 \psset{xunit=0.9cm}
 \psline[linestyle=dotted](0.5,-2.7)(0.5,-2.0)
 \rput(0.5,-2.6){$\scriptstyle \frac{1}{2}$}
 \psline[linestyle=dotted](1,-2.7)(1,-2.0)
 \rput(1,-2.6){$\scriptstyle 1$}
\end{pspicture}
\begin{pspicture}(-4.2,-2.6)(3.9,3.5)%
\beginpsclip 
 \psset{xunit=0.9cm}
 \psset{linewidth=0.6pt,fillstyle=solid}%
 \stara(0,0)
 \psset{linestyle=dashed,fillstyle=none}
 \cef(0,0)
 \psset{linewidth=1pt,linestyle=solid,fillstyle=none,linecolor=white}
 \rayd(0,0) 
 \psset{linestyle=dashed}
 \rayc(0,-0.06) 
 \raye(0,0) 
 \psset{linestyle=none,fillstyle=vlines,hatchcolor=white,hatchwidth=0.02}
 \ccd(0.1,-0.2)
 \cce(0.1,0.3)
 \psdot(0,0)
 \endpsclip
 \psset{xunit=0.9cm}
 \oputm(0,1)\oputmm(0,-1)\oput(-1,6)\oput(1,0)
\end{pspicture}}

\def \threefoldH{
\definerays(1,2)(1,0)(-3,-2)(-1,0)
\defineextrarays(1,0)(1,0)(2,1)(0,-1)

\begin{pspicture}(-3,-2.6)(3,3.5)%
\beginpsclip 
 \psset{linewidth=0.5pt,fillstyle=solid}%
 \subdivh(0,0)(0.66,0.33)
\psset{linestyle=dashed,fillstyle=none}
 \rayb(0,0) 
 \rayh(0,0) 
 \rayg(0.66,0.33) 
 \uput[130](0.66,0.33){$\scriptstyle (\frac{2}{3},\frac{1}{3})$}
\endpsclip
\end{pspicture}
\begin{pspicture}(-3.2,-2.6)(3,3.5)%
\beginpsclip 
 \psset{linewidth=0.5pt,fillstyle=solid}%
 \subdivv(0,0)(0,-0.33)
 \psset{linestyle=dashed,fillstyle=none}
 \rayb(0,0) 
 \rayg(0,0) 
 \rayh(0,-0.33) 
 \uput[-50](0,-0.33){$\scriptstyle (0,-\frac{1}{3})$}
\endpsclip
\end{pspicture}
\begin{pspicture}(-3.2,-2.6)(3,3.5)%
\beginpsclip 
 \psset{linewidth=0.5pt,fillstyle=solid}%
 \starb(-0.5,0)
 \psset{linestyle=dashed,fillstyle=none}
 \rayg(0,0) 
 \rayh(0,0) 
 \uput[135](-0.5,0){$\scriptstyle (-\frac{1}{2},0)$}
\endpsclip
\end{pspicture}
\begin{pspicture}(-4.2,-2.6)(3.9,3.5)%
\beginpsclip 
 \psset{linewidth=0.6pt,fillstyle=solid}%
 \starb(0,0)
 \psset{linestyle=dashed,fillstyle=none}
 \rayg(0,0) 
 \rayh(0,0) 
 \psset{linewidth=1pt,linestyle=dashed,fillstyle=none,linecolor=white}
 \rayb(0,-0.06) 
 \psdot(0,0)
\endpsclip
\oput(1,2)\oput(1,0)\oput(-3,-2)\oput(-1,0)\oput(2,1)\oputmm(0,-1)
\end{pspicture}}

\newcounter{itemnumber}

\begin{document}
\title[Multigraded Factorial Rings and Fano varieties]%
{Multigraded Factorial Rings \\
and Fano varieties with torus action}
\author[J.~Hausen]{J\"urgen Hausen} 
\address{Mathematisches Institut, Universit\"at T\"ubingen,
Auf der Morgenstelle 10, 72076 T\"ubingen, Germany}
\email{juergen.hausen@uni-tuebingen.de}
\author[E.~Herppich]{Elaine Herppich} 
\address{Mathematisches Institut, Universit\"at T\"ubingen,
Auf der Morgenstelle 10, 72076 T\"ubingen, Germany}
\email{elaine.herppich@uni-tuebingen.de}
\author[H.~S\"uss]{Hendrik S\"uss}
\address{Institut f\"ur Mathematik,
        LS Algebra und Geometrie,
        Brandenburgische Technische Universit\"at Cottbus,
        PF 10 13 44, 
        03013 Cottbus, Germany}
\email{suess@math.tu-cottbus.de}

\begin{abstract}
In a first result, 
we describe all finitely generated 
factorial algebras over an algebraically 
closed field of characteristic zero that
come with an effective multigrading of 
complexity one by means of generators 
and relations.
This enables us to construct 
systematically varieties 
with free divisor class group
and a complexity one torus action
via their Cox rings.
For the Fano varieties of this type 
that have a free 
divisor class group of rank one, 
we provide explicit bounds for the 
number of possible 
deformation types depending on the 
dimension and the index of the Picard group 
in the divisor class group.
As a consequence, one can produce classification 
lists for fixed dimension and Picard index.
We carry this out expemplarily in the following 
cases.
There are 15 non-toric
surfaces with Picard index at most six.
Moreover, there are 116 non-toric threefolds 
with Picard index at most two;
nine of them are locally factorial, i.e. 
of Picard index one, and among these
one is smooth, six have canonical singularities
and two have non-canonical singularities.
Finally, there are 67 non-toric 
locally factorial fourfolds and two 
one-dimensional families of non-toric locally 
factorial fourfolds.
In all cases, we list the Cox rings explicitly.
\end{abstract}

\maketitle
\maketitle

\section*{Introduction}

Let $\KK$ be an algebraically closed field
of characteristic zero.
A first aim of this paper is to 
determine all finitely generated factorial 
$\KK$-algebras $R$ with an effective complexity 
one multigrading 
$R = \oplus_{u \in M} R_u$ satisfying $R_0 = \KK$; 
here effective complexity one multigrading 
means that with $d := \dim \, R$ we have 
$M \cong \ZZ^{d-1}$ and 
the $u \in M$ with $R_u \ne 0$ generate $M$ 
as a $\ZZ$-module.
Our result extends work by 
Mori~\cite{Mo} and Ishida~\cite{Is},
who settled the cases $d=2$ and $d=3$.

An obvious class of multigraded factorial
algebras as above is given by polynomial rings. 
A much larger class is obtained as follows.
Take a sequence $A = (a_0, \ldots, a_r)$ 
of vectors $a_i \in \KK^2$
such that $(a_i,a_k)$ is linearly independent
whenever $k \ne i$,
a sequence $\mathfrak{n} = (n_0, \ldots, n_r)$ 
of positive integers
and a family $L = (l_{ij})$ 
of positive integers,
where $0 \le i \le r$ and 
$1 \le j \le n_i$.
For every $0 \le i \le r$, we
define a monomial
$$
f_i 
\  := \
T_{i1}^{l_{i1}} \cdots T_{in_i}^{l_{in_i}}
\ \in \
\KK[T_{ij}; \; 0 \le i \le r, \; 1 \le j \le n_i],
$$
for any two indices $0 \le i,j \le r$,
we set $\alpha_{ij} :=  \det(a_i,a_j)$,
and for any three indices $0 \le i < j < k \le r$,
we define a trinomial
$$
g_{i,j,k} 
\ := \ 
\alpha_{jk}f_i
\ + \ 
\alpha_{ki}f_j
\ + \ 
\alpha_{ij}f_k
\ \in \
\KK[T_{ij}; \; 0 \le i \le r, \; 1 \le j \le n_i].
$$ 
Note that the coefficients of $g_{i,j,k}$ are 
all nonzero. 
The triple $(A,\mathfrak{n},L)$ then defines a 
$\KK$-algebra
$$
R(A,\mathfrak{n},L)
\ := \ 
\KK[T_{ij}; \; 0 \le i \le r, \; 1 \le j \le n_i] 
\ / \
\bangle{g_{i,i+1,i+2}; \; 0 \le i \le r-2}.
$$
It turns out that $R(A,\mathfrak{n},L)$ is a normal 
complete intersection, 
see Proposition~\ref{prop:RAnLnormal}.
In particular, it is of dimension 
\begin{eqnarray*}
\dim \, R(A,\mathfrak{n},L)
& = &
 n_0 + \ldots + n_r \ - \ r \ + \ 1.
\end{eqnarray*}
If the triple $(A,\mathfrak{n},L)$ 
is {\em admissible\/}, i.e., 
the numbers $\gcd(l_{i1}, \ldots, l_{in_i})$,
where $0 \le i \le r$, are pairwise
coprime, then $R(A,\mathfrak{n},L)$ 
admits a canonical effective complexity 
one grading  by a lattice $K$,
see Construction~\ref{constr:Kgrading}.
Our first result is the following.

\goodbreak

\begin{introthm1}
Up to isomorphy, the finitely generated 
factorial $\KK$-algebras with an effective 
complexity one grading $R = \oplus_M R_u$ 
and $R_0 = \KK$ are 
\begin{enumerate}
\item 
the polynomial algebras $\KK[T_1, \ldots, T_d]$ 
with a grading $\deg(T_i) = u_i \in \ZZ^{d-1}$ 
such that $u_1, \ldots, u_d$ generate $\ZZ^{d-1}$ 
as a lattice and the convex cone on $\QQ^{d-1}$
generated by $u_1, \ldots, u_d$ is pointed,
\item
the $(K \times \ZZ^m)$-graded 
algebras $R(A,\mathfrak{n},L)[S_1,\ldots,S_m]$,
where $R(A,\mathfrak{n},L)$ is the $K$-graded
algebra defined by an admissible triple
$(A,\mathfrak{n},L)$ and $\deg \, S_j \in \ZZ^m$ 
is the $j$-th canonical base vector.
\end{enumerate}
\end{introthm1}

The further paper is devoted to
normal (possibly singular) $d$-dimensional 
Fano varieties $X$ with an effective action 
of an algebraic torus $T$.
In the case $\dim \, T = d$, we have the 
meanwhile extensively studied class of
toric Fano varieties, see~\cite{Bat1}, 
\cite{WaWa} and~\cite{Bat2} for the
initiating work.
Our aim is to show that the above 
Theorem provides an approach 
to classification results for the 
case $\dim \, T = d-1$, 
that means Fano varieties with a 
complexity one torus action.
Here, we treat the case of divisor class group
$\Cl(X) \cong \ZZ$; note that in the toric setting
this gives precisely the weighted projective 
spaces. The idea is to consider the Cox ring 
\begin{eqnarray*}
\mathcal{R}(X)
& = & 
\bigoplus_{D \in \Cl(X)} \Gamma(X, \mathcal{O}_X(D)).
\end{eqnarray*}
The ring $\mathcal{R}(X)$ is factorial,
finitely generated as a $\KK$-algebra and 
the $T$-action on $X$ gives rise to an effective 
complexity one multigrading of $\mathcal{R}(X)$
refining the $\Cl(X)$-grading, 
see~\cite{BeHa1} and~\cite{HaSu}.
Consequently, $\mathcal{R}(X)$ is one of the rings 
listed in the first Theorem.
Moreover, $X$ can be easily reconstructed from 
$\mathcal{R}(X)$; it is the homogeneous 
spectrum with respect to the $\Cl(X)$-grading
of $\mathcal{R}(X)$.
Thus, in order to construct Fano varieties,
we firstly have to figure out the Cox rings
among the rings occuring in the first Theorem
and then find those, which belong to a Fano variety;
this is done in Propositions~\ref{prop:coxchar} 
and~\ref{Prop:FanoPicard}.

In order to produce classification results
via this approach,
we need explicit bounds on the number 
of deformation types of Fano varieties with 
prescribed discrete invariants. 
Besides the dimension, in our setting,
a suitable invariant is the 
{\em Picard index\/} $[\Cl(X):\Pic(X)]$.
Denoting by $\xi(\mu)$ the number of 
primes less or equal to $\mu$,
we obtain the following bound,
see Corollary~\ref{cor:finitefanos}:
for any pair $(d,\mu) \in \ZZ^2_{>0}$,
the number $\delta(d,\mu)$ of different 
deformation types of $d$-dimensional 
Fano varieties with a complexity one 
torus action such that
$\Cl(X) \cong \ZZ$ and $\mu = [\Cl(X):\Pic(X)]$
hold is bounded by
\begin{eqnarray*} 
\delta(d,\mu)
& \le &
(6d\mu)^{2\xi(3d\mu)+d-2}\mu^{\xi(\mu)^2 + 2\xi((d+2)\mu)+2d+2}.
\end{eqnarray*}
In particular, we conclude that for fixed 
$\mu \in \ZZ_{>0}$, the number $\delta(d)$ 
of different deformation types of $d$-dimensional 
Fano varieties with a complexity one 
torus action $\Cl(X) \cong \ZZ$ 
and  Picard index $\mu$
is asymptotically bounded by 
$d^{Ad}$ with a constant~$A$ 
depending only on~$\mu$, 
see~Corollary~\ref{cor:asymptoticsmu}.

In fact, in Theorem~\ref{Th:FiniteIndex} we 
even obtain explicit bounds for the discrete input 
data of the rings $R(A,\mathfrak{n},L)[S_1,\ldots,S_m]$.
This allows us to construct all
Fano varieties $X$ with prescribed dimension
and Picard index that come with an effective 
complexity one torus action and have divisor 
class group $\ZZ$. 
Note that, by the approach, we get the Cox
rings of the resulting Fano varieties $X$ 
for free.
In Section~\ref{sec:tables}, we give 
some explicit classifications.
We list all non-toric surfaces $X$ with 
Picard index at most six and the 
non-toric threefolds~$X$ 
with Picard index up at most two.
They all have a Cox ring defined by a single 
relation; 
in fact, for surfaces the first 
Cox ring with more than one relation 
occurs for Picard index~29, and for the 
threefolds this happens with Picard index~3,
see Proposition~\ref{prop:fano22rel} 
as well as Examples~\ref{ex:fanosurf2rel} 
and~\ref{ex:fano32rel}.
Moreover, we determine all locally 
factorial fourfolds~$X$, i.e. 
those of Picard index one:
67 of them occur sporadic and 
there are two one-dimensional families.
Here comes the result on the locally factorial 
threefolds; in the table, we denote by $w_i$ 
the $\Cl(X)$-degree of the variable $T_i$.

\goodbreak

\begin{introthm2}
The following table lists the 
Cox rings $\mathcal{R}(X)$ 
of the three-dimensional 
locally factorial non-toric Fano 
varieties $X$ with an effective two torus 
action and $\Cl(X) = \ZZ$.
\begin{center}
\begin{longtable}[htbp]{llll}
\toprule
No.
&
$\mathcal{R}(X)$ 
& 
$(w_1,\ldots, w_5)$
& 
$(-K_X)^3$
\\
\midrule
1 
\hspace{.5cm}
&
$
\KK[T_1, \ldots, T_5] 
\ / \ 
\bangle{T_1T_2^5 + T_3^3 + T_4^2}
$
\hspace{.5cm}
&
$(1,1,2,3,1)$
\hspace{.5cm}
&
$8$
\\
\midrule
2
&
$
\KK[T_1, \ldots, T_5] 
\ / \ 
\bangle{T_1T_2T_3^4 + T_4^3 + T_5^2}
$
&
$(1,1,1,2,3)$
&
$8$
\\
\midrule
3
&
$
\KK[T_1, \ldots, T_5] 
\ / \ 
\bangle{T_1T_2^2T_3^3 + T_4^3 + T_5^2}
$
&
$(1,1,1,2,3)$
&
$8$
\\
\midrule
4
&
$
\KK[T_1, \ldots, T_5] 
\ / \ 
\bangle{T_1T_2 + T_3T_4 + T_5^2}
$
&
$(1,1,1,1,1)$
&
$54$
\\
\midrule
5
&
$
\KK[T_1, \ldots, T_5] 
\ / \ 
\bangle{T_1T_2^2 + T_3T_4^2 + T_5^3}
$
&
$(1,1,1,1,1)$
&
$24$
\\
\midrule
6
&
$
\KK[T_1, \ldots, T_5] 
\ / \ 
\bangle{T_1T_2^3 + T_3T_4^3 + T_5^4}
$
&
$(1,1,1,1,1)$
&
$4$
\\
\midrule
7
&
$
\KK[T_1, \ldots, T_5] 
\ / \ 
\bangle{T_1T_2^3 + T_3T_4^3 + T_5^2}
$
&
$(1,1,1,1,2)$
&
$16$
\\
\midrule
8
&
$
\KK[T_1, \ldots, T_5] 
\ / \ 
\bangle{T_1T_2^5 + T_3T_4^5 + T_5^2}
$
&
$(1,1,1,1,3)$
&
$2$
\\
\midrule
9
&
$
\KK[T_1, \ldots, T_5] 
\ / \ 
\bangle{T_1T_2^5 + T_3^3T_4^3 + T_5^2}
$
&
$(1,1,1,1,3)$
&
$2$
\\
\bottomrule 
\end{longtable}
\end{center}
\end{introthm2}

Note that each of these varieties $X$ is 
a hypersurface in the respective 
weighted projective space 
$\PP(w_1, \ldots, w_5)$.
Except number~4, none of them is 
quasismooth in the sense that 
$\Spec \, \mathcal{R}(X)$ is 
singular at most in the origin; 
quasismooth hypersurfaces of 
weighted projective spaces were studied 
in~\cite{JoKo} and~\cite{CCC}.
In Section~\ref{sec:geom3folds}, we take 
a closer look at the singularities 
of the threefolds listed above.
It turns out that number~1,3,5,7 and 9 are singular 
with only canonical singularities and all of them 
admit a crepant resolution.
Number~6 and 8 are  singular with  
non-canonical singularities but
admit a smooth relative minimal model.
Number two is singular with only canonical singularities,
one of them of type $\mathbf{cA_1}$,
and it admits only a singular relative minimal model.
Moreover, in all cases, we determine the Cox rings 
of the resolutions.


The authors would like to thank Ivan Arzhantsev 
for helpful comments and discussions and also the 
referee for valuable remarks and many references.

\section{UFDs with complexity one multigrading}
\label{sec:factrings}

As mentioned before, we work over an algebraically 
closed field $\KK$ of characteristic zero.
In Theorem~\ref{thm:factrings}, we describe all
factorial finitely generated $\KK$-algebras~$R$ 
with an effective complexity one grading and $R_0=\KK$.
Moreover, we characterize the possible Cox rings
among these algebras, see Proposition~\ref{prop:coxchar}.
First we recall the construction sketched in the 
introduction.

\begin{construction}
\label{constr:triple2ring}
Consider a sequence
$A = (a_0, \ldots, a_r)$ 
of vectors $a_i = (b_i,c_i)$ in $\KK^2$
such that any pair $(a_i,a_k)$ with
$k \ne i$ is linearly independent,
a sequence
$\mathfrak{n} = (n_0, \ldots, n_r)$ 
of positive integers
and a family $L = (l_{ij})$ 
of positive integers,
where $0 \le i \le r$ and 
$1 \le j \le n_i$.
For every $0 \le i \le r$, define a monomial
$$
f_i 
\  := \
T_{i1}^{l_{i1}} \cdots T_{in_i}^{l_{in_i}}
\ \in \
\KK[T_{ij}; \; 0 \le i \le r, \; 1 \le j \le n_i],
$$
for any two indices $0 \le i,j \le r$,
set $\alpha_{ij} :=  \det(a_i,a_j)  =  b_ic_j-b_jc_i$
and for any three indices 
$0 \le i < j < k \le r$ define 
a trinomial
$$
g_{i,j,k} 
\ := \ 
\alpha_{jk}f_i
\ + \ 
\alpha_{ki}f_j
\ + \ 
\alpha_{ij}f_k
\ \in \
\KK[T_{ij}; \; 0 \le i \le r, \; 1 \le j \le n_i].
$$ 
Note that the coefficients of this trinomial are 
all nonzero. 
The triple $(A,\mathfrak{n},L)$ then defines a ring
$$
R(A,\mathfrak{n},L)
\ := \ 
\KK[T_{ij}; \; 0 \le i \le r, \; 1 \le j \le n_i] 
\ / \
\bangle{g_{i,i+1,i+2}; \; 0 \le i \le r-2}.
$$
\end{construction}

\begin{proposition}
\label{prop:RAnLnormal}
For every triple $(A,\mathfrak{n},L)$ 
as in~\ref{constr:triple2ring},
the ring $R(A,\mathfrak{n},L)$ is a 
normal complete intersection of 
dimension
$$
\dim \, R(A,\mathfrak{n},L)
\ = \ 
 n-r+1,
\qquad\qquad
n \ := \ n_0 + \ldots + n_r.
$$
\end{proposition}

\begin{lemma}
\label{lem:alltrins}
In the setting of~\ref{constr:triple2ring},
one has for any $0 \le i < j < k < l \le r$
the identities
$$ 
g_{i,k,l}
\ = \ 
\alpha_{kl} \cdot g_{i,j,k} + \alpha_{ik} \cdot g_{j,k,l},
\qquad\qquad
g_{i,j,l}
\ = \ 
\alpha_{jl} \cdot g_{i,j,k} + \alpha_{ij} \cdot g_{j,k,l}.
$$
In particular, every trinomial $g_{i,j,k}$, 
where $0 \le i < j < k \le r$
is contained in the ideal 
$\bangle{g_{i,i+1,i+2}; \; 0 \le i \le r-2}$.
\end{lemma}

\begin{proof}
The identities are easily obtained by direct 
computation;
note that for this one may assume 
$a_j = (1,0)$ and $a_k = (0,1)$.
The supplement then follows by repeated
application of the identities.
\end{proof}

\begin{lemma}
\label{lem:twotrinszero}
In the notation of~\ref{constr:triple2ring}
and~\ref{prop:RAnLnormal},
set $X := V(\KK^n,g_0,\ldots, g_{r-2})$, and let
$z \in X$.
If we have $f_i(z) = f_j(z) = 0$ for 
two $0 \le i < j \le r$, 
then $f_k(z) = 0$ holds for all 
$0 \le k \le r$.
\end{lemma}

\begin{proof}
If $i<k<j$ holds, then, according to 
Lemma~\ref{lem:alltrins},
we have $g_{i,k,j}(z)=0$,
which implies $f_k(z)=0$.
The cases $k<i$ and $j<k$ are
obtained similarly.
\end{proof}

\begin{proof}[Proof of Proposition~\ref{prop:RAnLnormal}]
Set $X := V(\KK^n; g_0, \ldots, g_{r-2})$,
where $g_i := g_{i,i+1,i+2}$.
Then we have to show that $X$ is 
a connected complete intersection with 
at most normal singularities.
In order to see that $X$ is connected,
set $\ell := \prod n_i \prod l_{ij}$ and 
$\zeta_{ij} : = \ell n_i^{-1}l_{ij}^{-1}$.
Then $X \subseteq \KK^n$ is invariant 
under the $\KK^*$-action given by
\begin{eqnarray*}
t \cdot z
& := & 
(t^{\zeta_{ij}} z_{ij})
\end{eqnarray*} 
and the point $0 \in \KK^n$ lies in the closure of 
any orbit $\KK^* \mal x  \subseteq X$,
which implies connectedness.
To proceed, consider the Jacobian $J_g$  
of $g := (g_0, \ldots, g_{r-2})$.
According to Serre's criterion, 
we have to show that the set of points of 
$z \in X$ with $J_g(z)$ not of full rank
is of codimension at least two in $X$.
Note that the Jacobian $J_g$ is of the shape
\begin{eqnarray*}
J_g
& = & 
\left(
\begin{array}{rrrrrcrrrrr}
\delta_{0 \, 0} & \delta_{0 \, 1} & \delta_{0 \, 2} & 0 &  & &&&&& 0
\\
0 & \delta_{1 \, 1} & \delta_{1 \, 2} & \delta_{1 \, 3} & 0 & &&&&&
\\
&&&&&  \vdots  &&&&& 
\\
\\
&&&&& & 0 & \delta_{r-3 \, r-3} & \delta_{r-3 \, r-2} & \delta_{r-3 \, r-1} & 0
\\
0 &&&&& &  & 0 & \delta_{r-2 \, r-2} & \delta_{r-2 \, r-1} & \delta_{r-2 \, r}
\end{array}
\right)
\end{eqnarray*} 
where $\delta_{ti}$ is a nonzero multiple
of the gradient $\delta_i := \grad \, f_i$.
Consider $z \in X$ 
with $J_g(z)$ not of full rank.
Then $\delta_i(z) = 0 = \delta_k(z)$ 
holds with some $0 \le i < k \le r$.
This implies 
$z_{ij} = 0 = z_{kl}$ 
for some $1 \le j \le n_i$ and 
$1 \le l \le n_k$.
Thus, we have $f_i(z) = 0 = f_k(z)$.
Lemma~\ref{lem:twotrinszero} gives
$f_s(z) = 0$, for all $0 \le s \le r$.
Thus, some coordinate
$z_{st}$ must vanish
for every $0 \le s \le r$.
This shows that $z$ belongs to a 
closed subset of $X$ having 
codimension at least two in $X$.
\end{proof}

\begin{lemma}
\label{lem:tijprime}
Notation as in~\ref{constr:triple2ring}.
Then the variable $T_{ij}$ defines a prime ideal
in  $R(A,\mathfrak{n},L)$ if and only if 
the numbers $\gcd(l_{k1}, \ldots, l_{kn_k})$,
where $k \ne i$, are pairwise
coprime.
\end{lemma}

\begin{proof}
We treat exemplarily $T_{01}$.
Using Lemma~\ref{lem:alltrins}, 
we see that the ideal of relations of 
$R(A,\mathfrak{n},L)$ can be 
presented as follows
\begin{eqnarray*}
\bangle{g_{s,s+1,s+2}; \; 0 \le s \le r-2}
& = & 
\bangle{g_{0,s,s+1}; \; 1 \le s \le r-1}.
\end{eqnarray*}
Thus, the ideal 
$\bangle{T_{01}} \subseteq R(A,\mathfrak{n},L)$ 
is prime if and only if the following 
binomial ideal is prime
$$
\mathfrak{a}
\ := \ 
\bangle{\alpha_{s+1 \, 0}f_s  + \alpha_{0  s}f_{s+1}; \;   1 \le s \le r-1}
\ \subseteq \
\KK[T_{ij}; \; (i,j) \ne (0,1)].
$$
Set $l_i := (l_{i1}, \ldots, l_{in_i})$.
Then the ideal $\mathfrak{a}$ is prime if and 
only if the following family 
can be complemented to a lattice basis
$$ 
(l_1,-l_2,0,\ldots,0),
\
\ldots,
\
(0,\ldots,0,l_{r-1},-l_r).
$$
This in turn is equivalent to the 
statement that the numbers 
$\gcd(l_{k1}, \ldots, l_{kn_k})$,
where $1 \le k \le r$, are pairwise
coprime. 
\end{proof}

\begin{definition}
We say that a triple $(A,\mathfrak{n},L)$ 
as in~\ref{constr:triple2ring} is
{\em admissible\/} if the 
numbers $\gcd(l_{i1}, \ldots, l_{in_i})$,
where $0 \le i \le r$, are pairwise
coprime.
\end{definition}

\begin{construction}
\label{constr:Kgrading}
Let $(A,\mathfrak{n},L)$ be an admissible triple
and consider the following free abelian groups
$$ 
E
\quad := \quad
\bigoplus_{i=0}^r \bigoplus_{j=1}^{n_i}  \ZZ \mal e_{ij},
\qquad
\qquad
K
\quad := \quad
\bigoplus_{j=1}^{n_0}  \ZZ \mal u_{0j}
\ \oplus \
\bigoplus_{i=1}^r \bigoplus_{j=1}^{n_i-1}  \ZZ \mal u_{ij}
$$
and define vectors 
$u_{in_i} 
:= 
u_{01} + \ldots + u_{0r} - u_{i1} - \ldots - u_{in_i-1}
\in K$.
Then there is an epimorphism $\lambda \colon E \to K$ 
fitting into a commutative diagram with exact rows
$$ 
\xymatrix{
0
\ar[rr]
&&
E
\ar[rr]_{\alpha}^{e_{ij} \mapsto l_{ij} e_{ij}}
\ar[d]^{\eta}_{e_{ij} \mapsto u_{ij}}
&&
E
\ar[rr]^{e_{ij} \mapsto \b{e}_{ij} \qquad}
\ar[d]^{\lambda}
&&
{\bigoplus_{i,j} \ZZ / l_{ij} \ZZ}
\ar[rr]
\ar@{<->}[d]^{\cong}
&&
0
\\
0
\ar[rr]
&&
K
\ar[rr]_{\beta}
&&
K
\ar[rr]
&&
{\bigoplus_{i,j} \ZZ / l_{ij} \ZZ}
\ar[rr]
&&
0
}
$$
Define a $K$-grading of 
$\KK[T_{ij}; \; 0 \le i \le r, \; 1 \le j \le n_i]$
by setting $\deg \, T_{ij} := \lambda(e_{ij})$.
Then every 
$f_i = T_{i1}^{l_{i1}} \cdots T_{in_i}^{l_{in_i}}$ 
is $K$-homogeneous of degree
$$
\deg \, f_i
\ = \ 
l_{i1} \lambda(e_{i1}) + \ldots + l_{in_i}\lambda(e_{in_i}) 
\ = \ 
l_{01} \lambda(e_{01}) + \ldots + l_{0n_0}\lambda(e_{0n_0}) 
\ \in \ 
K.
$$
Thus, the polynomials $g_{i,j,k}$ 
of~\ref{constr:triple2ring}
are all $K$-homogeneous of the same degree
and we obtain an effective 
$K$-grading of complexity one
of $R(A,\mathfrak{n},L)$. 
\end{construction}

\begin{proof}
Only for the existence of the commutative diagram there 
is something to show.
Write for short 
$l_i := (l_{i1}, \ldots, l_{in_i})$.
By the admissibility condition, 
the vectors
$v_i := (0, \ldots, 0,l_i,-l_{i+1},0,\ldots,0)$,
where $0 \le i \le r-1$,
can be completed to a lattice basis for $E$.
Consequently, we find an epimorphism 
$\lambda \colon E \to K$ having precisely
$\lin(v_0, \ldots, v_{r-1})$ as its kernel.
By construction, $\ker(\lambda)$ equals 
$\alpha(\ker(\eta))$. 
Using this, we obtain the induced 
morphism $\beta \colon K \to K$ and 
the desired properties.
\end{proof}

\begin{lemma}
\label{lem:pwnonassoc}
Notation as in~\ref{constr:Kgrading}.
Then $R(A,\mathfrak{n},L)_0 = \KK$ 
and $R(A,\mathfrak{n},L)^* = \KK^*$ hold.
Moreover, the $T_{ij}$ define pairwise nonassociated 
prime elements in $R(A,\mathfrak{n},L)$.
\end{lemma}

\begin{proof}
The fact that all elements of degree zero
are constant is due to the fact that all
degrees $\deg \, T_{ij} = u_{ij} \in K$ 
are non-zero and generate a pointed convex
cone in $K_{\QQ}$.
As a consequence, we obtain that all units
in $R(A,\mathfrak{n},L)$ are constant.
The $T_{ij}$ are prime by the admissibility 
condition and Lemma~\ref{lem:tijprime},
and they are pairwise nonassociated
because they have pairwise different degrees
and all units are constant.
\end{proof}

\goodbreak

\begin{theorem}
\label{thm:factrings}
Up to isomorphy, the finitely generated 
factorial $\KK$-algebras with an effective 
complexity one grading $R = \oplus_M R_u$ 
and $R_0 = \KK$ are 
\begin{enumerate}
\item 
the polynomial algebras $\KK[T_1, \ldots, T_d]$ 
with a grading $\deg(T_i) = u_i \in \ZZ^{d-1}$ 
such that $u_1, \ldots, u_d$ generate $\ZZ^{d-1}$ 
as a lattice and the convex cone on $\QQ^{d-1}$
generated by $u_1, \ldots, u_d$ is pointed,
\item
the $(K \times \ZZ^m)$-graded 
algebras $R(A,\mathfrak{n},L)[S_1,\ldots,S_m]$,
where $R(A,\mathfrak{n},L)$ is the $K$-graded
algebra defined by an admissible triple
$(A,\mathfrak{n},L)$ as in~\ref{constr:triple2ring}
and~\ref{constr:Kgrading} 
and $\deg\, S_j \in \ZZ^m$ 
is the $j$-th canonical base vector.
\end{enumerate}
\end{theorem}

\begin{proof}
We first show that for any admissible triple
$(A,\mathfrak{n},L)$ the ring $R(A,\mathfrak{n},L)$
is a unique factorization domain.
If $l_{ij} = 1$ holds for any two $i,j$,
then, by~\cite[Prop.~2.4]{HaSu}, 
the ring $R(A,\mathfrak{n},L)$ 
is the Cox ring of a space $\PP_1(A,\mathfrak{n})$ 
and hence is a unique factorization domain.

Now, let $(A,\mathfrak{n},L)$ be arbitrary
admissible data and let $\lambda \colon E \to K$ 
be an epimorphism
as in~\ref{constr:Kgrading}. 
Set $n := n_0 + \ldots + n_r$
and consider the diagonalizable groups
$$
\TT^n \ := \ \Spec \, \KK[E],
\qquad
H \ := \ \Spec \, \KK[K],
\qquad
H_0 \ := \ \Spec \, \KK[\oplus_{i,j} \ZZ / l_{ij} \ZZ].
$$
Then $\TT^n = (\KK^*)^n$ is the standard $n$-torus
and $H_0$ is the direct product 
of the cyclic subgroups 
$H_{ij} := \Spec \, \KK[\ZZ / l_{ij} \ZZ]$.
Moreover, the diagram in~\ref{constr:Kgrading}
gives rise to a commutative diagram 
with exact rows
$$ 
\xymatrix{
0
&&
{\TT^n}
\ar[ll]
&&
{\TT^n}
\ar[ll]_{(t_{ij}^{l_{ij}}) \mapsfrom (t_{ij})}
&&
\ar[ll]
H_0
&&
0
\ar[ll]
\\
0
&&
{H}
\ar[ll]
\ar[u]^{\imath}
&&
{H}
\ar[ll]
\ar[u]^{\jmath}
&&
H_0
\ar[ll]
\ar@{<->}[u]_{\cong}
&&
0
\ar[ll]
}
$$
where $t_{ij} = \chi^{e_{ij}}$ are the coordinates 
of $\TT^n$ corresponding to the 
characters $e_{ij} \in E$
and the maps 
$\imath$, $\jmath$ are the closed embeddings
corresponding to the epimorphisms $\eta$, 
$\lambda$ respectively.

Setting $\deg \, T_{ij} := e_{ij}$ defines 
an action of $\TT^n$ on 
$\KK^n = \Spec \, \KK[T_{ij}]$;
in terms of the coordinates $z_{ij}$ 
corresponding to $T_{ij}$
this action is given by 
$t \mal z = (t_{ij} z_{ij})$.
The torus $H$ acts effectively on $\KK^n$ via 
the embedding $\jmath \colon H \to \TT^n$.
The generic isotropy group of $H$ along 
$V(\KK^n,T_{ij})$ is the subgroup 
$H_{ij} \subseteq H$
corresponding to 
$K \to K/\lambda(E_{ij})$, 
where $E_{ij} \subseteq E$ 
denotes the sublattice generated 
by all $e_{kl}$ with $(k,l) \ne (i,j)$;
recall that we have 
$K/\lambda(E_{ij}) \cong \ZZ / l_{ij}\ZZ$.

Now, set $l_{ij}' := 1$ 
for any two $i,j$ and 
consider the spectra
$X := \Spec \, R(A,\mathfrak{n},L)$
and 
$X' := \Spec \, R(A,\mathfrak{n},L')$.
Then the canonical surjections 
$\KK[T_{ij}] \to R(A,\mathfrak{n},L)$
and 
$\KK[T_{ij}] \to R(A,\mathfrak{n},L')$
define embeddings
$X \to \KK^n$ and $X' \to \KK^n$.
These embeddings fit into 
the following commutative diagram 
$$ 
\xymatrix{
{\KK^n}
\ar@{<-}[rrr]_{\pi}^{(z_{ij}^{l_{ij}}) \mapsfrom (z_{ij})}
& 
& 
&
{\KK^n}
\\
X'
\ar@{<-}[rrr]
\ar[u]
& 
&
& 
X
\ar[u]
}
$$
The action of $H$ leaves $X$ invariant
and the induced $H$-action on $X$ 
is the one given by the $K$-grading of 
$R(A,\mathfrak{n},L)$.
Moreover, $\pi \colon \KK^n \to \KK^n$
is the quotient map for the induced action 
of $H_0 \subseteq H$ on $\KK^n$,
we have $X = \pi^{-1}(X')$, and hence 
the restriction $\pi \colon X \to X'$
is a quotient map for the induced action 
of $H_0$ on $X$.

Removing all subsets $V(X;T_{ij},T_{kl})$,
where $(i,j) \ne (k,l)$ from $X$, 
we obtain an open subset $U \subseteq X$.
By Lemma~\ref{lem:pwnonassoc},
the complement $X \setminus U$
is of codimension at least two
and each $V(U,T_{ij})$ is irreducible.
By construction, the only isotropy groups 
of the $H$-action on $U$ are 
the groups $H_{ij}$ of the points of 
$V(U,T_{ij})$.
The image $U' := \pi(U)$ is open in
$X'$,  
the complement $X' \setminus U'$
is as well of codimension at least two
and $H/H_0$ acts freely on $U'$. 
According to~\cite[Cor.~5.3]{KKV},
we have two exact sequences fitting 
into the following diagram
$$ 
\xymatrix{
& 
& 
1 
\ar[d]
&
\\
& 
& 
{\Pic}(U') 
\ar[d]^{\pi^*}
&
\\
1 \ar[r]
&
{\Chi(H_0)} \ar[r]^{\alpha}
&
{\Pic_{H_0}}(U) \ar[r]^{\beta}
\ar[d]^{\delta}
&
{\Pic}(U)
\\
& 
& 
{\prod_{i,j}} \Chi(H_{ij}) 
&
}
$$
Since $X'$ is factorial, the Picard group
$\Pic(U')$ is trivial
and we obtain that $\delta$ is injective.
Since $H_0$ is the direct product 
of the isotropy groups $H_{ij}$
of the Luna strata $V(U,T_{ij})$, 
we see that 
$\delta \circ \alpha$ is an isomorphism.
It follows that $\delta$ is surjective 
and hence an isomorphism.
This in turn shows that $\alpha$ is an 
isomorphism.
Now, every bundle on $U$ is $H$-linearizable.
Since $H_0$ acts as a subgroup of $H$,
we obtain that every bundle is $H_0$-linearizable.
It follows that $\beta$ is surjective and hence
$\Pic(U)$ is trivial.   
We conclude $\Cl(X) = \Pic(U) = 0$,
which means that $R(A,\mathfrak{n},L)$ admits unique 
factorization.

The second thing we have to show is that 
any finitely generated factorial $\KK$-algebra
$R$ with an effective complexity one multigrading 
satisfying $R_0 = \KK$ is as claimed.
Consider the action of the torus $G$ on 
$X = \Spec\, R$ defined by the multigrading,
and let $X_0 \subseteq X$ be the set of 
points having finite isotropy $G_x$.
Then~\cite[Prop~3.3]{HaSu}
provides a graded splitting
\begin{eqnarray*}
R 
& \cong & 
R'[S_1, \ldots, S_m],
\end{eqnarray*}
where the variables $S_j$ are identified with 
the homogeneous functions defining the prime 
divisors $E_j$ inside the boundary 
$X \setminus X_0$ and $R'$ is the ring of 
functions of $X_0$, which are invariant 
under the subtorus $G_0 \subseteq G$ 
generated by the generic isotropy groups 
$G_j$ of $E_j$. 

Since $R'_0 = R_0 = \KK$ holds, the orbit 
space $X_0/G$ has only constant functions and 
thus is a space $\PP_1(A,\mathfrak{n})$
as constructed in~\cite[Section~2]{HaSu}.
This allows us to proceed exactly as in the 
proof of Theorem~\cite[Thm~1.3]{HaSu} and 
gives $R' = R(A,\mathfrak{n},L)$.
The admissibility condition follows
from Lemma~\ref{lem:tijprime} and the 
fact that each $T_{ij}$ defines a prime 
element in  $R'$.
\end{proof}

\begin{remark}
\label{rem:mori}
Let $(A,\mathfrak{n},L)$ be an admissible triple 
with $\mathfrak{n} =(1,\ldots, 1)$.
Then $K = \ZZ$ holds, the admissibility condition
just means that the numbers $l_{ij}$ are pairwise coprime
and we have 
$$ 
\dim \, R(A,\mathfrak{n},L)
\ = \ 
n_0 + \ldots + n_r - r + 1
\ = \ 
2. 
$$
Consequently, for two-dimensional rings, 
Theorem~\ref{thm:factrings} specializes to Mori's
description of almost geometrically graded 
two-dimensional 
unique factorization domains provided in~\cite{Mo}.
\end{remark}

\begin{proposition}
\label{prop:coxchar}
Let $(A,\mathfrak{n},L)$ be an admissible triple,
consider the associated
$(K \times \ZZ^m)$-graded ring
$R(A,\mathfrak{n},L)[S_1, \ldots, S_m]$
as in Theorem~\ref{thm:factrings}
and let $\mu \colon K \times \ZZ^m \to K'$ be a surjection 
onto an abelian group $K'$. 
Then the following statements are equivalent.
\begin{enumerate}
\item
The $K'$-graded ring 
$R(A,\mathfrak{n},L)[S_1, \ldots, S_m]$ 
is the Cox ring of a projective variety $X'$ with 
$\Cl(X') \cong K'$.
\item
For every pair $i,j$ with $0 \le i \le r$ and 
$1 \le j \le n_i$, the group $K'$ is generated
by the elements $\mu(\lambda(e_{kl}))$ and $\mu(e_s)$,
where $(i,j) \ne (k,l)$ and $1 \le s \le m$,
for every $1 \le t  \le m$, the group $K'$ is generated
by the elements $\mu(\lambda(e_{ij}))$ and $\mu(e_s)$,
where $0 \le i \le r$, $1 \le j \le n_i$ and $s \ne t$,
and, finally the following 
cone is of full dimension in $K'_{\QQ}$:
$$
\bigcap_{(k,l)} \cone(\mu(\lambda(e_{ij})), \mu(e_s); \; (i,j) \ne (k,l))
\ \cap \ 
\bigcap_{t} \cone(\mu(\lambda(e_{ij})), \mu(e_s); \; s \ne t).
$$
\end{enumerate}
\end{proposition}

\begin{proof}
Suppose that~(i) holds,
let $p \colon \rq{X}' \to X'$ 
denote the universal torsor
and let $X'' \subseteq X'$ be the set 
of smooth points.
According to~\cite[Prop.~2.2]{Ha2}, 
the group $H' = \Spec \, \KK[K']$ acts 
freely on $p^{-1}(X'')$, which 
is a big open subset of the total 
coordinate space 
$\Spec \, R(A,\mathfrak{n},L)[S_1, \ldots, S_m]$.
This implies the first condition of~(ii).
Moreover, by~\cite[Prop.~4.1]{Ha2}, the displayed cone 
is the moving cone of $X'$ and 
hence of full dimension.
Conversely, if~(ii) holds, 
then the $K'$-graded ring 
$R(A,\mathfrak{n},L)[S_1, \ldots, S_m]$ 
can be made into a bunched ring and 
hence is the Cox ring of a projective variety,
use~\cite[Thm.~3.6]{Ha2}.
\end{proof}

\section{Bounds for Fano varieties}

We consider  $d$-dimensional Fano varieties 
$X$ that come with a complexity one torus 
action and have divisor class group 
$\Cl(X) \cong \ZZ$.
Then the Cox ring $\mathcal{R}(X)$ of 
$X$ is factorial~\cite[Prop.~8.4]{BeHa1}
and has an effective complexity one 
grading, 
which refines the $\Cl(X)$-grading,
see~\cite[Prop.~2.6]{HaSu}.
Thus, according to Theorem~\ref{thm:factrings}, 
it is of the form
\begin{eqnarray*}
\mathcal R(X)
& \cong &
\KK[T_{ij}; \; 0 \le i \le r, \; 1 \le j \le n_i][S_1,\ldots, S_m]
\ / \
\bangle{g_{i,i+1,i+2}; \; 0 \le i \le r-2},
\\
g_{i,j,k}
& := & 
\alpha_{jk} T_{i1}^{l_{i1}} \cdots T_{in_i}^{l_{in_i}}
\ + \
\alpha_{ki} T_{j1}^{l_{j1}} \cdots T_{jn_{j}}^{l_{jn_{j}}}
\ + \
\alpha_{ij}T_{k1}^{l_{k1}} \cdots T_{kn_{k}}^{l_{kn_{k}}}.
\end{eqnarray*}
Here, we may (and will) assume $n_0 \ge \ldots \ge n_r \ge 1$.
With $n := n_0 + \ldots + n_r$, we have 
$n + m = d + r$.
For the degrees  of the variables in
$\Cl(X) \cong \ZZ$, we write
$w_{ij} := \deg \, T_{ij}$ for 
$0 \le i \leq r$, $1 \le j \le n_i$
and $u_k = \deg \, S_k$ for $1 \le k \le m$.
Moreover, for $\mu \in \ZZ_{>0}$, we denote 
by $\xi(\mu)$ the number of primes in
$\{2, \ldots, \mu\}$.
The following result provides bounds for the 
discrete data of the Cox ring.

\begin{theorem}
\label{Th:FiniteIndex}
In the above situation,
fix the dimension $d = \dim(X)$
and the Picard index 
$\mu = [\Cl(X):\Pic(X)]$.
Then we have 
$$
u_k \ \le \ \mu \quad \text{for } 1 \le k \le m.
$$
Moreover, for the degree $\gamma$ 
of the relations, the weights
$w_{ij}$ and the exponents $l_{ij}$,
where $0 \le i \le r$ and $1 \le j \le n_i$
one obtains the following.
\begin{enumerate}
\item
Suppose that $r = 0,1$ holds. 
Then $n + m \le d+1$ holds 
and one has the bounds
$$
w_{ij} \ \le \ \mu 
\quad\text{for } 0 \le i \le r \text{ and } 1 \le j \le n_i,
$$
and the Picard index is given by
$$
\mu 
\ = \ 
\mathrm{lcm}(w_{ij},u_k; \;0 \le i \le r,  1 \le j \le n_i, 1 \le k \le m ).
$$
\item
Suppose that $r \ge 2$ and $n_0=1$ hold.
Then $r \le \xi(\mu)-1$ and $n=r+1$ and 
$m=d-1$ hold and one has 
$$
w_{i1} \ \le \ \mu^r
\quad
\text{for } 0 \le i \le r,
\qquad
l_{01} \cdots l_{r1} \ \mid \ \mu,
\qquad
l_{01} \cdots l_{r1} \ \mid \ \gamma \ \le \ \mu^{r+1},
$$
and the Picard index is given by
$$
\mu 
\ = \
\mathrm{lcm}(\gcd(w_{j1}; \; j \neq i), u_k;\; 0 \le i \le r, 1\le k \le m).
$$
\item
Suppose that $r \ge 2$ and $n_0 > n_1=1$ hold.
Then we may assume $l_{11} > \ldots > l_{r1} \ge 2$,
we have $r \le \xi(3d\mu)-1$ and
$n_0+m = d$ and the bounds
$$
w_{01},\ldots,w_{0n_0} \ \le \ \mu,
\qquad
l_{01},\ldots,l_{0n_0} \ < \ 6d\mu,
$$
$$
w_{11},l_{21} \ < \ 2d\mu,
\qquad
w_{21},l_{11} \ < \ 3d\mu,
$$
$$
w_{i1} \ < \  6d\mu,
\quad
l_{i1} \ < \  2d\mu
\quad
\text{for } 2 \le i \le r,
$$
$$
l_{11} \cdots  l_{r1} \ \mid \ \gamma \ < \ 6d\mu,
$$
and the Picard index is given by
$$
\mu
\ = \
\mathrm{lcm}(w_{0j}, \gcd(w_{11},\ldots,w_{r1}), u_k; \; 
1 \le j \le n_0, 1 \le k \le m ).
$$
\item
Suppose that $n_1 > n_2 = 1$ holds.
Then we may assume $l_{21} > \ldots > l_{r1} \ge 2$,
we have $r \le \xi(2(d+1)\mu)-1$ 
and $n_0+n_1+m = d+1$ and the bounds
$$
w_{ij} \ \le \ \mu
\quad
\text{for } i=0,1 \text{ and }  1 \le j \le n_i,
\qquad
w_{21} \ < \ (d+1)\mu,
$$
$$
w_{ij}, l_{ij} \ < \ 2(d+1)\mu
\quad
\text{for } 0 \le i \le r \text{ and } 1 \le j \le n_i,
$$
$$
l_{21} \cdots l_{r1} \ \mid \ \gamma \ < \ 2(d+1)\mu,
$$
and the Picard index is given by
$$
\mu
\ = \
\mathrm{lcm}(w_{ij}, u_k; \; 0 \le i\le 1,  1 \le j \le n_i, 1 \le k \le m).
$$
\item
Suppose that $n_2 > 1$ holds
and let $s$ be the maximal number with $n_{s}>1$.
Then one may assume $l_{s+1,1} > \ldots > l_{r1} \ge 2$,
we have $r \le \xi((d+2)\mu)-1$ and 
$n_0+ \ldots + n_s+m = d+s$ and the bounds
$$
w_{ij} \ \le \ \mu, 
\quad \text{for } 0 \le i \le s,
$$
$$
w_{ij}, l_{ij} \ <  \ (d+2)\mu
\quad \text{for } 0 \le i \le r  \text{ and } 1 \le j \le n_i,
$$
$$
l_{s+1,1} \cdots l_{r1} \ \mid \ \gamma \ < \ (d+2)\mu,
$$
and the Picard index is given by
$$
\mu
\ = \
\mathrm{lcm}(w_{ij}, u_k; \; 0 \le i \le s, 1 \le j \le n_i, 1 \le k \le m).
$$
\end{enumerate} 
\end{theorem}

Putting all the bounds of the theorem together, 
we obtain the following (raw) bound for the number 
of deformation types.

\begin{corollary}
\label{cor:finitefanos}
For any pair $(d,\mu) \in \ZZ^2_{>0}$,
the number $\delta(d,\mu)$ of different 
deformation types of $d$-dimensional 
Fano varieties with a complexity one 
torus action such that
$\Cl(X) \cong \ZZ$ and $[\Cl(X):\Pic(X)]=\mu$
hold is bounded by
\begin{eqnarray*} 
\delta(d,\mu)
& \le &
(6d\mu)^{2\xi(3d\mu)+d-2}\mu^{\xi(\mu)^2 + 2\xi((d+2)\mu)+2d+2}.
\end{eqnarray*}
\end{corollary}

\begin{proof}
By Theorem~\ref{Th:FiniteIndex} the 
discrete data $r$, $\mathfrak{n}$, 
$L$ and $m$ occuring in $\mathcal{R}(X)$ 
are bounded as in the assertion. 
The continuous data  in $\mathcal{R}(X)$ 
are the coefficients $\alpha_{ij}$; 
they stem from the family 
$A = (a_0, \ldots, a_r)$ 
of points $a_i \in \KK^2$.
Varying the $a_i$ provides flat 
families of Cox rings and hence,
by passing to the homogeneous spectra, 
flat families of the resulting 
Fano varieties $X$.
\end{proof}

\begin{corollary}
\label{cor:asymptoticsd}
Fix $d \in \ZZ_{>0}$. Then 
the number $\delta(\mu)$ of different 
deformation types of $d$-dimensional 
Fano varieties with a complexity one 
torus action, $\Cl(X) \cong \ZZ$ 
and  Picard index
$\mu := [\Cl(X):\Pic(X)]$
is asymptotically bounded by 
$\mu^{A \mu^2 / \log^2 \mu}$
with a constant~$A$ depending only 
on~$d$.
\end{corollary}

\begin{corollary}
\label{cor:asymptoticsmu}
Fix $\mu \in \ZZ_{>0}$. Then 
the number $\delta(d)$ of different 
deformation types of $d$-dimensional 
Fano varieties with a complexity one 
torus action, $\Cl(X) \cong \ZZ$ 
and  Picard index
$\mu := [\Cl(X):\Pic(X)]$
is asymptotically bounded by 
$d^{Ad}$ with a constant~$A$ 
depending only on~$\mu$.
\end{corollary}

We first recall the necessary facts on 
Cox rings, for details, we refer 
to~\cite{Ha2}.
Let $X$ be a complete $d$-dimensional 
variety with divisor class group 
$\Cl(X) \cong \ZZ$.
Then the Cox ring $\mathcal{R}(X)$
is finitely generated and the total 
coordinate space $\b{X} := \Spec \, \mathcal{R}(X)$ 
is a factorial affine variety coming 
with an action of $\KK^*$ defined by 
the $\Cl(X)$-grading of $\mathcal{R}(X)$. 
Choose a system $f_1, \ldots, f_\nu$ of 
homogeneous pairwise nonassociated 
prime generators for $\mathcal{R}(X)$.
This provides an $\KK^*$-equivariant 
embedding
$$ 
\b{X} \ \to \ \KK^{\nu},
\qquad
\b{x} \ \mapsto \ (f_1(\b{x}), \ldots, f_{\nu}(\b{x})).
$$
where $\KK^*$ acts  diagonally with the 
weights $w_i = \deg(f_i) \in \Cl(X) \cong \ZZ$ 
on $\KK^{\nu}$.
Moreover, $X$ is the geometric 
$\KK^*$-quotient of  
$\rq{X} := \b{X} \setminus \{0\}$,
and the quotient map 
$p \colon \rq{X} \to X$ is a universal 
torsor.
By the local divisor class group $\Cl(X,x)$
of a point $x \in X$, we mean the group of 
Weil divisors $\WDiv(X)$ modulo those that 
are principal near~$x$.

\begin{proposition}
\label{Prop:FanoPicard}
For any 
$\b{x} =(\b{x}_1,\ldots,\b{x}_{\nu}) \in \rq{X}$
the local divisor class group $\Cl(X,x)$ 
of $x := p(\b{x})$ 
is finite of order $\gcd(w_i; \; \b{x}_i \ne 0)$.
The index of the Picard group $\Pic(X)$ in
$\Cl(X)$ is given by  
\begin{eqnarray*}
[\Cl(X):\Pic(X)]
& = &
\mathrm{lcm}_{x \in X}( |\Cl(X,x)| ).
\end{eqnarray*}
Suppose that the ideal of $\b{X} \subseteq \KK^{\nu}$
is generated by $\Cl(X)$-homogeneous
polynomials $g_1, \ldots, g_{\nu-d-1}$
of degree $\gamma_j := \deg(g_j)$.
Then one obtains
$$ 
-\mathcal{K}_X 
\ = \ 
\sum_{i=1}^{\nu} w_i - \sum_{j=1}^{\nu-d-1} \gamma_j,
\qquad
(-\mathcal{K}_X )^d
\ = \ 
\left(\sum_{i=1}^{\nu} w_i - \sum_{j=1}^{\nu-d-1} \gamma_j\right)^d
\frac{\gamma_1 \cdots \gamma_{\nu-d-1}}{w_1 \cdots w_\nu}
$$
for the anticanonical class 
$-\mathcal{K}_X \in \Cl(X) \cong \ZZ$.
In particular, $X$ is a Fano variety
if and only if the following inequality holds
\begin{eqnarray*}
 \sum_{j=1}^{\nu-d-1} \gamma_j
 & < &
 \sum_{i=1}^{\nu} w_i.
\end{eqnarray*}
\end{proposition}

\begin{proof}
Using~\cite[Prop.~2.2, Thm.~4.19]{Ha2}, we observe
that $X$ arises from the bunched ring 
$(R,\mathfrak{F},\Phi)$, 
where $R = \mathcal{R}(X)$, 
$\mathfrak{F} = (f_1, \ldots, f_\nu)$
and $\Phi = \{\QQ_{\ge 0}\}$.
The descriptions of local class groups, the 
Picard index and the anticanonical class are 
then special cases 
of~\cite[Prop.~4.7, Cor.~4.9 and Cor.~4.16]{Ha2}.
The anticanonical self-intersection number 
is easily computed in the ambient weighted 
projective space $\PP(w_1, \ldots, w_\nu)$,
use~\cite[Constr.~3.13, Cor.~4.13]{Ha2}. 
\end{proof}

\begin{remark}
If the ideal of $\b{X} \subseteq \KK^{\nu}$
is generated by $\Cl(X)$-homogeneous
polynomials $g_1, \ldots, g_{\nu-d-1}$,
then~\cite[Constr.~3.13, Cor.~4.13]{Ha2} 
show that $X$ is a well formed 
complete intersection in the weighted 
projective space $\PP(w_1, \ldots, w_\nu)$ in the
sense of~\cite[Def.~6.9]{IaFl}.
\end{remark}

We turn back to the case that $X$ comes
with a complexity one torus action as at
the beginning of this section.
We consider the case $n_0 = \ldots = n_r=1$,
that means that each relation $g_{i,j,k}$
of the Cox ring $\mathcal{R}(X)$ depends
only on three variables.
Then we may write $T_i$ instead of $T_{i1}$ 
and $w_i$ instead of $w_{i1}$, etc..
In this setting, 
we obtain the following bounds 
for the numbers of possible varieties~$X$
(Fano or not).

\begin{proposition}
\label{prop:Finite3Var}
For any pair $(d,\mu) \in \ZZ^2_{>0}$
there is, up to deformation, 
only a finite number of 
complete $d$-dimensional
varieties with divisor class group
$\ZZ$,
Picard index $[\Cl(X):\Pic(X)] = \mu$
and Cox ring 
$$ 
\KK[T_0,\ldots, T_r,S_1,\ldots, S_m]
\ / \
\bangle{
\alpha_{i+1,i+2} T_i^{l_i}
+
\alpha_{i+2,i} T_{i+1}^{l_{i+1}}
+
\alpha_{i,i+1} T_{i+2}^{l_{i+2}}; 
\; 0 \le i \le r-2}.
$$
In this situation we have $r \le \xi(\mu)-1$.
Moreover, for the weights $w_i := \deg\, T_i$,
where $0 \le i \le r$
and $u_k := \deg\, S_k$, where $1 \le k \le m$,
the exponents $l_i$ 
and the degree $\gamma := l_0w_0$ of the relation
one has
$$
l_0 \cdots l_r \ \mid \ \gamma,
\qquad
l_0 \cdots l_r \ \mid \ \mu,
\qquad 
w_i \ \le \ \mu^{\xi(\mu)-1},
\qquad
u_k \ \le \ \mu.
$$
\end{proposition}

\begin{proof}
Consider the total coordinate space 
$\b{X} \subseteq \KK^{r+1+n}$ and 
the universal torsor $p \colon \rq{X} \to X$ 
as discussed before. 
For each $0 \le i \le r$ fix a point 
$\b{x}(i) = (\b{x}_0, \ldots, \b{x}_r, 0, \ldots, 0)$
in $\rq{X}$ such that $\b{x}_i = 0$ and 
$\b{x}_j \ne 0$ for $j \ne i$ hold.
Then, denoting $x(i) := p(\b{x}(i))$, we obtain
$$ 
\gcd(w_j; j \ne i)
\ = \ 
\vert \Cl(X,x(i)) \vert 
\ \mid \ \mu.
$$
Consider $i,j$ with $j \ne i$. 
Since all relations are homogeneous
of the same degree,
we have $l_iw_i = l_jw_j$.
Moreover, by the admissibility condition, 
$l_i$ and $l_j$ are coprime.
We conclude $l_i \vert w_j$
for all $j \ne i$ and hence 
$l_i \vert \gcd(w_j; \; j \ne i)$.
This implies
$$ 
l_0 \cdots l_r \ \mid \ l_0w_0 \ = \ \gamma,
\qquad\qquad
l_0 \cdots l_r \ \mid \ \mu.
$$
We turn to the bounds for the $w_i$, 
and first verify $w_0 \le \mu^r$.
Using the relation $l_iw_i = l_0w_0$,
we obtain  for every $l_i$ a presentation
$$ 
l_i
\ = \ 
l_0 \cdot \frac{w_0 \cdots w_{i-1}}{w_1 \cdots w_i}
\ = \ 
\eta_i \cdot \frac{\gcd(w_0, \ldots ,w_{i-1})}{\gcd(w_0, \ldots, w_i)}
$$
with suitable integers $1 \le \eta_i \le \mu$.
In particular, the very last fraction 
is bounded by $\mu$. 
This gives the desired estimate:   
$$ 
w_0 
 = 
\frac{w_{0}}{\gcd(w_{0},w_{1})}
\cdot
\frac{\gcd(w_{0},w_1)}{\gcd(w_{0},w_{1},w_2)}
\cdots
\frac
{\gcd(w_{0},\ldots,w_{r-2})}
{\gcd(w_{0},\ldots,w_{r-1})}
\cdot
\gcd(w_{0},\ldots,w_{r-1})
 \le  
\mu^r.
$$
Similarly, we obtain $w_i \le \mu^r$ 
for $1 \le i \le r$.
Then we only have to show that 
$r+1$ is bounded by $\xi(\mu)$,
but this follows immediately from
the fact that $l_0, \ldots, l_r$ 
are pairwise coprime.

Finally, to estimate the $u_k$,
consider the points $\b{x}(k) \in \rq{X}$
having the $(r+k)$-th coordinate one and 
all others zero. Set $x(k) : =p(\b{x}(k))$.
Then $\Cl(X,x(k))$ is of order $u_k$, 
which implies $u_k \le \mu$.
\end{proof}

\goodbreak

\begin{lemma}
\label{Lem:1relation}
Consider the ring
$\KK[T_{ij}; \; 0 \le i \le 2, \; 1 \le j \le n_i][S_1,\ldots,S_k]
 / 
\bangle{g}
$
where $n_0 \ge n_1 \ge n_2 \ge 1$ holds.
Suppose that $g$ is homogeneous 
with respect to a $\ZZ$-grading 
of $\KK[T_{ij},S_k]$ given by 
$\deg \, T_{ij} = w_{ij} \in \ZZ_{>0}$ 
and $\deg \, S_k = u_k \in \ZZ_{>0}$,
and assume 
\begin{eqnarray*}
\deg \, g
& < & 
\sum_{i=0}^2\sum_{j=1}^{n_i}w_{ij}
\ + \
\sum_{i=1}^m u_i.
\end{eqnarray*} 
Let $\mu \in \ZZ_{>1}$, assume 
$w_{ij} \le \mu$ whenever $n_i > 1$,
$1 \le j \le n_i$ and $u_k \le \mu$ 
for $1 \le k \le m$ and 
set $d := n_0+n_1+n_2+m-2$.
Depending on the shape of $g$, 
one obtains the following bounds.
\begin{enumerate}
\item
Suppose that
$g
= 
\eta_0 T_{01}^{l_{01}}  \cdots  T_{0n_0}^{l_{0n_0}} 
+  
\eta_1 T_{11}^{l_{11}}
+
\eta_2 T_{21}^{l_{21}}$ 
with $n_0 > 1$ and coefficients 
$\eta_i \in \KK^*$ holds,
we have $l_{11} \ge l_{21} \ge 2$ 
and $l_{11}$, $l_{21}$ are coprime.
Then, one has
$$ 
\qquad\qquad
w_{11}, l_{21} \ < \ 2d\mu,
\qquad 
w_{21}, l_{11} \ < \ 3d\mu,
\qquad
\deg \, g  \ < \ 6d\mu.
$$
\item
Suppose that
$g
=
\eta_0 T_{01}^{l_{01}}  \cdots T_{0n_0}^{l_{0n_0}}
+
\eta_1 T_{11}^{l_{11}} \cdots  T_{1n_1}^{l_{1n_1}}
+
\eta_2 T_{21}^{l_{21}}$
with $n_1 > 1$ and coefficients 
$\eta_i \in \KK^*$ holds and 
we have $l_{21} \ge 2$.
Then one has
$$
\qquad\qquad
w_{21}
\ < \ 
(d+1)\mu,
\qquad\qquad
\deg \, g 
\ < \ 
2(d+1)\mu.
$$
\end{enumerate}
\end{lemma}

\begin{proof}
We prove~(i). Set for short
$c := (n_0+m)\mu = d\mu$.
Then, using homogeneity of $g$ 
and the assumed inequality, we obtain
$$
l_{11}w_{11}
\ = \ 
l_{21}w_{21}
\ = \ 
\deg \, g
\ < \ 
\sum_{i=0}^2\sum_{j=1}^{n_i}w_{ij}
+
\sum_{i=1}^m u_i 
\ \le \ 
c+w_{11}+w_{21}.
$$
Since $l_{11}$ and $l_{21}$ 
are coprime, we have 
$l_{11} > l_{21} \ge 2$.
Plugging this into the above inequalities, 
we arrive at 
$2 w_{11} < c + w_{21}$ and 
$w_{21} < c + w_{11}$.
We conclude $w_{11} < 2c$ and $w_{21} < 3c$.
Moreover, $l_{11}w_{11} = l_{21}w_{21}$ and 
$\gcd(l_{11},l_{21}) = 1$ imply 
$l_{11} \vert w_{21}$ and $l_{21} \vert w_{11}$.
This shows $l_{11} < 3c$ and $l_{21} < 2c$.
Finally, we obtain
$$
\deg \, g 
\ <  \ 
c + w_{11} + w_{21} 
\ <  \ 
6c.
$$
We prove (ii).
Here we set $c := (n_0+n_1+m)\mu = (d+1)\mu$.
Then the assumed inequality gives
$$
l_{21}w_{21}
\ = \ 
\deg g
\ < \
\sum_{i=0}^1\sum_{j=1}^{n_i}w_{ij}+
\sum_{i=1}^m u_i+ w_{21}
\ \le \
c+w_{21}.
$$
Since we assumed $l_{21} \geq 2$, 
we can conclude $w_{21} < c$.
This in turn gives us 
$\deg \, g < 2c$
for the degree of the relation.
\end{proof}

\begin{proof}
[Proof of Theorem~\ref{Th:FiniteIndex}]
As before, we denote by $\b{X} \subseteq \KK^{n+m}$
the total coordinate space and by 
$p \colon \rq{X} \to X$ the universal torsor.

We first consider the case that 
$X$ is a toric variety.
Then the Cox ring is a polynomial ring,
$\mathcal R(X) = \KK[S_1,\ldots,S_m]$.
For each $1 \le k \le m$, 
consider the point
$\overline x(k) \in \rq{X}$
having the $k$-th coordinate one 
and all others zero and 
set $x(k) := p(\b{x}(k))$.
Then, by~Proposition~\ref{Prop:FanoPicard},
the local class group
$\Cl(X,x(k))$ is of order $u_k$
where $u_k := \deg \, S_k$.
This implies 
$u_k \le \mu$ for $1 \le k \leq m$
and settles Assertion~(i).

Now we treat the non-toric case,
which means $r \ge 2$.
Note that we have $n \ge 3$.
The case $n_0=1$ is done in 
Proposition~\ref{prop:Finite3Var}.
So, we are left with $n_0>1$.
For every $i$ with $n_i > 1$
and every $1 \le j \le n_i$,
there is the point $\b{x}(i,j) \in \rq{X}$
with $ij$-coordinate $T_{ij}$ equal 
to one and all others equal to
zero, and thus we have the point
$x(i,j) := p(\b{x}(i,j)) \in X$.
Moreover, for every $1 \le k \le m$, we have 
the point $\b{x}(k) \in \rq{X}$
having the $k$-coordinate $S_k$ equal to one 
and all others zero; we 
set $x(k):=p(\b{x}(k))$.
Proposition~\ref{Prop:FanoPicard}
provides the bounds
$$ 
w_{ij}
\ = \ 
\deg \, T_{ij}
\ = \ 
\vert \Cl(X,x(i,j)) \vert 
\ \le \ 
\mu
\qquad
\text{for } 
n_i > 1, \, 1 \le j \le n_i,
$$
$$
u_k
\ = \ 
\deg \, S_k
\ = \ 
\vert \Cl(X,x(k)) \vert 
\ \le \
\mu
\qquad
\text{for } 
1 \le k \le m.
$$

Let $0 \le s \le r$ be the maximal number with 
$n_{s} > 1$. Then $g_{s-2,s-1,s}$ is the last 
polynomial such that each of its three monomials
depends on more than one variable.
For any $t \ge s$, we have the ``cut ring'' 
\begin{eqnarray*}
R_t
& := &
\KK[T_{ij}; \; 0 \le i \le t, \; 1 \le j \le n_i]
[S_1,\ldots,S_m] 
\ / \
\bangle{g_{i,i+1,i+2}; \; 0 \le i \le t-2}
\end{eqnarray*}
where the relations $g_{i,i+1,i+2}$ depend on
only three variables as soon as $i > s$ holds.
For the degree $\gamma$ of the relations we have 
\begin{eqnarray*}
(r-1)\gamma
& = & 
(t-1)\gamma \ + \ (r-t)\gamma
\\
& = &
(t-1)\gamma \ + \ l_{t+1,1}w_{t+1,1} + \ldots + l_{r1}w_{r1}
\\
& < & 
\sum_{i=0}^r\sum_{j=1}^{n_i}w_{ij}
\ + \ 
\sum_{i=1}^m u_i
\\
& = & 
\sum_{i=0}^t \sum_{j=1}^{n_i}w_{ij}
\ + \ 
w_{t+1,1}+ \ldots + w_{r1}
\ + \ 
\sum_{i=1}^m u_i.
\end{eqnarray*}
Since $l_{i1}w_{i1} > w_{i1}$ holds in 
particular for $t+1 \le i \le r$,
we derive from this the inequality
\begin{eqnarray*}
\gamma
& < &
\frac{1}{t-1}
\left(
\sum_{i=0}^t\sum_{j=1}^{n_i}w_{ij}
\ + \ 
\sum_{i=1}^m u_i
\right).
\end{eqnarray*}

To obtain the bounds in 
Assertions~(iii) and~(iv),
we consider the cut ring $R_t$ 
with $t=2$ and apply
Lemma~\ref{Lem:1relation};
note that we have 
$d = n_0+n_1+n_2+m-2$ 
for the dimension $d = \dim(X)$
and that $l_{22} \ge 0$ is due 
to the fact that $X$ is non-toric.
The bounds $w_{ij}, l_{0j} < 6d\mu$ in 
Assertion~(iii) follow from 
$l_{ij}w_{ij} = \gamma < 6 d\mu$
and $l_{i1} < 2d\mu$ follows from 
$l_{i1} \mid w_{21}$ for $3 \le i \le r$.
Moreover, $l_{i1} \mid w_{11}$ 
for $2 \le i \le r$ implies
$l_{11} \cdots l_{r1} \mid \gamma = l_{11}w_{11}$.
Similarly  $w_{ij},l_{ij} < 2(d+1)\mu$
in Assertion~(iv) follow from 
$l_{ij}w_{ij} = \gamma < 2(d+1) d\mu$
and $l_{21} \cdots l_{r1} \mid \gamma = l_{21}w_{21}$
follows from $l_{i1} \mid w_{21}$ for $3 \le i \le r$.
The bounds on $r$ in~(iii) in~(iv) are 
as well consequences of the admissibility 
condition.

To obtain the bounds in Assertion~(v),
we consider the cut ring $R_t$ 
with $t=s$. 
Using $n_i=1$ for $i \ge t+1$, we 
can estimate the degree of the relation
as follows:
$$
\gamma
\ \le \
\frac{(n_0 + \ldots + n_t + m) \mu}{t-1}
\ = \ 
\frac{(d + t) \mu}{t-1}
\ \le \ 
(d + 2) \mu.
$$
Since we have $w_{ij}l_{ij} \le \deg \, g_0$ 
for any $0 \le i \le r$ and any $1 \le j \le n_i$,
we see that all $w_{ij}$ and $l_{ij}$ 
are bounded by $(d+2)\mu$.
As before, $l_{s+1,1} \cdots l_{r1} \mid \gamma$ 
is a consequence of $l_{i1} \mid \gamma$ 
for $i = s+2, \ldots, r$ 
and also the bound on $r$ follows 
from the admissibility condition.

Finally, we have to express the Picard index
$\mu$ in terms of the weights $w_{ij}$ and $u_k$
as claimed in the Assertions.
This is a direct application of the formula of 
Proposition~\ref{Prop:FanoPicard}.
Observe that it suffices to work with the 
$p$-images of the following points:
For every $0 \le i \le r$ with $n_i > 1$ 
take a point $\b{x}(i,j) \in \rq{X}$
with $ij$-coordinate $T_{ij}$ equal 
to one and all others equal to
zero, 
for every $0 \le i \le r$ with $n_i = 1$ 
whenever $n_i=1$ take $\b{x}(i,j) \in \rq{X}$
with $ij$-coordinate $T_{ij}$ equal 
to zero, all other $T_{st}$ equal to
one and coordinates $S_k$ equal to zero,  
and, for every $1 \le k \le m$,
take a point $\b{x}(k) \in \rq{X}$
having the $k$-coordinate $S_k$ equal to one 
and all others zero.
\end{proof}

We conclude the section with discussing some 
aspects of the not necessarily Fano 
varieties of Proposition~\ref{prop:Finite3Var}.
Recall that we considered admissible triples 
$(A,\mathfrak{n},L)$ with 
$n_0 = \ldots = n_r =1$ and thus 
rings $R$ of the form
$$ 
\KK[T_0,\ldots, T_r,S_1,\ldots, S_m]
\ / \
\bangle{\alpha_{i+1,i+2} T_i^{l_i}
+
\alpha_{i+2,i} T_{i+1}^{l_{i+1}}
+
\alpha_{i,i+1} T_{i+2}^{l_{i+2}}; 
\; 0 \le i \le r-2}.
$$

\begin{proposition}
\label{prop:MoriCox}
Suppose that the ring $R$ as above 
is the Cox ring 
of a non-toric variety $X$ 
with $\Cl(X) = \ZZ$. 
Then we have $m \ge 1$ and 
$\mu := [\Cl(X):\Pic(X)] \ge 30$.
Moreover, if $X$ is a surface, then we 
have $m=1$ and $w_i= l_i^{-1} l_0 \cdots l_r$.
\end{proposition}

\begin{proof}
The homogeneity condition 
$l_{i}w_{i}=l_{j}w_{j}$
together with the
admissibility condition
$\gcd(l_{i},l_{j})=1$
for $0 \le  i \ne j\leq r$
gives us
$l_{i} \mid \gcd(w_{j}; j \ne i)$.
Moreover, by Proposition~\ref{prop:coxchar},
every set of $m+r$ weights $w_i$
has to generate the class group $\ZZ$,
so they must have
greatest common divisor one.
Since $X$ is non-toric,
$l_{i} \ge 2$ holds and we obtain
$m \ge 1$.
To proceed, we infer
$l_0 \cdots l_r \mid \mu$ 
and  
$l_0 \cdots l_r \mid \deg g_{ijk}$
from Proposition~\ref{Prop:FanoPicard}. 
As a consequence, 
the minimal value for $\mu$ and $\deg g_{ijk}$
is obviously $2\cdot3\cdot5=30$.
 what really can be received
as the following example shows.
Note that if $X$ is a surface we have $m=1$ and
$\gcd(w_{i}; 0\le i\le r )=1$.
Thus, $l_iw_i=l_jw_j$ gives us 
$\deg g_{ijk}= l_0 \cdots l_r$ and
$w_i= l_i^{-1} l_0 \cdots l_r$.
\end{proof}

The bound $[\Cl(X):\Pic(X)] \ge 30$ given 
in the above proposition is even sharp;
the surface discussed below realizes it.

\begin{example}
Consider $X$ with 
$\mathcal R(X)=
\KK[T_{0},T_{1},T_{2},T_{3}]/
\langle g\rangle$
with
$g=T_{0}^2+T_{1}^3+T_{2}^5$
and the grading
$$
\deg \, T_0 \ = \ 15,
\quad
\deg \, T_1 \ = \ 10,
\quad
\deg \, T_2 \ = \ 6,
\quad
\deg \, T_3 \ = \ 1.
$$
Then we have $\gcd(15,10)=5$, $\gcd(15,6)=3$
and $\gcd(10,6)=2$
and therefore $[\Cl(X):\Pic(X)]=30$.
Further $X$ is Fano because of
$$
\deg \, g
\ = \ 
30
\ < \
32
\ = \ 
\deg \, T_0 + \ldots + \deg \, T_3.
$$
\end{example}

Let us have a look at the geometric meaning 
of the condition $n_0 = \ldots = n_r = 1$.
For a variety $X$ with an action of a torus $T$,
we denote by $X_0 \subseteq X$ the 
union of all orbits with at most finite isotropy.
Then there is a possibly non-separated orbit space
$X_0/T$; we call it the maximal orbit space. 
From~\cite{HaSu}, we infer that 
$n_0 = \ldots = n_r = 1$ holds if and only 
if $X_0/T$ is separated. 
Combining this with Propositions~\ref{prop:Finite3Var}
and~\ref{prop:MoriCox} gives the following.

\begin{corollary}
For any pair $(d,\mu) \in \ZZ^2_{>0}$
there is, up to deformation, 
only a finite number of 
$d$-dimensional complete 
varieties $X$ 
with a complexity one torus action
having divisor class group $\ZZ$, 
Picard index $[\Cl(X):\Pic(X)] = \mu$
and maximal orbit space $\PP_1$
and for each of these varieties 
the complement
$X \setminus X_0$ contains divisors.
\end{corollary}

Finally, we present a couple of examples
showing that there are also non-Fano
varieties with a complexity one torus action
having divisor class group $\ZZ$
and maximal orbit space $\PP_1$.

\begin{example}
Consider $X$ with 
$\mathcal R(X)= \KK[T_{0},T_{1},T_{2},T_{3}]/ \langle g\rangle$
with $g=T_{0}^2+T_{1}^3+T_{2}^{7}$
and the grading
$$
\deg \, T_0 \ = \ 21,
\quad
\deg \, T_1 \ = \ 14,
\quad
\deg \, T_2 \ = \ 6,
\quad
\deg \, T_3 \ = \ 1.
$$
Then we have $\gcd(21,14)=7$, $\gcd(21,6)=3$
and $\gcd(14,6)=2$
and therefore $[\Cl(X):\Pic(X)]=42$.
Moreover, $X$ is not Fano, 
because its canonical class $\mathcal{K}_X$ is 
trivial
$$ 
\mathcal{K}_X
\ = \ 
\deg \, g - \deg \, T_0 - \ldots -\deg \, T_3 
\ = \ 
0.
$$
\end{example}
 
\begin{example}
Consider $X$ with 
$\mathcal R(X)= \KK[T_{0},T_{1},T_{2},T_{3}]/ \langle g\rangle$
with $g=T_{0}^2+T_{1}^3+T_{2}^{11}$
and the grading
$$
\deg \, T_0 \ = \ 33,
\quad
\deg \, T_1 \ = \ 22,
\quad
\deg \, T_2 \ = \ 6,
\quad
\deg \, T_3 \ = \ 1.
$$
Then we have $\gcd(22,33)=11$, $\gcd(33,6)=3$
and $\gcd(22,6)=2$
and therefore $[\Cl(X):\Pic(X)]=66$.
The canonical class $\mathcal{K}_X$ of 
$X$ is even ample:
$$ 
\mathcal{K}_X
\ = \ 
 \deg \, g - \deg \, T_0 - \ldots - \deg \, T_3 
\ = \ 
4.
$$
\end{example}

The following example shows 
that the Fano assumption is 
essential for the finiteness results 
in Theorem~\ref{Th:FiniteIndex}.

\begin{remark}
For any pair 
$p,q$ of coprime positive integers,
we obtain a locally factorial 
$\KK^*$-surface $X(p,q)$ with 
$\Cl(X) = \ZZ$ and Cox ring 
$$
\mathcal{R}(X(p,q))
\ = \ 
\KK[T_{01},T_{02},T_{11},T_{21}]
\ / \ 
\bangle{g},
\qquad\qquad
g \ = \ 
T_{01}T_{02}^{pq-1}
+T_{11}^{q}
+T_{21}^{p};
$$
the $\Cl(X)$-grading is given by 
$\deg \, T_{01} = \deg \, T_{02} =1 $, 
$\deg \, T_{11} = p$ and $\deg \, T_{21} = q$. 
Note that $\deg \, g =pq$ holds 
and for $p,q\geq 3$, the canonical class 
$\mathcal{K}_X$ satisfies  
$$
\mathcal{K}_X 
\ = \ 
\deg \, g - \deg \, T_{01} - \deg \, T_{02} - \deg \, T_{11} - \deg \, T_{21}
\ = \
pq - 2 - p - q
\ \ge \ 
0.
$$
\end{remark}

\section{Classification results}
\label{sec:tables}

In this section, we give classification results 
for Fano varieties~$X$ with $\Cl(X) \cong \ZZ$ 
that come with a complexity one torus action;
note that they are necessarily rational.
The procedure to obtain classification
lists for prescribed dimension $d = \dim \, X$ 
and Picard index $\mu = [\Cl(X) : \Pic(X)]$
is always the following.
By Theorem~\ref{thm:factrings}, we know that their 
Cox rings are of the form 
$\mathcal{R}(X) \cong R(A,\mathfrak{n},L)[S_1,\ldots,S_m]$
with admissible triples $(A,\mathfrak{n},L)$.
Note that for the family 
$A = (a_0, \ldots, a_r)$ of points $a_i \in \KK^2$,
we may assume 
$$ 
a_0 \ = \ (1,0),
\qquad
a_1 \ = \ (1,1),
\qquad
a_2 \ = \ (0,1).
$$
The bounds on the input data of $(A,\mathfrak{n},L)$
provided by Theorem~\ref{Th:FiniteIndex}
as well as the criteria of 
Propositions~\ref{prop:coxchar}
and~\ref{Prop:FanoPicard}
allow us to generate all the possible Cox rings
$\mathcal{R}(X)$ of the Fano varieties 
$X$ in question for fixed dimension~$d$
and Picard index~$\mu$.
Note that $X$ can be reconstructed from 
$\mathcal{R}(X)= R(A,\mathfrak{n},L)[S_1,\ldots,S_n]$ 
as the homogeneous spectrum 
with respect to the $\Cl(X)$-grading.
Thus $X$ is classified by its Cox ring
$\mathcal{R}(X)$.

In the following tables, we present the 
Cox rings as $\KK[T_1, \ldots, T_s]$ 
modulo relations and fix the 
$\ZZ$-gradings by giving the weight vector 
$(w_1, \ldots, w_s)$, where $w_i := \deg \, T_i$.
The first classification result concerns 
surfaces.

\begin{theorem}
Let $X$ be a non-toric Fano surface
with an effective $\KK^*$-action such that 
$\Cl(X)=\ZZ$ and $[\Cl(X):\Pic(X)]\leq 6$ hold.
Then its Cox ring is precisely one of the following.

\begin{center}
\begin{longtable}[htbp]{llll}
\multicolumn{4}{c}{\bf $[\Cl(X):\Pic(X)] = 1$}
\\[1ex]
\toprule
No.
&
$\mathcal{R}(X)$
&
$(w_1,\ldots,w_4)$
&
$(-K_X)^2$
\\
\midrule
1
\hspace{.5cm}
&
$\KK[{T_1,\ldots,T_4}]/ \bangle{T_1T_2^5+T_3^3+T_4^2}$
\hspace{.5cm}
&
$(1,1,2,3)$
\hspace{.5cm}
&
$1$
\\
\bottomrule
\\[2ex]
\multicolumn{4}{c}{\bf $[\Cl(X):\Pic(X)] = 2$}
\\[1ex]
\toprule
No.
&
$\mathcal{R}(X)$
&
$(w_1,\ldots,w_4)$
&
$(-K_X)^2$
\\
\midrule
2
&
$\KK[{T_1,\ldots,T_4}]/ \bangle{T_1^4T_2+T_3^3+T_4^2}$
&
$(1,2,2,3)$
&
$2$
\\
\bottomrule
\\[2ex]
\multicolumn{4}{c}{\bf $[\Cl(X):\Pic(X)] = 3$}
\\[1ex]
\toprule
No.
&
$\mathcal{R}(X)$
&
$(w_1,\ldots,w_4)$
&
$(-K_X)^2$
\\
\midrule
3
&
$\KK[{T_1,\ldots,T_4}]/ \bangle{T_1^3T_2+T_3^3+T_4^2}$
&
$(1,3,2,3)$
&
$3$
\\
\midrule
4
&
$\KK[{T_1,\ldots,T_4}] / \bangle{T_1T_2^3+T_3^5+T_4^2}$
&
$(1,3,2,5)$
&
$1/3$
\\
\midrule
5
&
$\KK[{T_1,\ldots,T_4}]/ \bangle{T_1^7T_2+T_3^5+T_4^2}$
&
$(1,3,2,5)$
&
$1/3$
\\
\bottomrule
\\[2ex]
\multicolumn{4}{c}{\bf $[\Cl(X):\Pic(X)] = 4$}
\\[1ex]
\toprule
No.
&
$\mathcal{R}(X)$
&
$(w_1,\ldots,w_4)$
&
$(-K_X)^2$
\\
\midrule
6
&
$\KK[{T_1,\ldots,T_4}] / \bangle{T_1^2T_2+T_3^3+T_4^2}$
&
$(1,4,2,3)$
&
$4$
\\
\midrule
7
&
$\KK[{T_1,\ldots,T_4}] / \bangle{T_1^6T_2+T_3^5+T_4^2}$
&
$(1,4,2,5)$
&
$1$
\\
\bottomrule
\\[2ex]
\multicolumn{4}{c}{\bf $[\Cl(X):\Pic(X)] = 5$}
\\[1ex]
\midrule
No.
&
$\mathcal{R}(X)$
&
$(w_1,\ldots,w_4)$
&
$(-K_X)^2$
\\
\midrule
8
&
$\KK[{T_1,\ldots,T_4}] / \bangle{T_1T_2+T_3^3+T_4^2}$
&
$(1,5,2,3)$
&
$5$
\\
\midrule
9
&
$\KK[{T_1,\ldots,T_4}] /  \bangle{T_1^5T_2+T_3^5+T_4^2}$
&
$(1,5,2,5)$
&
$9/5$
\\
\midrule
10
&
$\KK[{T_1,\ldots,T_4}] / \bangle{T_1^9T_2+T_3^7+T_4^2}$
&
$(1,5,2,7)$
&
$1/5$
\\
\midrule
11
&
$\KK[{T_1,\ldots,T_4}] / \bangle{T_1^7T_2+T_3^4+T_4^3}$
&
$(1,5,3,4)$
&
$1/5$
\\
\bottomrule
\\[2ex]
\multicolumn{4}{c}{\bf $[\Cl(X):\Pic(X)] = 6$}
\\[1ex]
\toprule
No.
&
$\mathcal{R}(X)$
&
$(w_1,\ldots,w_4)$
&
$(-K_X)^2$
\\
\midrule
12
&
$\KK[{T_1,\ldots,T_4}] / \bangle{T_1^4T_2+T_3^5+T_4^2}$
&
$(1,6,2,5)$
&
$8/3$
\\
\midrule
13
&
$\KK[{T_1,\ldots,T_4}] / \bangle{T_1^8T_2+T_3^7+T_4^2}$
&
$(1,6,2,7)$
&
$2/3$
\\
\midrule
14
&
$\KK[{T_1,\ldots,T_4}] / \bangle{T_1^6T_2+T_3^4+T_4^3}$
&
$(1,6,3,4)$
&
$2/3$
\\
\midrule
15
&
$\KK[{T_1,\ldots,T_4}] / \bangle{T_1^9T_2+T_3^3+T_4^2}$
&
$(1,3,4,6)$
&
$2/3$
\\
\bottomrule
\end{longtable}
\end{center}
\end{theorem}

\begin{proof}
As mentioned, Theorems~\ref{thm:factrings},
\ref{Th:FiniteIndex} and 
Propositions~\ref{prop:coxchar}, \ref{Prop:FanoPicard}
produce a list of all Cox rings of surfaces
with the prescribed data.
Doing this computation, we obtain the list 
of the assertion.
Note that none of the Cox rings listed 
is a polynomial ring and hence none of the 
resulting surfaces $X$ is a toric variety.
To show that different members of the 
list are not isomorphic to each other, 
we use the following two facts.
Firstly, observe that any two minimal 
systems of homogeneous generators of 
the Cox ring have (up to reordering)
the same list of degrees, and thus
the list of generator degrees is invariant 
under isomorphism  (up to reordering).
Secondly, by Construction~\ref{constr:Kgrading}, 
the exponents $l_{ij} >1$ are precisely the 
orders of the non-trivial isotropy groups 
of one-codimensional orbits of the action 
of the torus $T$ on $X$.
Using both principles and going through the 
list, we see that different members $X$ 
cannot be $T$-equivariantly isomorphic to 
each other.
Since all listed $X$ are non-toric,
the effective complexity one torus action
on each $X$ corresponds to a maximal torus in
the linear algebraic group $\Aut(X)$.
Any two maximal tori in the automorphism
group are conjugate, and thus we can conclude 
that two members are isomorphic if and only if they 
are $T$-equivariantly isomorphic.
\end{proof}

We remark that in~\cite[Section~4]{tfano}, 
log del Pezzo surfaces with an effective 
$\KK^*$-action and Picard number 1 and 
Gorenstein index less than 4 were classified.
The above list contains six such surfaces,
namely no. 1-4, 6 and~8;
these are exactly the ones where 
the maximal exponents of the monomials 
form a platonic triple, i.e., 
are of the form $(1,k,l)$, $(2,2,k)$, $(2,3,3)$, 
$(2,3,4)$ or $(2,3,5)$.
The remaining ones, i.e., no. 5, 7, and~9-15
have non-log-terminal and thus non-rational
singularities; to check this one may compute
the resolutions via resolution of the ambient 
weighted projective space as 
in~\cite[Ex.~7.5]{Ha2}.

With the same scheme of proof 
as in the surface case, one establishes 
the following classification results 
on Fano threefolds.

\goodbreak

\begin{theorem}
\label{thm:3fano}
Let $X$ be a three-dimensional 
locally factorial non-toric Fano 
variety with an effective two torus 
action such that $\Cl(X) = \ZZ$ holds.
Then its Cox ring is precisely 
one of the following.

\begin{center}
\begin{longtable}[htbp]{llll}
\toprule
No.
&
$\mathcal{R}(X)$ 
& 
$(w_1,\ldots, w_5)$
& 
$(-K_X)^3$
\\
\midrule
1 
\hspace{.5cm}
&
$
\KK[T_1, \ldots, T_5] 
\ / \ 
\bangle{T_1T_2^5 + T_3^3 + T_4^2}
$
\hspace{.5cm}
&
$(1,1,2,3,1)$
\hspace{.5cm}
&
$8$
\\
\midrule
2
&
$
\KK[T_1, \ldots, T_5] 
\ / \ 
\bangle{T_1T_2T_3^4 + T_4^3 + T_5^2}
$
&
$(1,1,1,2,3)$
&
$8$
\\
\midrule
3
&
$
\KK[T_1, \ldots, T_5] 
\ / \ 
\bangle{T_1T_2^2T_3^3 + T_4^3 + T_5^2}
$
&
$(1,1,1,2,3)$
&
$8$
\\
\midrule
4
&
$
\KK[T_1, \ldots, T_5] 
\ / \ 
\bangle{T_1T_2 + T_3T_4 + T_5^2}
$
&
$(1,1,1,1,1)$
&
$54$
\\
\midrule
5
&
$
\KK[T_1, \ldots, T_5] 
\ / \ 
\bangle{T_1T_2^2 + T_3T_4^2 + T_5^3}
$
&
$(1,1,1,1,1)$
&
$24$
\\
\midrule
6
&
$
\KK[T_1, \ldots, T_5] 
\ / \ 
\bangle{T_1T_2^3 + T_3T_4^3 + T_5^4}
$
&
$(1,1,1,1,1)$
&
$4$
\\
\midrule
7
&
$
\KK[T_1, \ldots, T_5] 
\ / \ 
\bangle{T_1T_2^3 + T_3T_4^3 + T_5^2}
$
&
$(1,1,1,1,2)$
&
$16$
\\
\midrule
8
&
$
\KK[T_1, \ldots, T_5] 
\ / \ 
\bangle{T_1T_2^5 + T_3T_4^5 + T_5^2}
$
&
$(1,1,1,1,3)$
&
$2$
\\
\midrule
9
&
$
\KK[T_1, \ldots, T_5] 
\ / \ 
\bangle{T_1T_2^5 + T_3^3T_4^3 + T_5^2}
$
&
$(1,1,1,1,3)$
&
$2$
\\
\bottomrule 
\end{longtable}
\end{center}
\end{theorem}

The singular threefolds listed in this theorem 
are rational degenerations of smooth Fano threefolds 
from~\cite{fano3}.
The (smooth) general Fano threefolds of 
the corresponding families are non-rational 
see~\cite{Gri} for no.~1-3,
\cite{CG} for no.~5,
\cite{IM} for no.~6, 
\cite{voi,tim}~for no.~7 
and \cite{Isk80} for no. 8-9.
Even if one allows certain mild singularities, 
one still has non-rationality in some cases, 
see \cite{Gri2}, \cite{Co,Pu}, \cite{CM}, \cite{CP}.

\begin{theorem}
\label{thm:3fano2}
Let $X$ be a three-dimensional non-toric Fano 
variety
with an effective two torus action such that 
$\Cl(X)=\ZZ$ and $[\Cl(X):\Pic(X)]=2$ hold.
Then its Cox ring is precisely one of the 
following.
\begin{center}
\begin{longtable}[htbp]{llll}
\toprule
No. \hspace{.5cm}
&
$\mathcal{R}(X)$ \hspace{.5cm}
&
$(w_1,\ldots,w_5)$ \hspace{.5cm}
& 
$(-K_X)^3$
\\
\midrule 1 & 
$\KK[{T_1,\ldots,T_5}]/
 \langle   T_1^4T_2+T_3^3+T_4^2  \rangle$
&
 $(1,2,2,3,1)$
&
$27/2$
\\
\midrule 2 &
$\KK[{T_1,\ldots,T_5}]/
 \langle   T_1^4T_2^3+T_3^5+T_4^2  \rangle$
&
 $(1,2,2,5,1)$
&
$1/2$
\\
\midrule 3 &
$\KK[{T_1,\ldots,T_5}]/
 \langle   T_1^8T_2+T_3^5+T_4^2  \rangle$
&
 $(1,2,2,5,1)$
&
$1/2$
\\
\midrule 4 &
$\KK[{T_1,\ldots,T_5}]/
 \langle   T_1^4T_2+T_3^3+T_4^2  \rangle$
&
 $(1,2,2,3,2)$
&
$16$
\\
\midrule 5 &
$\KK[{T_1,\ldots,T_5}]/
 \langle   T_1^4T_2^3+T_3^5+T_4^2  \rangle$
&
 $(1,2,2,5,2)$
&
$2$
\\
\midrule 6 &
$\KK[{T_1,\ldots,T_5}]/
 \langle   T_1^8T_2+T_3^5+T_4^2  \rangle$
&
 $(1,2,2,5,2)$
&
$2$
\\
\midrule 7 &
$\KK[{T_1,\ldots,T_5}]/
 \langle   T_1T_2^5+T_3^3+T_4^2  \rangle$
&
 $(1,1,2,3,2)$
&
$27/2$
\\
\midrule 8 &
$\KK[{T_1,\ldots,T_5}]/
 \langle   T_1T_2^9+T_3^5+T_4^2  \rangle$
&
 $(1,1,2,5,2)$
&
$1/2$
\\
\midrule 9 &
$\KK[{T_1,\ldots,T_5}]/
 \langle   T_1^3T_2^7+T_3^5+T_4^2  \rangle$
&
 $(1,1,2,5,2)$
&
$1/2$
\\
\midrule 10 &
$\KK[{T_1,\ldots,T_5}]/
 \langle   T_1T_2^{11}+T_3^3+T_4^2  \rangle$
&
 $(1,1,4,6,1)$
&
$1/2$
\\
\midrule 11 &
$\KK[{T_1,\ldots,T_5}]/
 \langle   T_1^5T_2^7+T_3^3+T_4^2  \rangle$
&
 $(1,1,4,6,1)$
&
$1/2$
\\
\midrule 12 &
$\KK[{T_1,\ldots,T_5}]/
 \langle   T_1T_2^{11}+T_3^3+T_4^2  \rangle$
&
 $(1,1,4,6,2)$
&
$2$
\\
\midrule 13 &
$\KK[{T_1,\ldots,T_5}]/
 \langle   T_1^5T_2^7+T_3^3+T_4^2  \rangle$
&
 $(1,1,4,6,2)$
 &
$2$
\\
\midrule 14 &
$\KK[{T_1,\ldots,T_5}]/
 \langle   T_1^2T_2^5+T_3^3+T_4^2  \rangle$
&
 $(1,2,4,6,1)$
&
$2$
\\
\midrule 15 &
$\KK[{T_1,\ldots,T_5}]/
 \langle   T_1^{10}T_2+T_3^3+T_4^2  \rangle$
&
 $(1,2,4,6,1)$
&
$2$
\\
\midrule 16 &
$\KK[{T_1,\ldots,T_5}]/
 \langle   T_1T_2^2+T_3^3+T_4^2  \rangle$
&
 $(2,2,2,3,1)$
&
$16$
\\
\midrule 17 &
$\KK[{T_1,\ldots,T_5}]/
 \langle   T_1T_2^4+T_3^5+T_4^2  \rangle$
&
 $(2,2,2,5,1)$
&
$2$
\\
\midrule 18 &
$\KK[{T_1,\ldots,T_5}]/
 \langle   T_1^2T_2^3+T_3^5+T_4^2  \rangle$
&
 $(2,2,2,5,1)$
&
$2$
\\
\midrule 19 &
$\KK[{T_1,\ldots,T_5}]/
 \langle   T_1T_2^2+T_3T_4+T_5^3  \rangle$
&
 $(1,1,1,2,1)$
&
$81/2$
\\
\midrule 20 &
$\KK[{T_1,\ldots,T_5}]/
 \langle   T_1T_2^4+T_3T_4^2+T_5^5  \rangle$
&
 $(1,1,1,2,1)$
&
$5/2$
\\
\midrule 21 &
$\KK[{T_1,\ldots,T_5}]/
 \langle   T_1^2T_2^3+T_3T_4^2+T_5^5  \rangle$
&
 $(1,1,1,2,1)$
&
$5/2$
\\
\midrule 22 &
$\KK[{T_1,\ldots,T_5}]/
 \langle   T_1T_2^3+T_3^2T_4+T_5^4  \rangle$
&
 $(1,1,1,2,1)$
&
$16$
\\
\midrule 23 &
$\KK[{T_1,\ldots,T_5}]/
 \langle   T_1T_2^4+T_3^3T_4+T_5^5  \rangle$
&
 $(1,1,1,2,1)$
&
$5/2$
\\
\midrule 24 &
$\KK[{T_1,\ldots,T_5}]/
 \langle   T_1^2T_2^3+T_3^3T_4+T_5^5  \rangle$
&
 $(1,1,1,2,1)$
&
$5/2$
\\
\midrule 25 &
$\KK[{T_1,\ldots,T_5}]/
 \langle   T_1T_2^3+T_3^2T_4+T_5^2  \rangle$
&
 $(1,1,1,2,2)$
&
$27$
\\
\midrule 26 &
$\KK[{T_1,\ldots,T_5}]/
 \langle   T_1T_2^5+T_3^2T_4^2+T_5^3  \rangle$
&
 $(1,1,1,2,2)$
&
$3/2$
\\
\midrule 27 &
$\KK[{T_1,\ldots,T_5}]/
 \langle   T_1T_2^5+T_3^4T_4+T_5^3  \rangle$
&
 $(1,1,1,2,2)$
&
$3/2$
\\
\midrule 28 &
$\KK[{T_1,\ldots,T_5}]/
 \langle   T_1^2T_2^4+T_3^4T_4+T_5^3  \rangle$
&
 $(1,1,1,2,2)$
&
$3/2$
\\
\midrule 29 &
$\KK[{T_1,\ldots,T_5}]/
 \langle   T_1T_2^5+T_3^4T_4+T_5^2  \rangle$
&
 $(1,1,1,2,3)$
&
$8$
\\
\midrule 30 &
$\KK[{T_1,\ldots,T_5}]/
 \langle   T_1^3T_2^3+T_3^4T_4+T_5^2  \rangle$
&
 $(1,1,1,2,3)$
&
$8$
\\
\midrule 31 &
$\KK[{T_1,\ldots,T_5}]/
 \langle   T_1T_2^7+T_3^2T_4^3+T_5^2  \rangle$
&
 $(1,1,1,2,4)$
&
$1$
\\
\midrule 32 &
$\KK[{T_1,\ldots,T_5}]/
 \langle   T_1^3T_2^5+T_3^2T_4^3+T_5^2  \rangle$
&
 $(1,1,1,2,4)$
&
$1$
\\
\midrule 33 &
$\KK[{T_1,\ldots,T_5}]/
 \langle   T_1T_2^7+T_3^6T_4+T_5^2  \rangle$
&
 $(1,1,1,2,4)$
&
$1$
\\
\midrule 34 &
$\KK[{T_1,\ldots,T_5}]/
 \langle   T_1^3T_2^5+T_3^6T_4+T_5^2  \rangle$
&
 $(1,1,1,2,4)$
&
$1$
\\
\midrule 35 &
$\KK[{T_1,\ldots,T_5}]/
 \langle   T_1T_2^3+T_3T_4+T_5^4  \rangle$
&
 $(1,1,2,2,1)$
&
$27$
\\
\midrule 36 &
$\KK[{T_1,\ldots,T_5}]/
 \langle   T_1T_2^5+T_3T_4^2+T_5^6  \rangle$
&
 $(1,1,2,2,1)$
&
$3/2$
\\
\midrule 37 &
$\KK[{T_1,\ldots,T_5}]/
 \langle   T_1T_2^3+T_3T_4+T_5^2  \rangle$
&
 $(1,1,2,2,2)$
&
$16$
\\
\midrule 38 &
$\KK[{T_1,\ldots,T_5}]/
 \langle   T_1T_2^5+T_3T_4^2+T_5^3  \rangle$
&
 $(1,1,2,2,2)$
&
$6$
\\
\midrule 39 &
$\KK[{T_1,\ldots,T_5}]/
 \langle   T_1^2T_2^4+T_3T_4^2+T_5^3  \rangle$
&
 $(1,1,2,2,2)$
&
$6$
\\
\midrule 40 &
$\KK[{T_1,\ldots,T_5}]/
 \langle   T_1^3T_2^3+T_3T_4^2+T_5^2  \rangle$
&
 $(1,1,2,2,2)$
&
$27/2$
\\
\midrule 41 &
$\KK[{T_1,\ldots,T_5}]/
 \langle   T_1^3T_2^5+T_3T_4^3+T_5^2  \rangle$
&
 $(1,1,2,2,2)$
&
$32$
\\
\midrule 42 &
$\KK[{T_1,\ldots,T_5}]/
 \langle   T_1T_2^5+T_3T_4^2+T_5^2  \rangle$
&
 $(1,1,2,2,3)$
&
$4$
\\
\midrule 43 &
$\KK[{T_1,\ldots,T_5}]/
 \langle   T_1T_2^7+T_3T_4^3+T_5^2  \rangle$
&
 $(1,1,2,2,4)$
&
$32$
\\
\midrule 44 &
$\KK[{T_1,\ldots,T_5}]/
 \langle   T_1T_2^9+T_3T_4^4+T_5^2  \rangle$
&
 $(1,1,2,2,5)$
&
$1/2$
\\
\midrule 45 &
$\KK[{T_1,\ldots,T_5}]/
 \langle   T_1T_2^9+T_3^2T_4^3+T_5^2  \rangle$
&
 $(1,1,2,2,5)$
&
$1/2$
\\
\midrule 46 &
$\KK[{T_1,\ldots,T_5}]/
 \langle   T_1^3T_2^7+T_3T_4^4+T_5^2  \rangle$
&
 $(1,1,2,2,5)$
&
$1/2$
\\
\midrule 47 &
$\KK[{T_1,\ldots,T_5}]/
 \langle   T_1^3T_2^7+T_3^2T_4^3+T_5^2  \rangle$
&
 $(1,1,2,2,5)$
&
$1/2$
\\
\midrule 48 &
$\KK[{T_1,\ldots,T_5}]/
 \langle   T_1^5T_2^5+T_3T_4^4+T_5^2  \rangle$
&
 $(1,1,2,2,5)$
&
$1/2$
\\
\midrule 49 &
$\KK[{T_1,\ldots,T_5}]/
 \langle   T_1^5T_2^5+T_3^2T_4^3+T_5^2  \rangle$
&
 $(1,1,2,2,5)$
&
$1/2$
\\
\midrule 50 &
$\KK[{T_1,\ldots,T_5}]/
 \langle   T_1T_2+T_3T_4+T_5^3  \rangle$
&
 $(1,2,1,2,1)$
&
$48$
\\
\midrule 51 &
$\KK[{T_1,\ldots,T_5}]/
 \langle   T_1^2T_2+T_3^2T_4+T_5^4  \rangle$
&
 $(1,2,1,2,1)$
&
$27$
\\
\midrule 52 &
$\KK[{T_1,\ldots,T_5}]/
 \langle   T_1T_2^2+T_3T_4^2+T_5^5  \rangle$
&
 $(1,2,1,2,1)$
&
$10$
\\
\midrule 53 &
$\KK[{T_1,\ldots,T_5}]/
 \langle   T_1T_2^2+T_3^3T_4+T_5^5  \rangle$
&
 $(1,2,1,2,1)$
&
$10$
\\
\midrule 54 &
$\KK[{T_1,\ldots,T_5}]/
 \langle   T_1^3T_2+T_3^3T_4+T_5^5  \rangle$
&
 $(1,2,1,2,1)$
&
$10$
\\
\midrule 55 &
$\KK[{T_1,\ldots,T_5}]/
 \langle   T_1^4T_2+T_3^4T_4+T_5^6  \rangle$
&
 $(1,2,1,2,1)$
&
$3/2$
\\
\midrule 56 &
$\KK[{T_1,\ldots,T_5}]/
 \langle   T_1^2T_2+T_3^2T_4+T_5^2  \rangle$
&
 $(1,2,1,2,2)$
&
$32$
\\
\midrule 57 &
$\KK[{T_1,\ldots,T_5}]/
 \langle   T_1^2T_2^2+T_3^4T_4+T_5^3  \rangle$
&
 $(1,2,1,2,2)$
&
$6$
\\
\midrule 58 &
$\KK[{T_1,\ldots,T_5}]/
 \langle   T_1^4T_2+T_3^4T_4+T_5^3  \rangle$
&
 $(1,2,1,2,2)$
&
$6$
\\
\midrule 59 &
$\KK[{T_1,\ldots,T_5}]/
 \langle   T_1^4T_2+T_3^4T_4+T_5^2  \rangle$
&
 $(1,2,1,2,3)$
&
$27/2$
\\
\midrule 60 &
$\KK[{T_1,\ldots,T_5}]/
 \langle   T_1^2T_2^3+T_3^2T_4^3+T_5^2  \rangle$
&
 $(1,2,1,2,4)$
&
$4$
\\
\midrule 61 &
$\KK[{T_1,\ldots,T_5}]/
 \langle   T_1^2T_2^3+T_3^6T_4+T_5^2  \rangle$
&
 $(1,2,1,2,4)$
&
$4$
\\
\midrule 62 &
$\KK[{T_1,\ldots,T_5}]/
 \langle   T_1^6T_2+T_3^6T_4+T_5^2  \rangle$
&
 $(1,2,1,2,4)$
&
$4$
\\
\midrule 63 &
$\KK[{T_1,\ldots,T_5}]/
 \langle   T_1^4T_2^3+T_3^4T_4^3+T_5^2  \rangle$
&
 $(1,2,1,2,5)$
&
$1/2$
\\
\midrule 64 &
$\KK[{T_1,\ldots,T_5}]/
 \langle   T_1^8T_2+T_3^4T_4^3+T_5^2  \rangle$
&
 $(1,2,1,2,5)$
&
$1/2$
\\
\midrule 65 &
$\KK[{T_1,\ldots,T_5}]/
 \langle   T_1^8T_2+T_3^8T_4+T_5^2  \rangle$
&
 $(1,2,1,2,5)$
&
$1/2$
\\
\midrule 66 &
$\KK[{T_1,\ldots,T_5}]/
 \langle   T_1^2T_2+T_3T_4+T_5^4  \rangle$
&
 $(1,2,2,2,1)$
&
$32$
\\
\midrule 67 &
$\KK[{T_1,\ldots,T_5}]/
 \langle   T_1^4T_2+T_3T_4^2+T_5^6  \rangle$
&
 $(1,2,2,2,1)$
&
$6$
\\
\midrule 68 &
$\KK[{T_1,\ldots,T_5}]/
 \langle   T_1^4T_2+T_3T_4^2+T_5^2  \rangle$
&
 $(1,2,2,2,3)$
&
$16$
\\
\midrule 69 &
$\KK[{T_1,\ldots,T_5}]/
 \langle   T_1^4T_2^3+T_3T_4^4+T_5^2  \rangle$
&
 $(1,2,2,2,5)$
&
$2$
\\
\midrule 70 &
$\KK[{T_1,\ldots,T_5}]/
 \langle   T_1^4T_2^3+T_3^2T_4^3+T_5^2  \rangle$
&
 $(1,2,2,2,5)$
&
$2$
\\
\midrule 71 &
$\KK[{T_1,\ldots,T_5}]/
 \langle   T_1^8T_2+T_3T_4^4+T_5^2  \rangle$
&
 $(1,2,2,2,5)$
&
$2$
\\
\midrule 72 &
$\KK[{T_1,\ldots,T_5}]/
 \langle   T_1^8T_2+T_3^2T_4^3+T_5^2  \rangle$
&
 $(1,2,2,2,5)$
&
$2$
\\
\midrule 73 &
$\KK[{T_1,\ldots,T_5}]/
 \langle   T_1T_2T_3^{10}+T_4^3+T_5^2  \rangle$
&
 $(1,1,1,4,6)$
&
$1/2$
\\
\midrule 74 &
$\KK[{T_1,\ldots,T_5}]/
 \langle   T_1T_2^2T_3^9+T_4^3+T_5^2  \rangle$
&
 $(1,1,1,4,6)$
&
$1/2$
\\
\midrule 75 &
$\KK[{T_1,\ldots,T_5}]/
 \langle   T_1T_2^3T_3^8+T_4^3+T_5^2  \rangle$
&
 $(1,1,1,4,6)$
&
$1/2$
\\
\midrule 76 &
$\KK[{T_1,\ldots,T_5}]/
 \langle   T_1T_2^4T_3^7+T_4^3+T_5^2  \rangle$
&
 $(1,1,1,4,6)$
&
$1/2$
\\
\midrule 77 &
$\KK[{T_1,\ldots,T_5}]/
 \langle   T_1T_2^5T_3^6+T_4^3+T_5^2  \rangle$
&
 $(1,1,1,4,6)$
&
$1/2$
\\
\midrule 78 &
$\KK[{T_1,\ldots,T_5}]/
 \langle   T_1^2T_2^3T_3^7+T_4^3+T_5^2  \rangle$
&
 $(1,1,1,4,6)$
&
$1/2$
\\
\midrule 79 &
$\KK[{T_1,\ldots,T_5}]/
 \langle   T_1^2T_2^5T_3^5+T_4^3+T_5^2  \rangle$
&
 $(1,1,1,4,6)$
&
$1/2$
\\
\midrule 80 &
$\KK[{T_1,\ldots,T_5}]/
 \langle   T_1^3T_2^4T_3^5+T_4^3+T_5^2  \rangle$
&
 $(1,1,1,4,6)$
&
$1/2$
\\
\midrule 81 &
$\KK[{T_1,\ldots,T_5}]/
 \langle   T_1T_2T_3^2+T_4^3+T_5^2  \rangle$
&
 $(1,1,2,2,3)$
&
$27/2$
\\
\midrule 82 &
$\KK[{T_1,\ldots,T_5}]/
 \langle   T_1T_2^3T_3+T_4^3+T_5^2  \rangle$
&
 $(1,1,2,2,3)$
&
$27/2$
\\
\midrule 83 &
$\KK[{T_1,\ldots,T_5}]/
 \langle   T_1^2T_2^2T_3+T_4^3+T_5^2  \rangle$
&
 $(1,1,2,2,3)$
&
$27/2$
\\
\midrule 84 &
$\KK[{T_1,\ldots,T_5}]/
 \langle   T_1T_2T_3^4+T_4^5+T_5^2  \rangle$
&
 $(1,1,2,2,5)$
&
$1/2$
\\
\midrule 85 &
$\KK[{T_1,\ldots,T_5}]/
 \langle   T_1T_2^3T_3^3+T_4^5+T_5^2  \rangle$
&
 $(1,1,2,2,5)$
&
$1/2$
\\
\midrule 86 &
$\KK[{T_1,\ldots,T_5}]/
 \langle   T_1T_2^5T_3^2+T_4^5+T_5^2  \rangle$
&
 $(1,1,2,2,5)$
&
$1/2$
\\
\midrule 87 &
$\KK[{T_1,\ldots,T_5}]/
 \langle   T_1T_2^7T_3+T_4^5+T_5^2  \rangle$
&
 $(1,1,2,2,5)$
&
$1/2$
\\
\midrule 88 &
$\KK[{T_1,\ldots,T_5}]/
 \langle   T_1^2T_2^2T_3^3+T_4^5+T_5^2  \rangle$
&
 $(1,1,2,2,5)$
&
$1/2$
\\
\midrule 89 &
$\KK[{T_1,\ldots,T_5}]/
 \langle   T_1^2T_2^6T_3+T_4^5+T_5^2  \rangle$
&
 $(1,1,2,2,5)$
&
$1/2$
\\
\midrule 90 &
$\KK[{T_1,\ldots,T_5}]/
 \langle   T_1^3T_2^3T_3^2+T_4^5+T_5^2  \rangle$
&
 $(1,1,2,2,5)$
&
$1/2$
\\
\midrule 91 &
$\KK[{T_1,\ldots,T_5}]/
 \langle   T_1^3T_2^5T_3+T_4^5+T_5^2  \rangle$
&
 $(1,1,2,2,5)$
&
$1/2$
\\
\midrule 92 &
$\KK[{T_1,\ldots,T_5}]/
 \langle   T_1^4T_2^4T_3+T_4^5+T_5^2  \rangle$
&
 $(1,1,2,2,5)$
&
$1/2$
\\
\midrule 93 &
$\KK[{T_1,\ldots,T_5}]/
 \langle   T_1T_2T_3^5+T_4^3+T_5^2  \rangle$
&
 $(1,1,2,4,6)$
&
$2$
\\
\midrule 94 &
$\KK[{T_1,\ldots,T_5}]/
 \langle   T_1T_2^3T_3^4+T_4^3+T_5^2  \rangle$
&
 $(1,1,2,4,6)$
&
$2$
\\
\midrule 95 &
$\KK[{T_1,\ldots,T_5}]/
 \langle   T_1T_2^5T_3^3+T_4^3+T_5^2  \rangle$
&
 $(1,1,2,4,6)$
&
$2$
\\
\midrule 96 &
$\KK[{T_1,\ldots,T_5}]/
 \langle   T_1T_2^7T_3^2+T_4^3+T_5^2  \rangle$
&
 $(1,1,2,4,6)$
&
$2$
\\
\midrule 97 &
$\KK[{T_1,\ldots,T_5}]/
 \langle   T_1T_2^9T_3+T_4^3+T_5^2  \rangle$
&
 $(1,1,2,4,6)$
&
$2$
\\
\midrule 98 &
$\KK[{T_1,\ldots,T_5}]/
 \langle   T_1^2T_2^4T_3^3+T_4^3+T_5^2  \rangle$
&
 $(1,1,2,4,6)$
&
$2$
\\
\midrule 99 &
$\KK[{T_1,\ldots,T_5}]/
 \langle   T_1^2T_2^8T_3+T_4^3+T_5^2  \rangle$
&
 $(1,1,2,4,6)$
&
$2$
\\
\midrule 100 &
$\KK[{T_1,\ldots,T_5}]/
 \langle   T_1^3T_2^5T_3^2+T_4^3+T_5^2  \rangle$
&
 $(1,1,2,4,6)$
&
$2$
\\
\midrule 101 &
$\KK[{T_1,\ldots,T_5}]/
 \langle   T_1^3T_2^7T_3+T_4^3+T_5^2  \rangle$
&
 $(1,1,2,4,6)$
&
$2$
\\
\midrule 102 &
$\KK[{T_1,\ldots,T_5}]/
 \langle   T_1^4T_2^6T_3+T_4^3+T_5^2  \rangle$
&
 $(1,1,2,4,6)$
&
$2$
\\
\midrule 103 &
$\KK[{T_1,\ldots,T_5}]/
 \langle   T_1^5T_2^5T_3+T_4^3+T_5^2  \rangle$
&
 $(1,1,2,4,6)$
&
$2$
\\
\midrule 104 &
$\KK[{T_1,\ldots,T_5}]/
 \langle   T_1^2T_2T_3+T_4^3+T_5^2  \rangle$
&
 $(1,2,2,2,3)$
&
$16$
\\
\midrule 105 &
$\KK[{T_1,\ldots,T_5}]/
 \langle   T_1^2T_2T_3^3+T_4^5+T_5^2  \rangle$
&
 $(1,2,2,2,5)$
&
$2$
\\
\midrule 106 &
$\KK[{T_1,\ldots,T_5}]/
 \langle   T_1^4T_2T_3^2+T_4^5+T_5^2  \rangle$
&
 $(1,2,2,2,5)$
&
$2$
\\
\midrule 107 &
$\KK[{T_1,\ldots,T_5}]/
 \langle   T_1^6T_2T_3+T_4^5+T_5^2  \rangle$
&
 $(1,2,2,2,5)$
&
$2$
\\
\bottomrule

\end{longtable}
\end{center}
\end{theorem}

The varieties no. 2,3 and 25, 26 are rational 
degenerations of quasismooth varieties from 
the list in \cite{IaFl}. 
In \cite{CPR} the non-rationality of a general 
(quasismooth) element of the corresponding family 
was proved.

The varieties listed so far might
suggest that we always obtain only 
one relation in the Cox ring.
We discuss now some examples, showing 
that for a Picard index big enough, 
we need in general more than one 
relation, where this refers always 
to a presentation as in 
Theorem~\ref{thm:factrings}~(ii).

\begin{example} 
\label{ex:fanosurf2rel}
A Fano $\KK^*$-surface $X$ with $\Cl(X)=\ZZ$
such that the Cox ring $\mathcal{R}(X)$ needs
two relations.
Consider the $\ZZ$-graded ring 
\begin{eqnarray*}
R & = & 
\KK[T_{01},T_{02},T_{11},T_{21},T_{31}]/\bangle{g_0,g_1},
\end{eqnarray*}
where the degrees of $T_{01},T_{02},T_{11},T_{21},T_{31}$
are $29,1,6,10,15$, respectively, 
and the relations $g_0,g_1$ are given by 
$$
g_0 \ := \ T_{01}T_{02}+T_{11}^5+T_{21}^3,
\qquad\qquad
g_1 \ := \ 
\alpha_{23} T_{11}^5+\alpha_{31} T_{21}^3+\alpha_{12}T_{31}^2
$$
Then $R$ is the Cox ring of a Fano $\KK^*$-surface.
Note that the Picard index is given by 
$
[\Cl(X):\Pic(X)]= \mathrm{lcm}(29,1)=29.
$
\end{example}

\begin{proposition}
\label{prop:fano22rel}
Let $X$ be a non-toric Fano surface with 
an effective $\KK^*$-action such that 
$\Cl(X) \cong \ZZ$ and $[\Cl(X):\Pic(X)] < 29$ 
hold.
Then the Cox ring of $X$ is of the form 
\begin{eqnarray*}
\mathcal{R}(X)
& \cong & 
\KK[T_1, \ldots, T_4]/\bangle{T_1^{l_1}T_2^{l_2} + T_3^{l_3} + T_4^{l_4}}.
\end{eqnarray*}
\end{proposition}

\begin{proof}
The Cox ring $\mathcal{R}(X)$ is as in 
Theorem~\ref{thm:factrings}, and, in the 
notation used there, we have 
$n_0 + \ldots + n_r + m = 2+r$. 
This leaves us with the possibilities
$n_0=m=1$ and $n_0=2$, $m=0$.
In the first case, Proposition~\ref{prop:MoriCox}
tells us that the Picard index of $X$ 
is at least $30$.

So, consider the case $n_0=2$ and $m=0$.
Then, according to Theorem~\ref{thm:factrings},
the Cox ring $\mathcal{R}(X)$ is 
$\KK[T_{01},T_{02},T_1 \ldots, T_r]$
divided by relations
$$
g_{0,1,2}=T_{01}^{l_{01}}T_{02}^{l_{02}} + T_1^{l_1} + T_2^{l_2},
\quad
g_{i,i+1,i+2}=
\alpha_{i+1,i+2}T_i^{l_i} + 
\alpha_{i+2,i}T_{i+1}^{l_{i+1}} + 
\alpha_{i,i+1}T_{i+2}^{l_{i+2}},
$$
where $1 \le i \le r-2$.
We have to show that $r=2$ holds. 
Set $\mu := [\Cl(X):\Pic(X)]$ and 
let $\gamma \in \ZZ$ denote the degree
of the relations. Then we have 
$\gamma = w_il_i$ for $1 \le i \le r$,
where $w_i := \deg \, T_i$.
With $w_{0i} := \deg \, T_{0i}$,
Proposition~\ref{Prop:FanoPicard} gives us
\begin{eqnarray*}
(r-1) \gamma 
& < & 
w_{01} + w_{02} + w_1 + \ldots + w_r.
\end{eqnarray*}
We claim that $w_{01}$ and $w_{02}$ are coprime.
Otherwise they had a common prime divisor $p$.
This $p$ divides $\gamma = l_iw_i$.
Since $l_1,\ldots,l_r$ are pairwise coprime,
$p$ divides at least $r-1$ of the weights 
$w_1,\ldots, w_r$.
This contradicts the Cox ring condition that 
any $r+1$ of the $r+2$ weights generate the class 
group $\ZZ$. 
Thus,  $w_{01}$ and $w_{02}$ are coprime and 
we obtain
$$
\mu \ \ge \  \rm{lcm}(w_{01},w_{02})
    \ = \ w_{01}\cdot w_{02} 
    \ \ge \ w_{01}+w_{02}-1.
$$
Now assume that $r \ge 3$ holds. Then we can conclude
$$ 
2 \gamma 
\ < \ 
w_{01} + w_{02} + w_1 + w_2 + w_3
\ \le \ 
\mu + 1 + 
\gamma \left( \frac{1}{l_1} +  \frac{1}{l_2} +  \frac{1}{l_3}  \right) 
$$
Since the numbers $l_i$ are pairwise coprime, 
we obtain $l_1 \ge 5$, $l_2 \ge 3$ and $l_3 \ge 2$.
Moreover, $l_iw_i = l_jw_j$ implies $l_i \mid w_j$ 
and hence $l_1l_2l_3 \mid \gamma$. Thus, we have 
$\gamma \ge 30$. Plugging this in the above 
inequality gives
$$ 
\mu 
\ \ge \ 
\gamma\left(2- \frac{1}{l_1} -  \frac{1}{l_2} -  \frac{1}{l_3}  \right)-1 
\ = \ 
29.
$$
\end{proof}

The Fano assumption is essential in this result;
if we omit it, then we may even construct locally 
factorial surfaces with a Cox ring that needs more
then one relation.

\begin{example} 
A locally factorial 
$\KK^*$-surface $X$ with $\Cl(X)=\ZZ$
such that the Cox ring $\mathcal{R}(X)$
needs two relations.
Consider the $\ZZ$-graded ring 
\begin{eqnarray*}
R & = & 
\KK[T_{01},T_{02},T_{11},T_{21},T_{31}]/\bangle{g_0,g_1},
\end{eqnarray*}
where the degrees of $T_{01},T_{02},T_{11},T_{21},T_{31}$
are $1,1,6,10,15$, respectively, 
and the relations $g_0,g_1$ are given by 
$$
g_0 \ := \ T_{01}^7T_{02}^{23}+T_{11}^5+T_{21}^3,
\qquad\qquad
g_1 \ := \ 
\alpha_{23} T_{11}^5+\alpha_{31}T_{21}^3+\alpha_{12}T_{31}^2
$$
Then $R$ is the Cox ring of a non Fano $\KK^*$-surface~$X$
of Picard index one, i.e, $X$ is locally factorial.
\end{example}

For non-toric Fano threefolds~$X$ 
with an effective 2-torus action 
$\Cl(X) \cong \ZZ$, 
the classifications~\ref{thm:3fano}
and~\ref{thm:3fano2} 
show that for Picard indices one and 
two we only obtain hypersurfaces as 
Cox rings. 
The following example shows that 
this stops at Picard index three.

\begin{example}
\label{ex:fano32rel}
A Fano threefold $X$ with $\Cl(X)=\ZZ$ 
and a 2-torus action such that 
the Cox ring $\mathcal{R}(X)$ 
needs two relations.
Consider  
\begin{eqnarray*}
R
& = & 
\KK[T_{01},T_{02},T_{11},T_{12},T_{21},T_{31}]/
\bangle{g_0,g_1}
\end{eqnarray*}
where the degrees of $T_{01},T_{02},T_{11},T_{12},T_{21},T_{31}$
are  $1,1,3,3,2,3$, respectively, 
and the relations are given by 
$$
g_0
\ = \ 
T_{01}^5T_{02}+T_{11}T_{12}+T_{21}^3,
\qquad
g_1
\ = \ 
\alpha_{23} T_{11}T_{12}+\alpha_{31}T_{21}^3+\alpha_{12}T_{31}^2. 
$$
Then $R$ is the Cox ring of a Fano threefold 
with a 2-torus action.
Note that the Picard index is given by 
$$
[\Cl(X):\Pic(X)]
\ = \ 
\mathrm{lcm}(1,1,3,3)
\ = \ 3.
$$
\end{example}

Finally, we turn to locally factorial 
Fano fourfolds.
Here we observe more than one relation 
in the Cox ring even
in the locally factorial case.

\begin{theorem}
Let $X$ be a four-dimensional locally factorial non-toric 
Fano variety with an effective three torus action such that 
$\Cl(X)=\ZZ$ holds.
Then its Cox ring is precisely one of the following.
\begin{center}
\begin{longtable}[htbp]{llll}
\toprule
No.
&
$\mathcal{R}(X)$
&
 $(w_1,\ldots,w_6)$
& 
$(-K_X)^4$
\\
\midrule  1 &
$\KK[{T_1,\ldots,T_6}]/
 \langle   T_1T_2^5+T_3^3+T_4^2  \rangle$
&
 $(1,1,2,3,1,1)$
&
$81$
\\
\midrule  2 &
$\KK[{T_1,\ldots,T_6}]/
 \langle   T_1T_2^9+T_3^2+T_4^5  \rangle$
&
 $(1,1,2,5,1,1)$
&
$1$
\\
\midrule  3 &
$\KK[{T_1,\ldots,T_6}]/
 \langle   T_1^3T_2^7+T_3^2+T_4^5  \rangle$
&
 $(1,1,2,5,1,1)$
&
$1$
\\
\midrule  4 &
$\KK[{T_1,\ldots,T_6}]/
 \langle   T_1T_2T_3^4+T_4^3+T_5^2  \rangle$
&
 $(1,1,1,2,3,1)$
&
$81$
\\
\midrule  5 & 
$\KK[{T_1,\ldots,T_6}]/
 \langle   T_1T_2^2T_3^3+T_4^3+T_5^2  \rangle$
&
 $(1,1,1,2,3,1)$
&
$81$
\\
\midrule  6 &
$\KK[{T_1,\ldots,T_6}]/
 \langle   T_1T_2T_3^8+T_4^5+T_5^2  \rangle$
&
 $(1,1,1,2,5,1)$
&
$1$
\\
\midrule  7 &
$\KK[{T_1,\ldots,T_6}]/
 \langle   T_1T_2^2T_3^7+T_4^5+T_5^2  \rangle$
&
 $(1,1,1,2,5,1)$
&
$1$
\\
\midrule  8 &
$\KK[{T_1,\ldots,T_6}]/
 \langle   T_1T_2^3T_3^6+T_4^5+T_5^2  \rangle$
&
 $(1,1,1,2,5,1)$
&
$1$
\\
\midrule  9 &
$\KK[{T_1,\ldots,T_6}]/
 \langle   T_1T_2^4T_3^5+T_4^5+T_5^2  \rangle$
&
 $(1,1,1,2,5,1)$
&
$1$
\\
\midrule  10 &
$\KK[{T_1,\ldots,T_6}]/
 \langle   T_1^2T_2^3T_3^5+T_4^5+T_5^2  \rangle$
&
 $(1,1,1,2,5,1)$
&
$1$
\\
\midrule  11 &
$\KK[{T_1,\ldots,T_6}]/
 \langle   T_1^3T_2^3T_3^4+T_4^5+T_5^2  \rangle$
&
 $(1,1,1,2,5,1)$
&
$1$
\\
\midrule  12 &
$\KK[{T_1,\ldots,T_6}]/
 \langle   T_1T_2+T_3T_4+T_5^2  \rangle$
&
 $(1,1,1,1,1,1)$
&
$512$
\\
\midrule  13 &
$\KK[{T_1,\ldots,T_6}]/
 \langle   T_1T_2^2+T_3T_4^2+T_5^3  \rangle$
&
 $(1,1,1,1,1,1)$
&
$243$
\\
\midrule  14 &
$\KK[{T_1,\ldots,T_6}]/
 \langle   T_1T_2^3+T_3T_4^3+T_5^4  \rangle$
&
 $(1,1,1,1,1,1)$
&
$64$
\\
\midrule  15 &
$\KK[{T_1,\ldots,T_6}]/
 \langle   T_1T_2^4+T_3T_4^4+T_5^5  \rangle$
&
 $(1,1,1,1,1,1)$
&
$5$
\\
\midrule  16 &
$\KK[{T_1,\ldots,T_6}]/
 \langle   T_1T_2^4+T_3^2T_4^3+T_5^5  \rangle$
&
 $(1,1,1,1,1,1)$
&
$5$
\\
\midrule  17 &
$\KK[{T_1,\ldots,T_6}]/
 \langle   T_1^2T_2^3+T_3^2T_4^3+T_5^5  \rangle$
&
 $(1,1,1,1,1,1)$
&
$5$
\\
\midrule  18 &
$\KK[{T_1,\ldots,T_6}]/
 \langle   T_1T_2^3+T_3T_4^3+T_5^2  \rangle$
&
 $(1,1,1,1,2,1)$
&
$162$
\\
\midrule  19 &
$\KK[{T_1,\ldots,T_6}]/
 \langle   T_1T_2^5+T_3T_4^5+T_5^3  \rangle$
&
 $(1,1,1,1,2,1)$
&
$3$
\\
\midrule  20 &
$\KK[{T_1,\ldots,T_6}]/
 \langle   T_1T_2^5+T_3^2T_4^4+T_5^3  \rangle$
&
 $(1,1,1,1,2,1)$
&
$3$
\\
\midrule  21 &
$\KK[{T_1,\ldots,T_6}]/
 \langle   T_1T_2^5+T_3T_4^5+T_5^2  \rangle$
&
 $(1,1,1,1,3,1)$
&
$32$
\\
\midrule  22 &
$\KK[{T_1,\ldots,T_6}]/
 \langle   T_1T_2^5+T_3^3T_4^3+T_5^2  \rangle$
&
 $(1,1,1,1,3,1)$
&
$32$
\\
\midrule  23 &
$\KK[{T_1,\ldots,T_6}]/
 \langle   T_1T_2^7+T_3T_4^7+T_5^2  \rangle$
&
 $(1,1,1,1,4,1)$
&
$2$
\\
\midrule  24 &
$\KK[{T_1,\ldots,T_6}]/
 \langle   T_1T_2^7+T_3^3T_4^5+T_5^2  \rangle$
&
 $(1,1,1,1,4,1)$
&
$2$
\\
\midrule  25 &
$\KK[{T_1,\ldots,T_6}]/
 \langle   T_1^3T_2^5+T_3^3T_4^5+T_5^2  \rangle$
&
 $(1,1,1,1,4,1)$
&
$2$
\\
\midrule  26 &
$\KK[{T_1,\ldots,T_6}]/
 \langle   T_1T_2T_3T_4^3+T_5^3+T_6^2  \rangle$
&
 $(1,1,1,1,2,3)$
&
$81$
\\
\midrule  27 &
$\KK[{T_1,\ldots,T_6}]/
 \langle   T_1T_2T_3^2T_4^2+T_5^3+T_6^2  \rangle$
&
 $(1,1,1,1,2,3)$
&
$81$
\\
\midrule  28 &
$\KK[{T_1,\ldots,T_6}]/
 \langle   T_1T_2T_3T_4^7+T_5^5+T_6^2  \rangle$
&
 $(1,1,1,1,2,5)$
&
$1$
\\
\midrule  29 &
$\KK[{T_1,\ldots,T_6}]/
 \langle   T_1T_2T_3^2T_4^6+T_5^5+T_6^2  \rangle$
&
 $(1,1,1,1,2,5)$
&
$1$
\\
\midrule  30 &
$\KK[{T_1,\ldots,T_6}]/
 \langle   T_1T_2T_3^3T_4^5+T_5^5+T_6^2  \rangle$
&
 $(1,1,1,1,2,5)$
&
$1$
\\
\midrule  31 &
$\KK[{T_1,\ldots,T_6}]/
 \langle   T_1T_2T_3^4T_4^4+T_5^5+T_6^2  \rangle$
&
 $(1,1,1,1,2,5)$
&
$1$
\\
\midrule  32 &
$\KK[{T_1,\ldots,T_6}]/
 \langle   T_1T_2^2T_3^2T_4^5+T_5^5+T_6^2  \rangle$
&
 $(1,1,1,1,2,5)$
&
$1$
\\
\midrule  33 &
$\KK[{T_1,\ldots,T_6}]/
 \langle   T_1T_2^2T_3^3T_4^4+T_5^5+T_6^2  \rangle$
&
 $(1,1,1,1,2,5)$
&
$1$
\\
\midrule  34 &
$\KK[{T_1,\ldots,T_6}]/
 \langle   T_1T_2^3T_3^3T_4^3+T_5^5+T_6^2  \rangle$
&
 $(1,1,1,1,2,5)$
&
$1$
\\
\midrule  35 &
$\KK[{T_1,\ldots,T_6}]/
 \langle   T_1^2T_2^2T_3^3T_4^3+T_5^5+T_6^2  \rangle$
&
 $(1,1,1,1,2,5)$
&
$1$
\\
\midrule  36 &
$\KK[{T_1,\ldots,T_6}]/
 \langle   T_1T_2T_3+T_4T_5^2+T_6^3  \rangle$
&
 $(1,1,1,1,1,1)$
&
$243$
\\
\midrule  37 &
$\KK[{T_1,\ldots,T_6}]/
 \langle   T_1T_2T_3^2+T_4T_5^3+T_6^4  \rangle$
&
 $(1,1,1,1,1,1)$
&
$64$
\\
\midrule  38 &
$\KK[{T_1,\ldots,T_6}]/
 \langle   T_1T_2T_3^3+T_4T_5^4+T_6^5  \rangle$
&
 $(1,1,1,1,1,1)$
&
$5$
\\
\midrule  39 &
$\KK[{T_1,\ldots,T_6}]/
 \langle   T_1T_2T_3^3+T_4^2T_5^3+T_6^5  \rangle$
&
 $(1,1,1,1,1,1)$
&
$5$
\\
\midrule  40 &
$\KK[{T_1,\ldots,T_6}]/
 \langle   T_1T_2^2T_3^2+T_4T_5^4+T_6^5  \rangle$
&
 $(1,1,1,1,1,1)$
&
$5$
\\
\midrule  41 &
$\KK[{T_1,\ldots,T_6}]/
 \langle   T_1T_2^2T_3^2+T_4^2T_5^3+T_6^5  \rangle$
&
 $(1,1,1,1,1,1)$
&
$5$
\\
\midrule  42 &
$\KK[{T_1,\ldots,T_6}]/
 \langle   T_1T_2T_3^2+T_4T_5^3+T_6^2  \rangle$
&
 $(1,1,1,1,1,2)$
&
$162$
\\
\midrule  43 &
$\KK[{T_1,\ldots,T_6}]/
 \langle   T_1T_2T_3^4+T_4T_5^5+T_6^3  \rangle$
&
 $(1,1,1,1,1,2)$
&
$3$
\\
\midrule  44 &
$\KK[{T_1,\ldots,T_6}]/
 \langle   T_1T_2T_3^4+T_4^2T_5^4+T_6^3  \rangle$
&
 $(1,1,1,1,1,2)$
&
$3$
\\
\midrule  45 &
$\KK[{T_1,\ldots,T_6}]/
 \langle   T_1T_2^2T_3^3+T_4T_5^5+T_6^3  \rangle$
&
 $(1,1,1,1,1,2)$
&
$3$
\\
\midrule  46 &
$\KK[{T_1,\ldots,T_6}]/
 \langle   T_1T_2^2T_3^3+T_4^2T_5^4+T_6^3  \rangle$
&
 $(1,1,1,1,1,2)$
&
$3$
\\
\midrule  47 &
$\KK[{T_1,\ldots,T_6}]/
 \langle   T_1^2T_2^2T_3^2+T_4T_5^5+T_6^3  \rangle$
&
 $(1,1,1,1,1,2)$
&
$3$
\\
\midrule  48 &
$\KK[{T_1,\ldots,T_6}]/
 \langle   T_1T_2^2T_3^3+T_4^3T_5^3+T_6^2  \rangle$
&
 $(1,1,1,1,1,3)$
&
$32$
\\
\midrule  49 &
$\KK[{T_1,\ldots,T_6}]/
 \langle   T_1T_2^2T_3^3+T_4T_5^5+T_6^2  \rangle$
&
 $(1,1,1,1,1,3)$
&
$32$
\\
\midrule  50 &
$\KK[{T_1,\ldots,T_6}]/
 \langle   T_1T_2T_3^4+T_4^3T_5^3+T_6^2  \rangle$
&
 $(1,1,1,1,1,3)$
&
$32$
\\
\midrule  51 &
$\KK[{T_1,\ldots,T_6}]/
 \langle   T_1T_2T_3^4+T_4T_5^5+T_6^2  \rangle$
&
 $(1,1,1,1,1,3)$
&
$32$
\\
\midrule  52 &
$\KK[{T_1,\ldots,T_6}]/
 \langle   T_1T_2T_3^6+T_4T_5^7+T_6^2  \rangle$
&
 $(1,1,1,1,1,4)$
&
$2$
\\
\midrule  53 &
$\KK[{T_1,\ldots,T_6}]/
 \langle   T_1T_2T_3^6+T_4^3T_5^5+T_6^2  \rangle$
&
 $(1,1,1,1,1,4)$
&
$2$
\\
\midrule  54 &
$\KK[{T_1,\ldots,T_6}]/
 \langle   T_1T_2^2T_3^5+T_4T_5^7+T_6^2  \rangle$
&
 $(1,1,1,1,1,4)$
&
$2$
\\
\midrule  55 &
$\KK[{T_1,\ldots,T_6}]/
 \langle   T_1T_2^2T_3^5+T_4^3T_5^5+T_6^2  \rangle$
&
 $(1,1,1,1,1,4)$
&
$2$
\\
\midrule  56 &
$\KK[{T_1,\ldots,T_6}]/
 \langle   T_1T_2^3T_3^4+T_4T_5^7+T_6^2  \rangle$
&
 $(1,1,1,1,1,4)$
&
$2$
\\
\midrule  57 &
$\KK[{T_1,\ldots,T_6}]/
 \langle   T_1T_2^3T_3^4+T_4^3T_5^5+T_6^2  \rangle$
&
 $(1,1,1,1,1,4)$
&
$2$
\\
\midrule  58 &
$\KK[{T_1,\ldots,T_6}]/
 \langle   T_1^2T_2^3T_3^3+T_4T_5^7+T_6^2  \rangle$
&
 $(1,1,1,1,1,4)$
&
$2$
\\
\midrule  59 &
$\KK[{T_1,\ldots,T_6}]/
 \langle   T_1^2T_2^3T_3^3+T_4^3T_5^5+T_6^2  \rangle$
&
 $(1,1,1,1,1,4)$
&
$2$
\\
\midrule  60 &
$\KK[{T_1,\ldots,T_6}]/
 \langle   T_1T_2+T_3T_4+T_5T_6  \rangle$
&
 $(1,1,1,1,1,1)$
&
$512$
\\
\midrule  61 &
$\KK[{T_1,\ldots,T_6}]/
 \langle   T_1T_2^2+T_3T_4^2+T_5T_6^2  \rangle$
&
 $(1,1,1,1,1,1)$
&
$243$
\\
\midrule  62 & 
$\KK[{T_1,\ldots,T_6}]/
 \langle   T_1T_2^3+T_3T_4^3+T_5T_6^3  \rangle$
&
 $(1,1,1,1,1,1)$
&
$64$
\\
\midrule  63 & 
$\KK[{T_1,\ldots,T_6}]/
 \langle   T_1T_2^3+T_3T_4^3+T_5^2T_6^2  \rangle$
&
 $(1,1,1,1,1,1)$
&
$64$
\\
\midrule  64 &
$\KK[{T_1,\ldots,T_6}]/
 \langle   T_1T_2^4+T_3T_4^4+T_5T_6^4  \rangle$
&
 $(1,1,1,1,1,1)$
&
$5$
\\
\midrule  65 & 
$\KK[{T_1,\ldots,T_6}]/
 \langle   T_1T_2^4+T_3T_4^4+T_5^2T_6^3  \rangle$
&
 $(1,1,1,1,1,1)$
&
$5$
\\
\midrule  66 &
$\KK[{T_1,\ldots,T_6}]/
 \langle   T_1T_2^4+T_3^2T_4^3+T_5^2T_6^3  \rangle$
&
 $(1,1,1,1,1,1)$
&
$5$
\\
\midrule  67 &
$\KK[{T_1,\ldots,T_6}]/
 \langle   T_1^2T_2^3+T_3^2T_4^3+T_5^2T_6^3  \rangle$
&
 $(1,1,1,1,1,1)$
&
$5$
\\
\midrule 68 
&
$\KK[T_1,\ldots,T_7] /
\left\langle
\begin{smallmatrix}
T_1T_2 + T_3T_4 + T_5T_6,
\\
\alpha T_3T_4 +  T_5T_6 + T_7^2
\end{smallmatrix}
\right\rangle
$
&
 $(1,1,1,1,1,1,1)$
&
$324$
\\
\midrule 69 
&
$\KK[T_1,\ldots,T_7] /
\left\langle
\begin{smallmatrix}
T_1T_2^2 + T_3T_4^2 + T_5T_6^2,
\\
\alpha T_3T_4^2 +  T_5T_6^2 + T_7^3
\end{smallmatrix}
\right\rangle
$
&
 $(1,1,1,1,1,1,1)$
&
$9$
\\
\bottomrule

\end{longtable}
\end{center}
where in the last two rows of the table the parameter 
$\alpha$ can be any element from $\KK^* \setminus \{1\}$.
\end{theorem}

By the result of \cite{Pu2}, the singular 
quintics of this list are rational degenerations
of smooth non-rational Fano fourfolds.

\section{Geometry of the locally factorial threefolds}
\label{sec:geom3folds}

In this section, we take a closer look 
at the (factorial) singularities of the
Fano varieties~$X$ listed in Theorem~\ref{thm:3fano}.
Recall that the discrepancies of a resolution 
$\varphi \colon \t{X} \to X$ of a
singularity are the coefficients 
of $K_{\t{X}} - \varphi^* K_X$, where
$K_X$ and $K_{\t{X}}$ are canonical divisors
such that $K_{\t{X}} - \varphi^* K_X$
is supported on the exceptional locus 
of $\varphi$.
A resolution is called crepant, if its
discrepancies vanish and a singularity 
is called canonical (terminal), 
if it admits a resolution with 
nonnegative (positive) discrepancies.
By a relative minimal model we mean a 
projective morphism $\t{X} \to X$ such that 
$\t{X}$ has at most terminal singularities
and its relative canonical divisor is 
relatively nef.

\begin{theorem}
\label{prop:3foldsing}
For the nine 3-dimensional Fano varieties 
listed in Theorem~\ref{thm:3fano}, we have 
the following statements.
\begin{enumerate}
\item 
No.~4 is a smooth quadric in $\PP^4$.
\item 
Nos.~1,3,5,7 and 9 are singular with only 
canonical singularities and all admit 
a crepant resolution.
\item 
Nos.~6 and 8 are  singular with  
non-canonical singularities but
admit a smooth relative minimal model.
\item 
No.~2 is singular with only canonical singularities,
one of them of type $\mathbf{cA_1}$,
and admits only a singular relative minimal model.
   \end{enumerate}
The Cox ring of the relative minimal model $\widetilde{X}$ 
as well as the the Fano degree of $X$ itself are 
given in the following table. 
\begin{center}
\begin{longtable}[htbp]{llc}
\toprule
No.
& 
\hspace{3.5cm}$\mathcal{R}(\widetilde{X})$ 
&
$(-K_X)^3$
\\
\midrule
1
&
$\KK[T_1,\ldots,T_{14}]/(T_1T_2T_3^2T_4^3T_5^4T_6^5+T_7^3T_8^2T_9+T_{10}^2T_{11}\rangle$
&
$8$
\\
\cmidrule{1-3}
2
&
$\KK[T_1,\ldots,T_9]/\langle T_1T_2T_3^2T_4^4+T_5T_6^2T_7^3+T_8^2 \rangle$
&
$8$
\\
\cmidrule{1-3}
3
&
$\KK[T_1,\ldots,T_8]/ \langle T_1T_2^2T_3^3+T_4T_5^3+T_6T_7^2\rangle$
&
$8$
\\
\cmidrule{1-3}
4
&
$\KK[T_1, \ldots, T_5]/\bangle{T_1T_2 + T_3T_4 + T_5^2}$
&
$54$
\\
\cmidrule{1-3}
5
&
$\KK[T_1,\ldots,T_6]/\langle T_1T_2^2+T_3T_4^2+T_5^3T_6\rangle$
&
$24$
\\
\cmidrule{1-3}
6
&
$\KK[T_1,\ldots,T_6]/ \langle T_1T_2^3+T_3T_4^3+T_5^4T_6 \rangle$
&
$4$
\\
\cmidrule{1-3}
7
&
$\KK[T_1, \ldots, T_7]/\bangle{T_1T_2^3 + T_3T_4^3 + T_5^2T_6}$
&
$16$
\\
\cmidrule{1-3}
8
&
$\KK[T_1, \ldots, T_7]/\bangle{T_1T_2^5 + T_3T_4^5 + T_5^2T_6}$
&
$2$
\\
\cmidrule{1-3}
9
&
$\displaystyle \KK[T_1,\ldots,T_{46}]/
\left\langle
\begin{smallmatrix}
   T_1T_2T_3T_4^2T_5^2T_6^3T_7^3T_8^4T_9^4T_{10}^5
\;+\; \\ \;+\; T_{11} \cdots T_{18} T_{19}^2\cdots T_{24}^2 T_{25}^3 T_{26}^3
\;+\; T_{27}\cdots T_{32} T_{33}^2
\end{smallmatrix}
\right\rangle$
&
$2$\\
\bottomrule 
\end{longtable}
\end{center}
\end{theorem}

For the proof, it is convenient to work 
in the language of polyhedral divisors introduced 
in~\cite{MR2207875} and~\cite{divfans}.
As we are interested in rational varieties
with a complexity one torus action,
we only have to consider polyhedral divisors
on the projective line $Y = \PP^1$.
This considerably simplifies the general 
definitions and allows us to give a 
short summary.
In the sequel, $N \cong \ZZ^n$ denotes a lattice 
and $M = \Hom(N,\ZZ)$ its dual. 
For the associated rational vector spaces we write
$N_\QQ$ and $M_\QQ$.
A {\em polyhedral divisor\/} on the projective line 
$Y := \PP^1$ is a formal sum
\begin{eqnarray*}
\D 
& = & 
\sum_{y \in Y} \D_y \cdot y,
\end{eqnarray*}
where the coefficients $\D_y \subseteq N_{\QQ}$ 
are (possibly empty) convex polyhedra 
all sharing the same tail (i.e.~recession) 
cone $\D_Y = \sigma \subseteq N_\QQ$,
and only finitely many $\D_y$ differ from
$\sigma$. 
The {\em locus\/} of $\D$ is the open subset 
$Y(\D) \subseteq Y$ obtained by removing all 
points $y \subseteq Y$ with $\D_y = \emptyset$.
For every $u \in \sigma^{\vee} \cap M$ we have the 
{\em evaluation\/}
\begin{eqnarray*}
\D(u)
& := & 
\sum_{y \in Y} \min_{v \in \D_y} \bangle{u ,v} \mal y,
\end{eqnarray*}
which is a usual rational divisor on $Y(\D)$.
We call the polyhedral divisor $\D$ on $Y$ 
{\em proper\/} if $\deg \, \D \subsetneq \sigma$
holds, where the {\em polyhedral degree\/} 
is defined by 
\begin{eqnarray*}
\deg \, \D 
& := & 
\sum_{y \in Y} \D_y.
\end{eqnarray*}
Every proper polyhedral divisor 
$\D$ on $Y$ defines a
normal affine variety $X(\D)$ 
of dimension $\rk(N)+1$
coming with an effective action of
the torus $T = \Spec \, \KK[M]$:
set $X(\D) := \Spec \, A(\D)$, where
$$
A(\D)
\ := \ 
\bigoplus_{u \in \sigma^\vee \cap M} \Gamma(Y(\D),\mathcal{O}(\D(u)))
\ \subseteq \ 
\bigoplus_{u \in M} \KK(Y) \cdot \chi^u.
$$

A {\em divisorial fan\/}, 
is a finite set $\fan$ of polyhedral divisors 
$\D$ on $Y$, all having their 
polyhedral coefficients $\D_y$ in 
the same $N_{\QQ}$ 
and fulfilling certain compatibility 
conditions, see~\cite{divfans}.
In particular, for every point $y \in Y$,
the {\em slice\/}
\begin{eqnarray*}
\fan_y  
& := & 
\left\{\D_y; \; \D \in \fan \right\}
\end{eqnarray*}
must be a polyhedral subdivision.
The {\em tail fan\/}  is 
the set $\fan_Y$ of the tail cones $\D_Y$
of the $\D \in \fan$; it is a fan in the 
usual sense. 
Given a divisorial fan $\fan$,
the affine varieties $X(\D)$, where 
$\D \in \fan$, glue equivariantly
together to a normal variety $X(\fan)$,
and we obtain every rational normal 
variety with a 
complexity one torus action this way.

Smoothness of $X = X(\fan)$ 
is checked locally.
For a proper polyhedral divisor $\D$ on $Y$,
we infer the following from~\cite[Theorem~3.3]{tfano}.
If $Y(\D)$ is affine, 
then $X(\D)$ is smooth
if and only if 
$\cone(\{1\} \times \D_y) \subseteq \QQ \times N_{\QQ}$,
the convex, polyhedral cone 
generated  by $\{1\} \times \D_y$,
is regular for every $y \in Y(\D)$.
If $Y(\D) = Y$ holds, then
$X(\D)$ is smooth
if and only if there are $y,z \in Y$ 
such that $\D = \D_y y + \D_z z$ 
holds and 
$\cone(\{1\} \times \D_y)
+
\cone(\{-1\} \times \D_z)$
is a regular cone in  
$\QQ \times N_{\QQ}$.
Similarly to toric geometry, singularities 
of $X(\D)$ are resolved by means of 
subdividing~$\D$. 
This means to consider divisorial fans $\Xi$ 
such that for any $y \in Y$, the slice $\Xi_y$ 
is a subdivision of $\D_y$. 
Such a $\fan$ defines a dominant morphism 
$X(\Xi) \rightarrow X(\D)$
and a slight generalization 
of~\cite[Thm.~7.5.]{divfans} 
yields that this morphism is proper. 

\goodbreak

\begin{proposition}
\label{sec:prop-divfans}
The 3-dimensional Fano varieties 
No. 1-8 listed in Theorem~\ref{thm:3fano} 
and their relative minimal models
arise from divisorial fans having the 
following slices and tail cones.

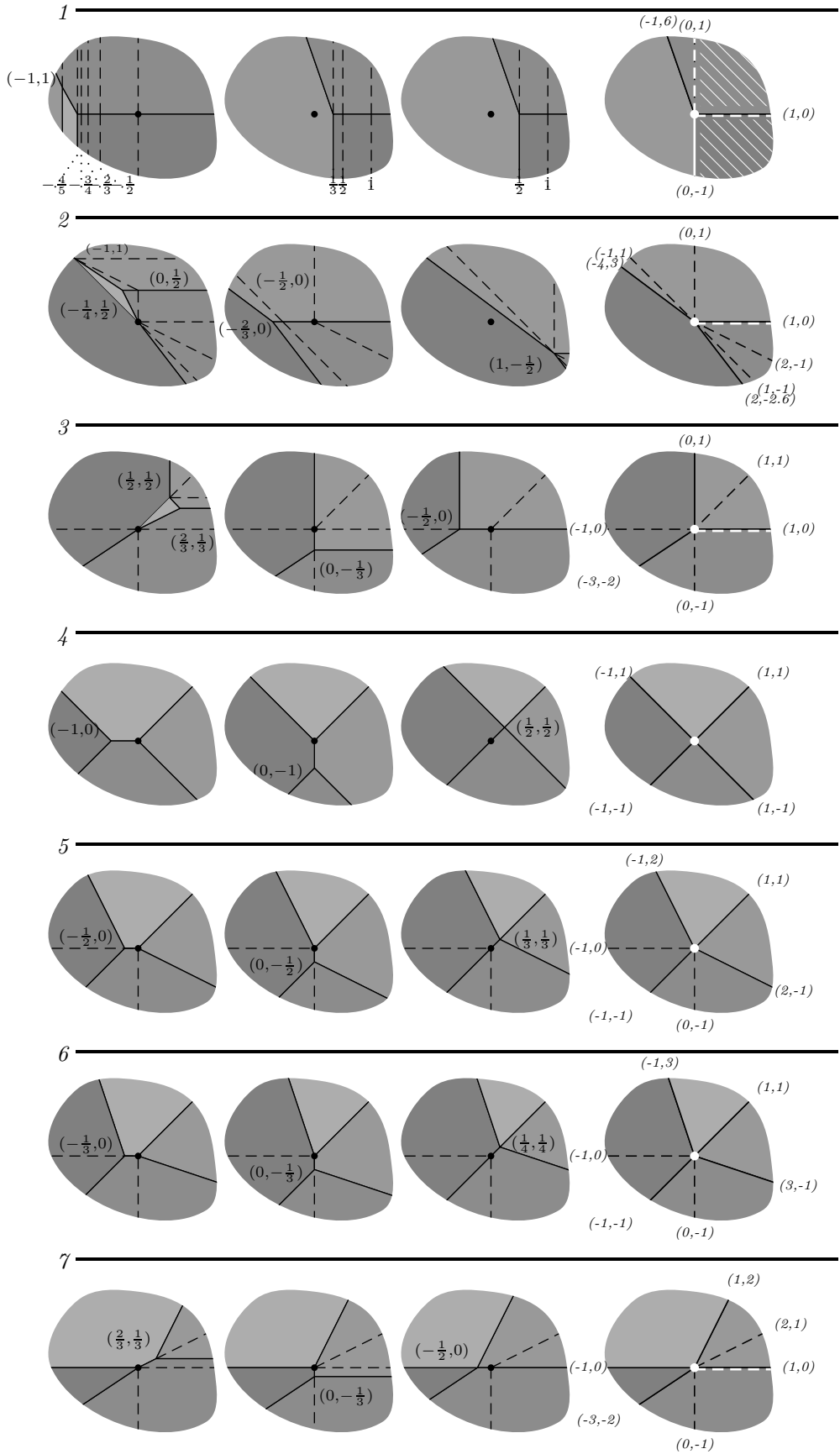
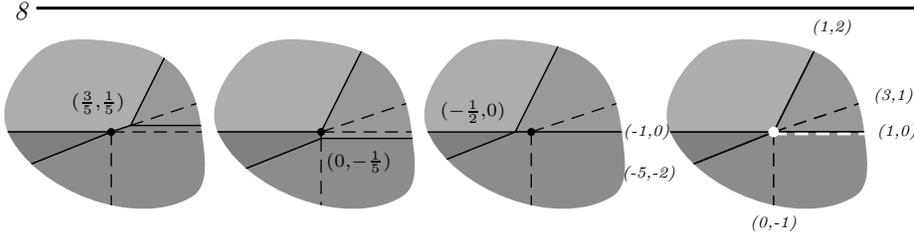
\myrule{1}{
\threefoldG
}

\myrule{2}{
\threefoldE
}

\myrule{3}{
\threefoldF
}

\myrule{4}{
\threefoldA
}

\myrule{5}{
\threefoldB
}

\myrule{6}{
\threefoldC
}

\myrule{7}{
\threefoldH
}

\myrule{8}{
\threefoldD
}
\noindent
\end{proposition}

The above table should be interpreted as follows. 
The first three pictures in each row are the slices at 
$0$, $1$ and $\infty$ and the last one is the tail fan. 
The divisorial fan of the fano variety itself is 
given by the solid polyhedra in the pictures. 
Here, all polyhedra of the same gray scale 
belong to the same polyhedral divisor.  
The subdivisions for the relative minimal models 
are sketched with dashed lines. 
In general, polyhedra with the same tail cone belong 
all to a unique polyhedral divisor with complete locus. 
For the white cones inside the tail fan we have another rule: 
for every polyhedron $\Delta \in \fan_y$ with the given 
white cone as its tail there is a polyhedral divisor 
$\Delta \cdot y + \emptyset \cdot z \in \fan$, 
with $z \in \{0,1,\infty\} \setminus \{y\}$. 
Here, different choices of $z$ lead to isomorphic 
varieties, only the affine covering given by the $X(\D)$ changes.

In order to prove Theorem~\ref{prop:3foldsing},
we also have to understand invariant divisors 
on $X = X(\fan)$ in terms of $\fan$, 
see~\cite[Prop.~4.11 and~4.12]{HaSu}
for details.
A first type of invariant prime divisors,
is in bijection $D_{y,v} \leftrightarrow (y,v)$
with the vertices $(y,v)$, where $y \in Y$ 
and $v \in \fan_y$ is of dimension zero. 
The order of the generic isotropy group
along $D_{y,v}$ equals the minimal positive integer
$\mu(v)$ with $\mu(v) v \in N$.
A second type of invariant prime divisors,
is in
$D_{\varrho} \leftrightarrow \varrho$
with the extremal rays $\varrho \in \fan_Y$,
where a ray $\varrho \in \fan_Y$ 
is called extremal if there is a 
$\D \in \fan$ such that 
$\varrho \subseteq \D_Y$ 
and $\deg \, \D \cap \varrho = \emptyset$
holds. 
The set of extremal rays is denoted by $\fan_Y^\times$.
The divisor of a semi-invariant function $f \cdot \chi^u \in \KK(X)$ 
is then given by 
\begin{eqnarray*}
\div(f \cdot \chi^u)
& = & 
- \sum_{y \in Y} \sum_{v \in \Xi_y^{(0)}} 
\mu(v) \cdot (\langle v, u \rangle + \ord_y f) \cdot D_{y,v} 
\ - \ 
\sum_{\varrho \in \fan_Y^\times} \langle n_\varrho, u \rangle \cdot D_\varrho.
\end{eqnarray*}
Next we describe the canonical divisor.
Choose a point $y_0 \in Y$ such that 
$\fan_{y_0} = \fan_Y$ holds.
Then a canonical divisor on $X = X(\fan)$ is given by
\begin{eqnarray*}
K_X 
& = & 
(s - 2) \cdot y_0 
\ - \ 
\sum_{\fan_y \ne \fan_Y} \sum_{v \in \fan_i^{(0)}}  D_{y,v} 
\ - \ 
\sum_{\varrho \in \fan_Y^\times} E_\varrho.
\end{eqnarray*}

\begin{proposition}
\label{prop:discrepancies}
Let $\D$ be a proper polyhedral divisor 
with $Y(\D) = \PP_1$,
let $\fan$ be a refinement of $\D$
and denote by $y_1, \ldots, y_s \in Y$ 
the points with $\fan_{y_i} \ne \fan_Y$.
Then the associated morphism 
$\varphi \colon X(\fan) \to X(\D)$
satisfies the following.
\begin{enumerate}
\item
The prime divisors in the exceptional 
locus of $\varphi$ are the divisors
$D_{y_i,v}$ and $D_{\varrho}$ 
corresponding to 
$v \in \fan_{y_i}^{(0)} \setminus \D_{y_i}^{(0)}$ 
and 
$\varrho \in \fan_Y^\times \setminus \D^\times$
respectively.
\item 
Then the discrepancies along 
the prime divisors 
$D_{y_i,v}$ and $D_{\varrho}$ 
of~(i) are computed as
\[
d_{y_i,v} 
\ = \ 
-\mu(v)\cdot (\langle v, u' \rangle + \alpha_y) - 1,
\qquad\qquad
 d_{\varrho} 
\ = \ 
-\langle v_\varrho ,  u' \rangle - 1,
\]
where the numbers $\alpha_i$ are determined by
\begin{eqnarray*}
\begin{pmatrix}
  -1 & -1 & \ldots &  -1& 0 \\
  \hline
  \mu(v_{1}^1) & 0 & \ldots & 0 & \mu(v_{1}^1) v_{1}^1 \\
  \vdots & \vdots  & & \vdots &\vdots \\
  \mu(v_{1}^{r_1})& 0 & \ldots & 0 & \mu(v_{1}^{r_1}) v_{1}^{r_1} \\
         &         & \ddots & & \\
  0 & 0 & \ldots &  \mu(v_{s}^1)    &  \mu(v_{s}^1) v_{s}^{1} \\
  \vdots & \vdots  & & \vdots &\vdots \\
  0 & 0 & \ldots &  \mu(v_{s}^{r_s}) & \mu(v_{s}^{r_s}) v_{s}^{r_s} \\
  \hline
  0 & 0 & \ldots & 0 & n_{\varrho_1} \\
  \vdots & \vdots  &  & \vdots &\vdots \\
  0 & 0 & \ldots & 0 & n_{\varrho_{r}}   
\end{pmatrix}
\ \cdot \
\begin{pmatrix}
  \alpha_{y_1}\\
  \vdots\\
  \alpha_{y_s}\\
  u
\end{pmatrix}
& = & 
\begin{pmatrix}
2-s\\
1 \\
\vdots \\
1\\
\hline
1\\
\vdots\\
1\\
\end{pmatrix}
\end{eqnarray*}
\end{enumerate}
\end{proposition}

\begin{proof}
The first claim is obvious by the characterization 
of invariant prime divisors. 
For the second claim note that by~\cite[Theorem~3.1]{tidiv}
every Cartier divisor on $X(\D)$ is principal. 
Hence, we may assume 
$$
\ell\cdot K_X
\ = \ 
\div(f \cdot \chi^{u}),
\qquad\qquad
\div(f) 
\ = \
\sum_y \alpha_y \cdot y. 
$$
Then our formul{\ae} for $\div(f \cdot \chi^{u})$ 
and $K_X$ provide a row for every vertex 
$v_{i}^j \in \fan_{y_i}$, $i=0,\ldots,s$, 
and for every extremal ray $\varrho_i \in \fan^\times$,
and ${\ell}^{-1}(\alpha,u)$
is the (unique) solution of the above system.
\end{proof}

Note, that in the above Proposition, 
the variety $X(\D)$ is $\QQ$-Gorenstein if and only 
if the linear system of equations has a solution.

\begin{proof}[Proof of Theorem~\ref{prop:3foldsing} 
and Proposition~\ref{sec:prop-divfans}]
We exemplarily discuss variety number
eight.
Recall that its Cox ring is given as
\begin{eqnarray*}
\mathcal{R}(X)
& = & 
\mathbb{K}[T_1,\ldots,T_5]/(T_1T_2^5+T_3T_4^5+T_5^2)
\end{eqnarray*}
with the degrees $1,1,1,1,3$. 
In particular, $X$ is a hypersurface 
of degree $6$ in $\PP(1,1,1,1,3)$,
and the self-intersection of the anti-canonical 
divisor can be calculated as
$$
(-K_X^3) 
\ = \ 
6 \cdot \frac{(1+1+1+1+3-6)^3}{1\cdot 1\cdot 1\cdot 1\cdot 3} 
\ = \ 
2.
$$

The embedding $X \subseteq \PP(1,1,1,1,3)$ is equivariant,
and thus we can use the technique described 
in~\cite[Sec.~11]{MR2207875} to calculate a divisorial 
fan $\fan$ for $X$.
The result is the following divisorial fan;
we draw its slices and indicate the polyhedral 
divisors with affine locus by colouring their
tail cones $\D_Y \in \fan_Y$ white: 

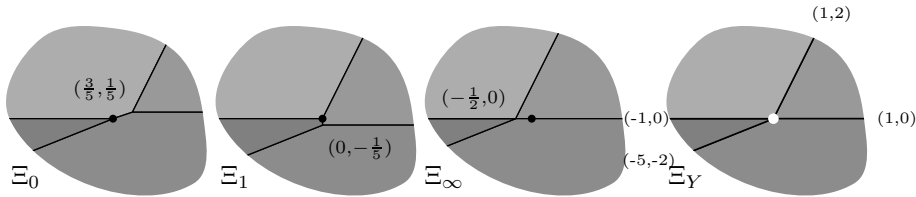
\threefoldDplain

\noindent
One may also use~\cite[Cor.~4.9.]{HaSu} to verify 
that  $\fan$ is the right divisorial fan: 
it computes the Cox ring in terms of $\fan$, 
and, indeed, we obtain again $\mathcal{R}(X)$. 
Now we subdivide and obtain a divisorial fan 
having the refined slices as indicated 
in the following picture.

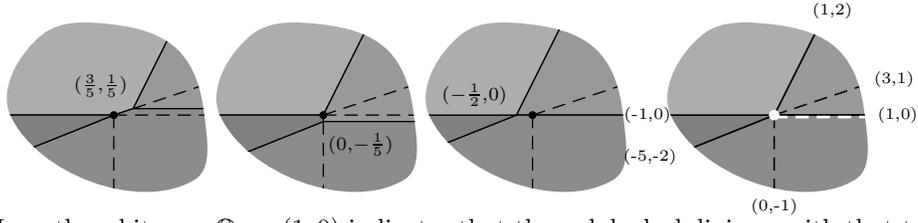
\threefoldD

\noindent
Here, the white ray $\QQ_{\geq 0}\cdot (1,0)$ indicates 
that the polyhedral divisors with that tail have affine loci. 
According to~\cite[Cor.~4.9.]{HaSu}, the corresponding Cox ring 
is given by
\begin{eqnarray*} 
\mathcal{R}(\widetilde{X})
& = & 
\KK[T_1, \ldots, T_7]/\bangle{T_1T_2^5 + T_3T_4^5 + T_5^2T_6}.
\end{eqnarray*} 

We have to check that $\widetilde{X}$ is smooth. 
Let us do this explicitly for the affine chart 
defined by the polyhedral divisor $\D$ 
with tail cone $\D_Y = \cone((1,2),(3,1))$. 
Then $\D$ is given by 
\begin{eqnarray*}
\D 
& = & 
\left(\left(\frac{3}{5},\frac{1}{5}\right) + \sigma\right) \cdot \{0\} 
\ + \ 
\left(\left[-\frac{1}{2},0\right] \times 0 + \sigma\right)\cdot \{\infty\}.
\end{eqnarray*}
Thus, $\cone(\{1\} \times \D_0) + \cone(\{-1\} \times \D_\infty)$
is generated by $(5,3,1)$, $(-2,-1,0)$ and $(-1,0,0)$; 
in particular, it is a regular cone.
This implies smoothness of the affine chart $X(\D)$.
Furthermore, we look at the affine charts 
defined by the polyhedral divisors 
$\D$ with tail cone 
$\D_Y = \cone(1,0)$. 
Since they have affine locus, 
we have to check 
$\cone(\{1\} \times \D_y)$, 
where $y \in Y$.
For $y \neq 0, 1$,
we have $\D_y = \D_Y$.
In this case, 
$\cone(\{1\} \times \D_y)$ is 
generated by $(1,1,0)$, $(0,1,0)$
and thus is regular.
For $y=0$, we obtain
that $\cone(\{1\} \times \D_y)$
is generated by $(5,3,1)$, $(1,0,0)$, $(0,1,0)$ 
and this is regular.
For $y=1$ we get the same result. 
Hence, the polyhedral divisors 
with tail cone $\D_y = \cone(1,0)$ 
give rise to smooth affine charts. 

Now we compute the discrepancies according to
Proposition~\ref{prop:discrepancies}.
The resolution has two exceptional divisors 
$D_{\infty, \mathbf{0}}$ and $E_{(1,0)}$. 
We work in the chart defined by 
the divisor $\D \in \fan$ with tail cone 
$\D_Y = \cone((1,2),(1,0))$.
The resulting system of linear equations 
and its unique solution are given by 
\[
\left(\begin{array}{ccccc|c}
  -1 & -1 & -1 & 0 & 0    & -1\\
  5 & 0 & 0 & 3 & 1    & 1\\
  0 & 1 & 0 & 0 & 0    & 1\\
  0 & 5 & 0 & 0 & -1   & 1\\
  0 & 0 & 2 & -1 & 0   & 1
\end{array}\right), 
\qquad\qquad
\begin{pmatrix}
  \alpha_0\\
  \alpha_1\\
  \alpha_\infty\\
  \hline
  u
\end{pmatrix}
\ = \
\begin{pmatrix}
  0\\
  1\\
  0 \\
  \hline
  -1\\
  4
\end{pmatrix}.
\]
The formula for the discrepancies yields 
$d_{\infty,\mathbf{0}}= -1$ and $d_{(1,0)}= -2$. 
In particular, $X$ has non-canonical singularities. 
By a criterion from~\cite[Sec.~3.4.]{tidiv},
we know that $D_{\infty, \mathbf{0}} + 2 \cdot E_{(1,0)}$ 
is a nef divisor. 
It follows that $\t{X}$ is a minimal model over $X$.
\end{proof}


\begin{thebibliography}{}%
%
\bibitem{MR2207875}
K. Altmann, J. Hausen:
Polyhedral divisors and algebraic torus actions.
Math. Ann. 334 (2006), no. 3, 557--607.
%
\bibitem{divfans} 
K. Altmann, J. Hausen, H. S\"u{\ss}:
Gluing Affine Torus Actions Via Divisorial Fans.
Transformation Groups 13 (2008), no. 2, 215--242.
%
\bibitem{Bat1}
V.V.~Batyrev:
Toric Fano threefolds. 
Izv. Akad. Nauk SSSR Ser. Mat.  45  (1981), 
no. 4, 704--717, 927. 
%
\bibitem{Bat2}
V.V.~Batyrev: 
On the classification of toric Fano 4-folds. 
J. Math. Sci. New York 94, 1021--1050 (1999).
%
\bibitem{BeHa1} 
F.~Berchtold, J.~Hausen: 
Homogeneous coordinates for algebraic varieties.  
J. Algebra  266  (2003),  no. 2, 636--670.
%
\bibitem{CP}
Cheltsov, Ivan; Park, Jihun
Sextic double solids.
In. Bogomolov, Fedor (ed.) et al.: 
Cohomological and geometric approaches 
to rationality problems. 
New Perspectives. Boston, MA: Birkh\"{a}user. 
Progress in Mathematics 282, 75--132 (2010).
%
\bibitem{CCC}
J.J.~Chen, J.A.~Chen, M.~Chen:
On quasismooth weighted complete 
intersections.
Preprint, arXiv:0908.1439.
%
\bibitem{CG}
C.H.~Clemens, P.A.~Griffiths:
The intermediate Jacobian of the cubic threefold. 
Ann. Math. (2) 95 (1972), 281--356.
%
\bibitem{Co}
A.~Corti:
Singularities of linear systems and 
3-fold birational geometry.
 Explicit birational geometry of 3-folds,  
259–-312, 
London Math. Soc. Lecture Note Ser., 281, 
Cambridge Univ. Press, Cambridge, 2000.
%
\bibitem{CM}
A.~Corti, M.~Mella:
Birational geometry of terminal quartic 3-folds I. 
Am. J. Math. 126 (2004), No. 4, 739--761.
%
\bibitem{CPR}
A.~Corti, A.~Pukhlikov, M.~Reid:
Fano 3-fold hypersurfaces. 
Explicit birational geometry of 3-folds,  
175--258, 
London Math. Soc. Lecture Note Ser., 281, 
Cambridge Univ. Press, Cambridge, 2000.
%
\bibitem{Gri}
M.M.~Grinenko:
Mori structures on a Fano threefold of index 2 and degree 1. 
Tr. Mat. Inst. Steklova  246  (2004),  
Algebr. Geom. Metody, Svyazi i Prilozh., 116--141;  
translation in  Proc. Steklov Inst. Math.  2004,  
no. 3 (246), 103--128.
%
\bibitem{Gri2}
M.M.~Grinenko:
Birational automorphisms of a three-dimensional double cone. 
Mat. Sb.  189  (1998),  no. 7, 37--52;  
translation in  Sb. Math.  189  (1998),  no. 7-8, 991--1007. 
%
\bibitem{Ha2}
J.~Hausen:
Cox rings and combinatorics II.
Mosc. Math. J., 8 (2008),  711--757.
%
\bibitem{HaSu}
J.~Hausen, H.~S\"u{\ss}:
The Cox ring of an algebraic variety with torus action.
Adv. Math.  225  (2010),  no. 2, 977–-1012.
%
\bibitem{IaFl}
A.R.~Iano-Fletcher:
Working with weighted complete intersections.  
Explicit birational geometry of 3-folds,  101–-173, 
London Math. Soc. Lecture Note Ser., 281, 
Cambridge Univ. Press, Cambridge, 2000.  
%
\bibitem{Is}
M.-N.~Ishida:
Graded factorial rings of dimension $3$ of a restricted type.  
J. Math. Kyoto Univ.  17  (1977), no. 3, 441--456.
%
\bibitem{fano3}
V.A.~Iskovskih:
Fano threefolds. II. 
Izv. Akad. Nauk SSSR Ser. Mat. 42 (1978), no. 3, 506-–549.
%
\bibitem{Isk80}
V.A.~Iskovskih:
Birational automorphisms of three-dimensional algebraic varieties.
Current problems in mathematics, Vol. 12 (Russian),  
pp. 159--236, 239 (loose errata), VINITI, Moscow, 1979. 
%
\bibitem{IM}
V.A.~Iskovskih, Yu.I.~Manin:
Three-dimensional quartics and counterexamples 
to the L\"uroth problem. 
Mat. Sb. (N.S.) 86 (128)  (1971), 140--166. 
%
\bibitem{JoKo}
J.M.~Johnson, J.~Koll\'{a}r:
Fano hypersurfaces in weighted projective 4-spaces.  
Experiment. Math.  10  (2001),  no. 1, 151--158. 
%
\bibitem{KKV}
F.~Knop, H~Kraft, T.~Vust: 
The Picard group of a $G$-variety.
In: Algebraic Transformation Groups 
and Invariant Theory, 
DMV Seminar, Band~ 3, Birkh\"auser.
%
\bibitem{Mo}
S.~Mori:
Graded factorial domains.
Japan J. Math. 3 (1977),
no. 2, 223--238.
%
\bibitem{tidiv} 
L. Petersen, H. S\"u{\ss}:
Torus invariant divisors.
Preprint, arXiv:0811.0517 (2008),
to appear in Israel J. Math.
%
\bibitem{Pu}
A.V.~Pukhlikov: 
Birational automorphisms of a three-dimensional quartic 
with a simple singularity. 
Mat. Sb. (N.S.)  135(177)  (1988),  no. 4, 472--496, 559;  
translation in  Math. USSR-Sb.  63  (1989),  no. 2, 457--482.
%
\bibitem{Pu2}
A.V.~Pukhlikov:
Birational isomorphisms of four-dimensional quintics. 
Invent. Math. 87, no.~2 (1987), 303--329.
%
\bibitem{WaWa} 
K.~Watanabe, M.~Watanabe: 
The classification of Fano 3-folds with torus embeddings. 
Tokyo J. Math. 5, 37--48 (1982).
%
\bibitem{tfano} 
H. S\"u{\ss}:
Canonical divisors on $T$-varieties.
Preprint, arXiv:0811.0626v1. 
%
\bibitem{tim}
A.S.~Tikhomirov:
The intermediate Jacobian of double ${\bf P}^{3}$ 
that is branched in a quartic. 
Izv. Akad. Nauk SSSR Ser. Mat.  44  (1980), no. 6, 
1329--1377, 1439. 
%
\bibitem{voi}
Voisin, Claire:
Sur la jacobienne interm\'{e}diaire du double solide d'indice deux. 
Duke Math. J.  57  (1988),  no. 2, 629--646. 
\end{thebibliography}
\end{document}